
\input epsf


 \hsize = 5.0in
 \hoffset = 0.75in
 \parskip = 0pt
 \topskip=18pt


\def \Buzbee {1}
\def \Chan {2}
\def \Dongarra {3}
\def \Dryja {4}
\def \Duff {5}
\def \George {6}
\def \Golub {7}
\def \Kee {8}
\def \Kellera {9}
\def \Kellerb {10}
\def \Householder {11}
\def \Parter {12}


 \newbox\partpage 
 \newdimen\mydimen 
 \newdimen\bighsize \bighsize = \hsize
 \newdimen\bigvsize \bigvsize = \vsize
 \newdimen\smallhsize \smallhsize = \hsize 
 \advance\smallhsize by-0.25in \divide\smallhsize by2 

\def \makeheadline {\vbox to0pt{\vskip -22.5pt\hbox
to\bighsize {\vbox to8.5pt{}\the \headline}\vss }\nointerlineskip}

\def \makefootline {\baselineskip = 24pt\hbox to\bighsize {\the
\footline }}

\def \begindoublecolumns{\begingroup
 \topskip = 12pt
 \def \pagebody{
  \splittopskip = \topskip \splitmaxdepth = \maxdepth
  \setbox0 = \vsplit255 to\mydimen 
  \setbox2 = \vsplit255 to\mydimen 
  \vbox {\unvbox\partpage\hbox to\bighsize{\box0\hfil\box2}}
  \global \vsize = 2\bigvsize
  \global \advance \vsize by-2\topskip}
 \output = {\global\setbox\partpage = \vbox{\unvbox255}}\eject
 \output = {\plainoutput}\tolerance = 1000 \hsize = \smallhsize 
 \mydimen = \bigvsize \advance\mydimen by-\ht\partpage
 \advance\mydimen by-\topskip
 \global\vsize = 2\mydimen}

\def \enddoublecolumns{
 \output = {
  \setbox0 = \vbox{\unvbox255}
  \mydimen = \ht0 \advance \mydimen by\topskip
  \advance \mydimen by-\baselineskip \divide \mydimen by2
  \splittopskip = \topskip
  {\vbadness = 10000 \loop \global \setbox3 = \copy0
  \global \setbox1 = \vsplit3 to\mydimen \ifdim \ht3>\mydimen
  \global \advance \mydimen by1pt \repeat}
  \setbox0 = \vbox to\mydimen {\unvbox1} 
  \setbox2 = \vbox to\mydimen {\unvbox3}
  \unvbox\partpage \hbox to\bighsize{\box0 \hfil \box2}}
 \eject \endgroup 
 \global \vsize = \bigvsize \pagegoal = \bigvsize}
 

 \font\bigbold=cmbx24
 \font\medbold=cmbx18
 \font\smallbold=cmbx12
 \font\helvetica=cmss10

 \def\bold{\bf}
 


 \def \one(#1) {\raise0.25ex\hbox{#1}} 
 \def \two(#1) {\lower0.25ex\hbox{#1}}  
 \headline={\sevenrm \ifnum \pageno=1{\hfil}\else {\ifodd
\pageno {\hfil MATRIX \one(S) \two(T) \one(R) \two(E) \one(T)
\two(C) \one(H) \two(I) \one(N) \two(G) \hfil}\else {\hfil JOSEPH F.
GRCAR\hfil}\fi}\fi}
 
 \def\thismonth{\ifcase \month \or JANUARY\or FEBRUARY\or
MARCH\or APRIL\or MAY\or JUNE\or JULY\or AUGUST\or
SEPTEMBER\or OCTOBER\or NOVEMBER\or DECEMBER\fi}

 \def\footremark{{\sevenrm \the \day \ \thismonth \ \the \year}}
 \def\pagenumber{{\tenrm \the \pageno}}
 \footline={\ifodd \pageno {\noindent \hfil \footremark \hfil
\pagenumber }\else {\noindent \pagenumber \hfil \footremark
\hfil}\fi}


 \outer\def\beginsection#1\par{\vskip0pt plus2.0in\penalty
-250\vskip0pt plus-2.0in\bigskip \bigskip \message {#1}\leftline
{\medbold #1}\nobreak \smallskip \noindent}

 \outer\def\beginsubsection#1\par{\vskip0pt plus2.0in\penalty
-250\vskip0pt plus-2.0in\bigskip \message{#1}\leftline {\smallbold
#1}\nobreak \smallskip \noindent}


 \outer \def \Algorithm #1. #2\par {\vskip \abovedisplayskip
{\noindent \narrower {\bf #1.} #2\par }\vskip \belowdisplayskip}


 \outer \def \figure #1. #2\par {{\parindent = 0pt \parskip=0pt
\ifdim \lastskip <\bigskipamount \removelastskip \penalty55
\bigskip \fi \bold Figure #1. \sl #2\par \bigskip}}

\outer \def \Picture #1. (#2 #3 #4) #5\par {{\parindent=0pt
\parskip=0pt 
\centerline {\epsfbox {#2}}
\smallskip
\bold Figure #1. \sl #5\par}} 

 \outer \def \Proclaim #1. #2\par {\bigbreak {\narrower \noindent
{\bold #1. }{\sl #2}\par} \ifdim \lastskip <\bigskipamount
\removelastskip \penalty55 \bigskip \fi}

 \outer \def \table #1. #2 #3\par {{\narrower \noindent {\bold
Table #1. }{\sl#2}\par \smallskip {\noindent \hfil{#3}}}}

 \outer\def\Table #1. #2 #3\par {\vfil \eject \ \vfil {\narrower
{\noindent \bold Table #1. }{\sl#2}\par} \bigskip {
\hfil#3\hfil} \vglue1.0in \vfil \eject}


 \def \Answer #1 {\raise 1pt\hbox{$\Leftarrow$}~}

 \def \Ddots{\lower0pt \vbox {\baselineskip=2pt
\lineskiplimit=0pt \kern2pt \hbox{.} \hbox{\kern0.4em.}
\hbox{\kern0.8em.}}}

 \def \filled (#1) (#2){\hbox to \hsize {#1\leaders \hbox to
2em{\hss.\hss}\hfill #2}}

 \def \listtwo #1 #2\endlist{#1\Ls#2}
 \def \listthree #1 #2 #3\endlist{#1\Ls#2\Ls#3}
 \def \listfour #1 #2 #3 #4\endlist{#1\Ls#2\Ls#3\Ls#4}
 \def \listfive #1 #2 #3 #4 #5\endlist{#1\Ls#2\Ls#3\Ls#4\Ls#5}
 \def \listsix #1 #2 #3 #4 #5 #6\endlist{#1\Ls#2\Ls#3\Ls#4\Ls#5\Ls#6}
 \def \listseven #1 #2 #3 #4 #5 #6 #7\endlist{#1\Ls#2\Ls#3\Ls#4\Ls#5\Ls#6\Ls#7}

 \def \Ls {,{\kern 0.5em}}

 \def \hquad {\hskip 0.5em\relax}

 \def \Question #1 {#1.~}

 \def \Vdots{\lower2pt \vbox {\baselineskip=4pt
\lineskiplimit=0pt \kern2pt \hbox{.} \hbox{.} \hbox{.}}}

 \def \R{{\rm \rlap I\kern 0.15emR}}
 
 \def \row{\mathop {\hbox {\rm row}}}


 \def \A{A^{SS \cdots S}}
 \def \AS{{\overline {A^S}}}

 \def \bars {{\rlap {$\kern 0.06em\overline {\phantom
{\jmath}}$}s}\vphantom s}
 \def \barx {{\rlap {$\kern 0.08em\overline {\phantom
{c}}$}x}\vphantom x}
 \def \barz {{\smash {\rlap {$\kern 0.09em\overline {\phantom
{\jmath}}$}z}\vphantom z}}
 \def \barA {{\smash {\rlap {$\kern 0.20em\overline {\phantom
{I}}$}A}\vphantom A}}
 \def \barI {{\smash {\rlap {$\kern 0.08em\overline {\phantom 
{I}}$}I}}\vphantom I}

 \def \Hbox {\hbox to 2.25in}

 \def \Matrix #1{\null \,\vcenter {\normalbaselines \ialign {\vrule
depth1.0ex height2.0ex width0em \hfil $##$\hfil&& \hfil $\;\;
##$\hfil\crcr \mathstrut \crcr \noalign {\kern-\baselineskip}
#1\crcr \mathstrut \crcr \noalign {\kern-\baselineskip}}}\,}

 \def\newi{{\tilde \imath}}
 \def\newj{{\tilde \jmath}} 

 \def\rightarrowfill{$\mathsurround=0pt \mathord- \mkern-6mu
 \cleaders\hbox{$\mkern-2mu \mathord- \mkern-2mu$}\hfill
 \mkern-6mu \mathord\rightarrow$}

 \def\S{{S \vrule depth0pt height1.5ex width0pt}}

 \def \vstrut {{\vrule depth3pt height12pt width0pt}}
 \def \Vstrut {{\vrule depth12pt height12pt width0pt}}

 \def\X{{^{-}}{\kern-.1667em}X}


 \headline={\hfil}


 \def\footremark{{\ }}
 \def\pagenumber{{\tenrm \the \pageno}}
 \footline={\ifodd \pageno {\noindent \hfil
\footremark \hfil \pagenumber}\else {\noindent \pagenumber
\hfil \footremark \hfil}\fi}

 



 \pageno = 3
 
 {\parindent=0pt \parskip=0pt 

 \line {\hphantom {UC??}\hfil SAND90-8723\hfil \hphantom {UC??}}
 \centerline {Unlimited Release}
 \centerline {Printed November 1990}

 \vglue1.25in
 \centerline {\bigbold Matrix \one(S) \two(T) \one(R) \two(E)
 \one(T) \two(C) \one(H) \two(I) \one(N) \two(G) }
 \vglue0.1in
 \centerline {\bigbold for Linear Equations*\footnote {}{\noindent
\tenrm *~Prepared for submission to SIAM Review.}}
 \vglue0.5in
 {\baselineskip=13pt \bold
 \centerline {\medbold Joseph F. Grcar} 
 \smallskip
 \centerline {Scientific Computing and Applied Math Division}
 \centerline {Sandia National Laboratories}
 \centerline {Livermore, CA 94551-0969 USA}}

 \vglue 1.0in
 \centerline {\medbold Abstract}
 \vglue 0.25in 
{\parindent=20pt Stretching is a new sparse matrix method that
makes matrices sparser by making them larger. Stretching has
implications for computational complexity theory and
applications in scientific and parallel computing. It changes
matrix sparsity patterns to render linear equations more easily
solved by parallel and sparse techniques. Some stretchings
increase matrix condition numbers only moderately, and thus
solve linear equations stably. For example, these stretchings
solve arrow equations with accuracy and expense preferable to
other solution methods.}

 \vfil \ \eject}

 {

 \ \vfil 
 \centerline {\medbold Acknowledgements}
 \vglue 0.33in

I thank Dr.~R.~J.~Kee and Prof.~M.~D.~Smooke for introducing me
to the need and practice of solving two-point boundary value
problems. Their use of the analytic transformations in Section~7
prompted this work. Thanks also to Dr.~C.~Ascraft for suggesting
Figure~5. I thank Dr.~J.~F.~Lathrop and Dr.~T.~H.~Jefferson for
reading the preliminary manuscript. And thanks to
Prof.~D.~E.~Knuth for {\TeX} and to Blue Sky Research for
Macintosh {\sl Textures}.

Some of this research was presented at the NASIG Conference at
Northstar Village in September 1987, at the SIAM Conference on
Linear Algebra in Madison in May 1988, and at the SIAM
Conference on Sparse Matrices at Salishan Lodge in May 1989.

 \vfil \vglue 2.0in \ \eject}

 {\parskip=0pt \parindent=0pt

 \ \vfil 
 \centerline {\medbold Contents}
 \vglue 0.5in
 \hfil \vbox {\hsize=4.5in \openup 2\jot
 \filled (1. Introduction) (7)
 \filled (2. Stretching Equations) (13)
 \filled (3. Simple Stretchings) (16)
 \filled (4. Numerical Stability) (21)
 \filled (\quad 4a. Condition Numbers) (22)
 \filled (\quad 4b. A Priori Accuracy) (24)
 \filled (5. Arrow Matrices) (26)
 \filled (6. Deflated Block Elimination) (29)
 \filled (7. Antecedents) (33)
 \filled (References) (37)
 \filled (Appendix 1. Proofs) (38)
 \filled (Appendix 2. Figure Explanations) (56)
 }
 \vglue 2.0in \vfil \eject}

 \ \vfil \eject 

\beginsection {1. Introduction}

Many matrices of computational interest contain mostly zeroes
and so are called {\sl sparse}. Yet sparse algorithms exploit both
the quantity of zeroes and the placement of nonzeroes, so in a
practical sense sparse matrices are those with a few nonzeroes in
the right place. {\sl Stretching\/} is a new sparse matrix method
that increases sparsity by rearranging the nonzeroes into larger
matrices. Some systems of linear equations can be solved more
easily by stretching them first. Stretching thereby addresses two
fundamental issues in scientific computing.

\Proclaim Question 1. What are the limits of easy parallelism?

Computations with uniform data dependencies lend themselves to
parallel execution, but small changes to regular dependency
structures inhibit parallelism. Figure~1 shows an uniform
structure with a disastrous perturbation. The irregularity may
represent a globally synchronized task or globally shared data.
Both are troublesome to parallel machines of various kinds.
Question~1 asks whether these irregularities necessarily block
parallelism.

\midinsert {\Picture 1. (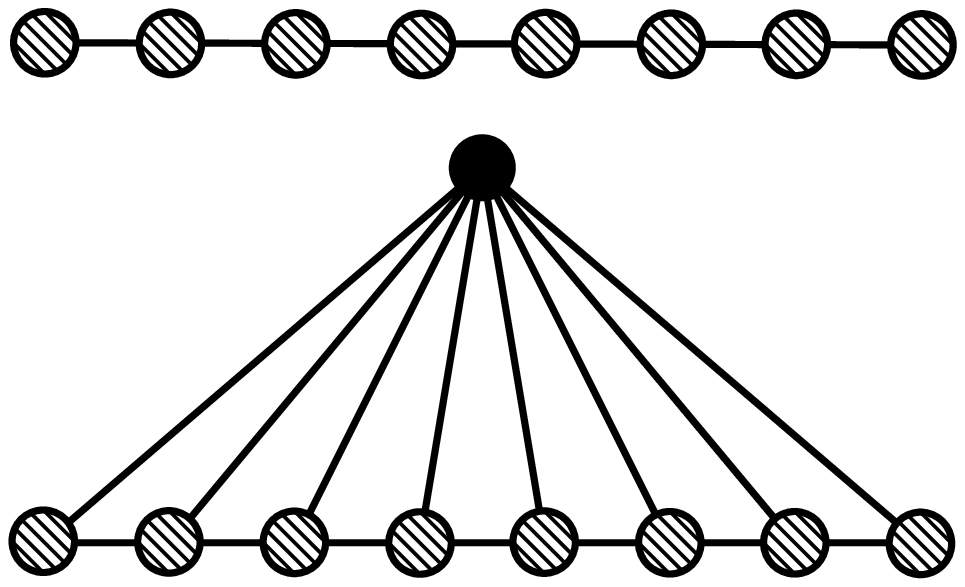 3.81 2.31) Dependency structures
affording easy and uneasy parallelism. \par }\endinsert

Figure~1's dependency structures occur in scientific computing as
{\sl occupancy graphs\/} for sparse matrices. This terminology is
new but the concept is well known. For a system of linear
equations written in matrix notation, $Ax = y$, the matrix
diagonal positions become the graph's vertices, and
if the row of one vertex has a nonzero entry in the column of
another, then an edge connects the two vertices.\footnote {$^1$}
{An edge connects vertices $j$ and $k$ when a nonzero occupies
matrix entry $(j,k)$. Parter [\Parter ] originated the study of
Gaussian elimination using these graphs [\Duff,~p.~4] but didn't
name them. The sparse matrix literature now prefers ``the graph
of the matrix'' or the {\sl matrix graph}, but still doesn't award the
concept a formal definition or a separate place in the index. See
also [\George ]. Conversely, the matrix whose entry $(j,k)$ is
nonzero when an edge connects vertices $j$ and $k$ is well known
in combinatorial mathematics as the {\sl adjacency matrix\/} of a
graph. Both concepts extend to directed graphs, and may include
loop edges $(j,j)$.} Figure~2 displays a matrix whose occupancy
graph is the distorted one of Figure~1. The dense column
represents a variable that appears in every equation, the dense
row represents an equation that includes all the variables, and
both are common in problems from linear programming and
differential equations. Dense rows and columns entail slow global
communication on computers with massive parallelism and
distributed memories. They complicate load balancing on
computers with limited parallelism and shared memories. 

\topinsert {
 {\Picture 2. (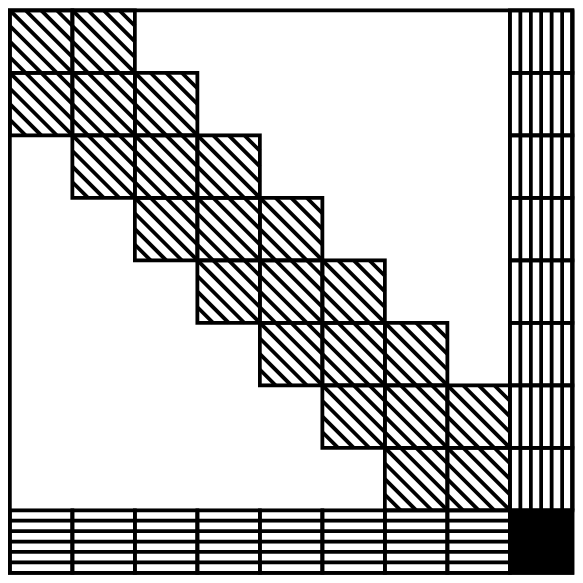 2.29 2.29) Matrix whose occupancy graph is
the irregular one in Figure~1. This matrix is visually sparse but
functionally dense. \par}
 \bigskip 
 {\Picture 3. (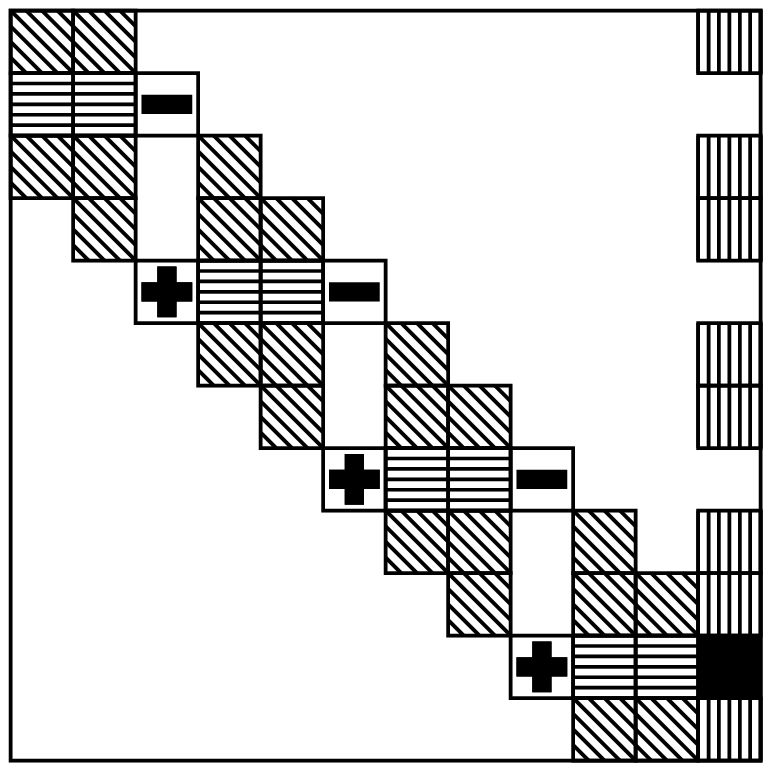 3.04 3.04) Sparser form of the matrix in
Figure~2 obtained by row stretching. The stretched rows sum to
the original dense row. \par}
 }\endinsert 

 \topinsert {\vfil
 {\Picture 4. (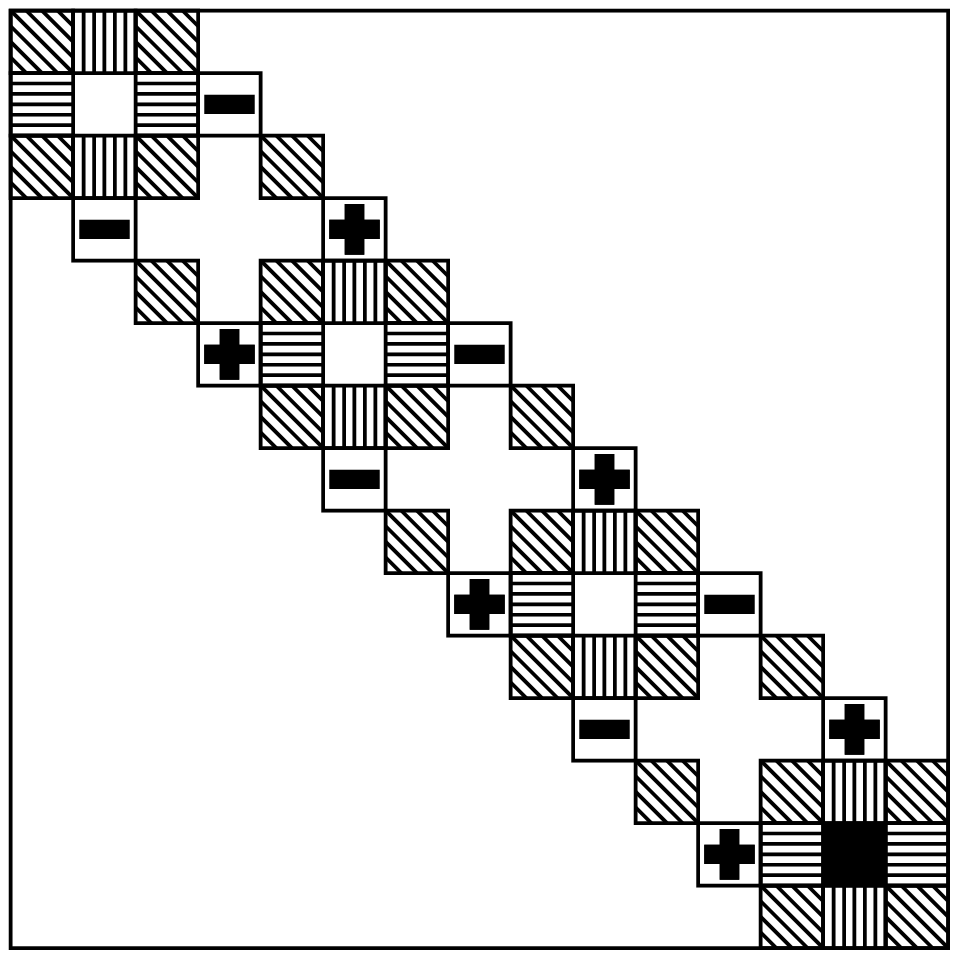 3.79 3.79) Sparsest form of the matrix in
Figure~2 obtained by row and column stretching. The stretched
rows and columns sum to the original dense row and column,
respectively. \par }
 \vfil
 {\Picture 5. (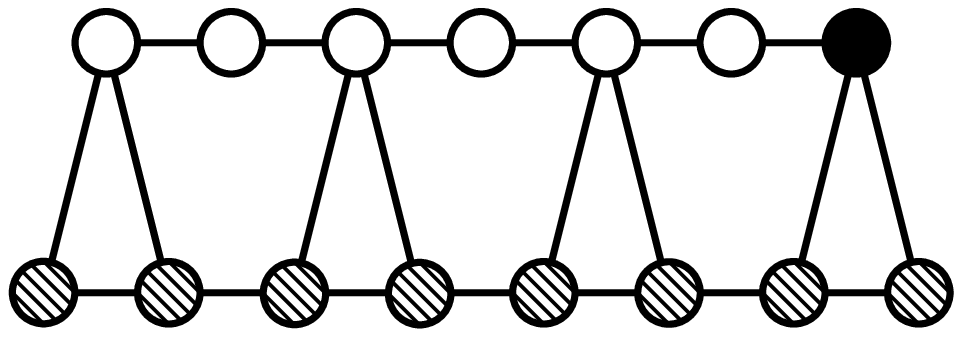 3.81 1.31) Occupancy graph for the matrix in
Figure~4. Stretching makes the matrix in Figure~2 sparse, and
makes the irregular dependency graph in Figure~1 uniform. \par}
 }\endinsert 

Stretching removes dense rows and columns that frustrate
parallel processing. Figures~3 and~4 exhibit stretched versions
of Figure~2's matrix, and Figure~5 shows the altered occupancy
graph. These particular stretchings move entries of dense rows
and columns into new, sparser rows and columns. They glue the
scattered pieces together by introducing some new nonzeroes.
Compared to the original matrices, the stretched matrices are
larger and have the same nonzeroes in different places. Whence
the name {\sl stretching}. 

\Proclaim Question 2. What is the price of accuracy?

Computational complexity theory usually treats a single
algorithm and so overlooks a central concern in scientific
computing. More complex algorithms may be needed to maintain
accuracy when a problem's data changes.\footnote{$^2$}{The
serial time complexity of a calculation is the number of
operations it performs, the space complexity is the number of
memory cells it touches. For systems of linear equations solved
by matrix factorization, space complexity is roughly the nonzero
population of the factors.}

The complexity of finding accurate solutions can be a strongly
discontinuous function of the problem. This is illustrated by linear
equations with the irregular dependencies of Figures~1 and~2
whose coefficient matrices vary with a parameter. Figure~6
shows the matrices are well-conditioned so it is feasible to ask
for accurate solutions. Figures~7 and~8 show the accuracy and
complexity vary greatly. Some parameter values demand much
more complex solution algorithms. 

\midinsert {\Picture 6. (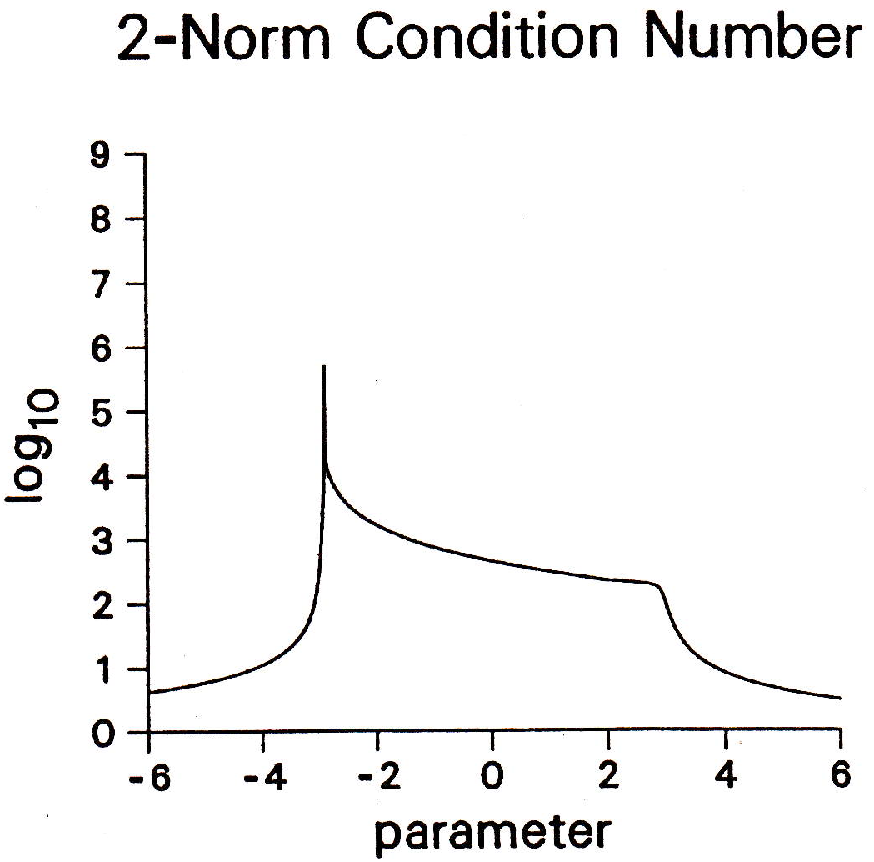 xxx 3.5) $2$-norm condition numbers for
parameterized matrices of order $51$ with sparsity patterns like
the matrix in Figure~2. Appendix~2 and Section~1 explain the
calculations. \par \medskip \bigskip}\endinsert

 \pageinsert {\vfil
 {\Picture 7. (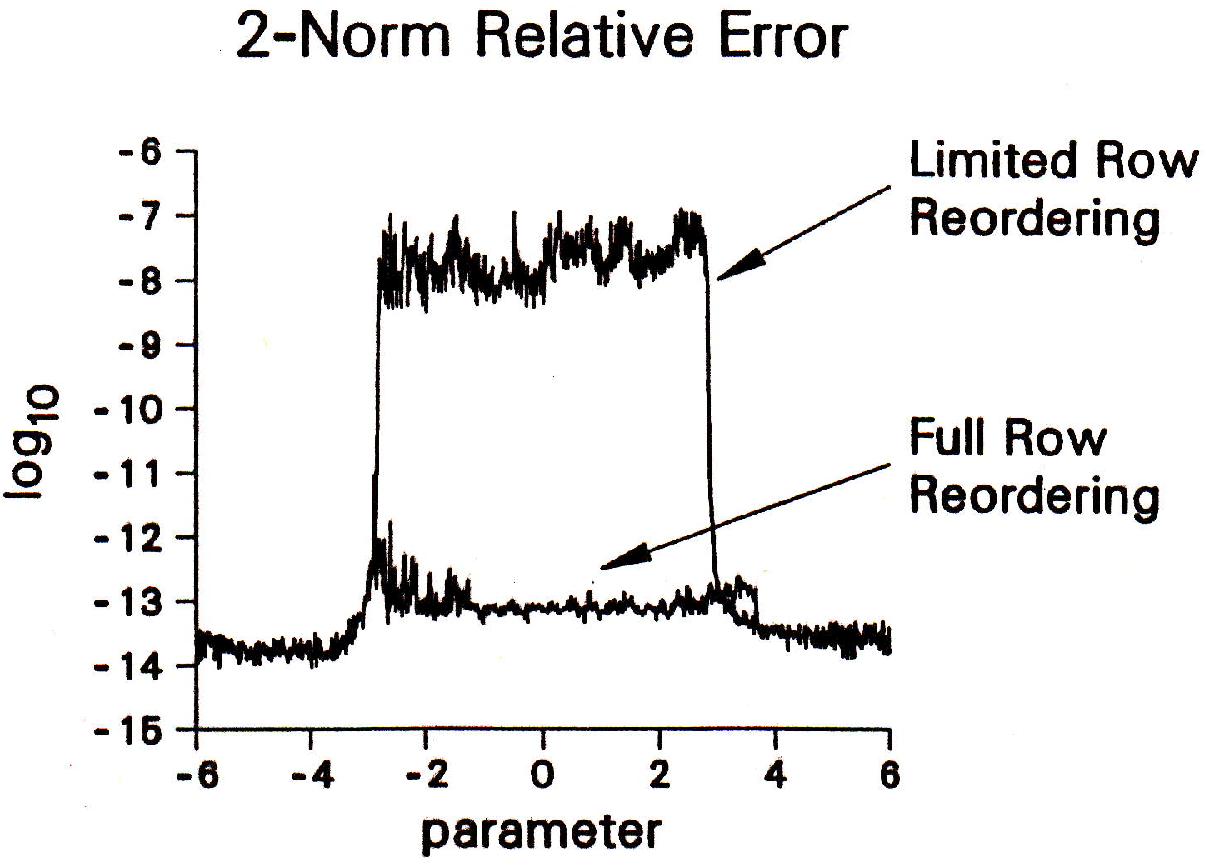 xxx 3.5) Maximum $2$-norm relative errors for equations
$A x = y$, with $20$ different $y$'s and the parameterized
matrices $A$ of Figure~6, solved by triangular factorization. The
lower curve allows full row reordering. The upper curve restricts
row reordering to the tridiagonal band.  Appendix~2 and
Section~1 explain the calculations. \par}
 \vfil
 {\Picture 8. (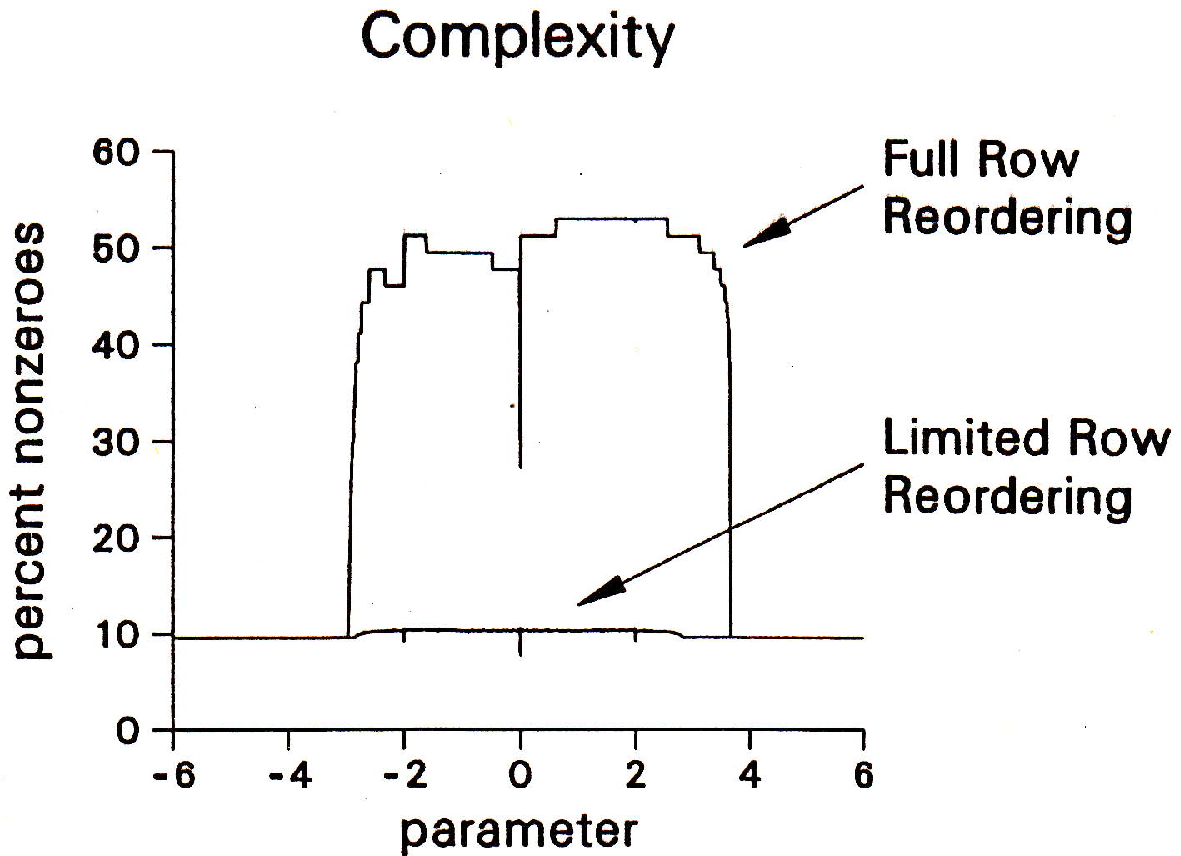 xxx 3.5) Percent of non-zeroes in the triangular factors of
the matrices of Figure~6. The upper curve allows full row
reordering. The lower curve restricts row reordering to the
tridiagonal band. Appendix~2 and Section~1 explain the
calculations. \par }
 \vfil}\endinsert
 
The increased complexity stems from the reordering algorithms
that stabilize matrix computations. The complexity in Figure~8
jumps when reordering is needed to maintain uniformly low errors
in Figure~7, as follows. If reordering selects a dense row to
participate at an early stage of the factorization, it engenders
more of the same, and increases the likelihood that additional
dense rows will be selected, and created. So many zeroes may be
lost in this way that the factors become completely dense and
the complexity becomes very high.

Stretching removes dense rows and columns that make
reordering expen\-sive. Stretched matrices have only sparse
rows and columns and therefore have fewer or no reorderings that
entail many nonzeroes. Although stretched matrices are larger,
they are likely factored more easily. Figures~9 and~10 display
the accuracy and complexity when the matrices of Figure~6
stretch in the manner of Figure~3. The accuracy matches
Figure~7's best; the complexity almost matches Figure~8's
lowest. Stretching achieves high accuracy and low complexity.

 \pageinsert {\vfil
 {\Picture 9. (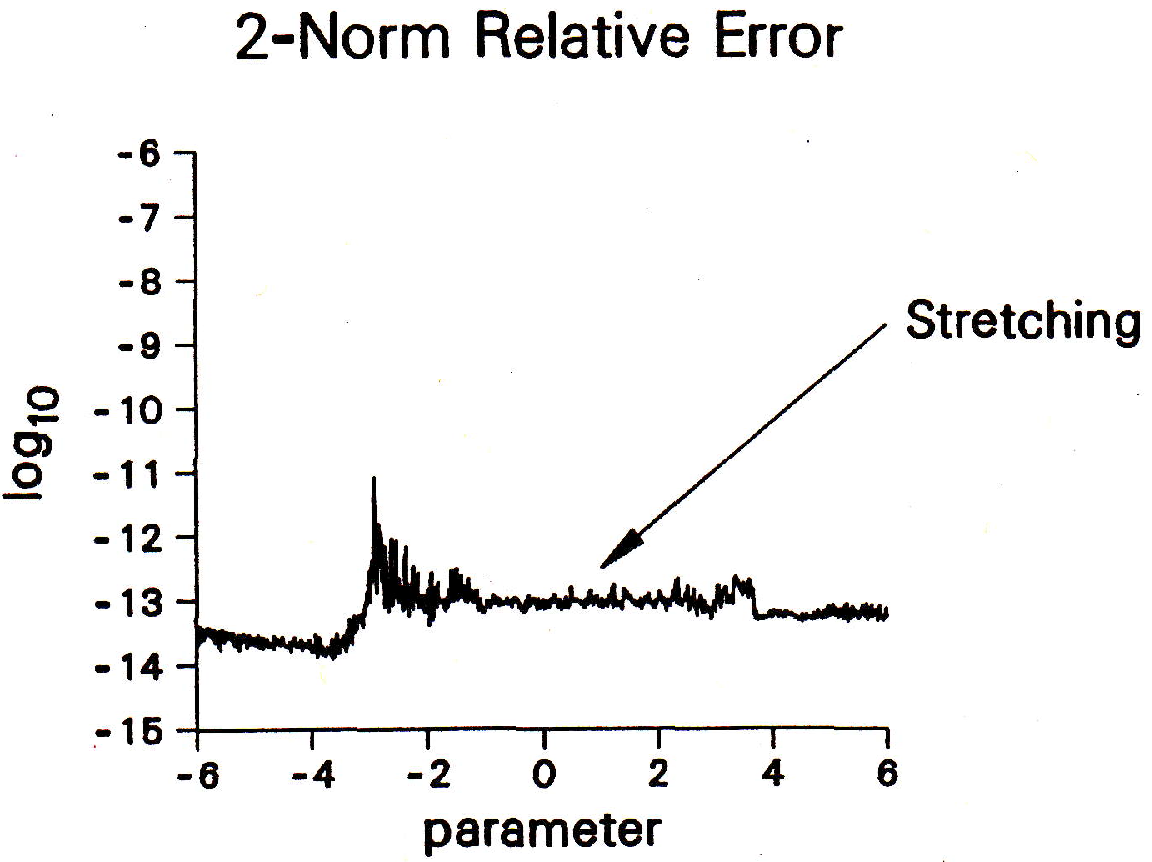 xxxx 3.5) $2$-norm relative errors for the equations of
Figure~7 solved by triangular factorization with full row
reordering after stretching in the manner of Figure~3.
Appendix~2 and Section~1 explain the calculations. \par}
 \vfil
 {\Picture 10. (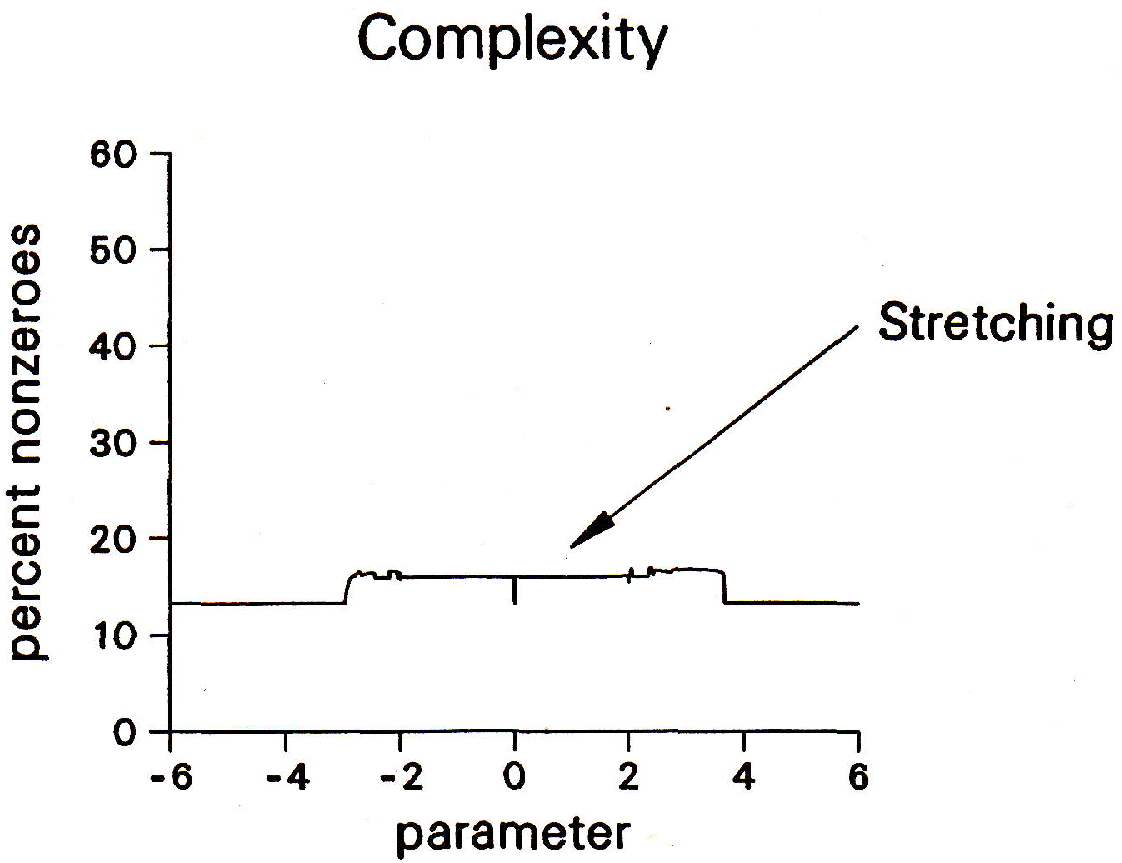 xxx 3.5) Percent of non-zeroes in the triangular factors of
the stretched matrices of Figure~9. The percentages are relative
to the size of the unstretched matrices. Appendix~2 and
Section~1 explain the calculations. \par}
 \vfil}\endinsert

Many scientific calculations implicitly avoid matrices with
inconvenient sparsity patterns. The final section of this paper
describes a precedent that inspires matrix stretching: analytic
transformations that ease numerical solution of some differential
equations. This paper is the first step toward making stretching a
purely algebraic---and therefore a broadly applicable---tool of
scientific computation. 

These are the paper's major results. First, stretching is
recognized as a sparse matrix method with implications beyond
numerical linear algebra and with potentially widespread
applications. It does not appear in the sparse matrix literature,
but it has been used indirectly to prepare some differential
equations for numerical solution. Second, some stretchings are
shown to increase matrix condition numbers moderately. The
proof of this is different from others in linear algebra and may
have independent interest. Third, the a priori error bounds for
solving linear equations are proved to increase only slightly with
stretching. Fourth, stretching's reliability and economy are
demonstrated by the special class of arrow equations for which
stretching is found preferable to other solution methods.

This paper is organized as follows. Section~2 presents a general
framework for constructing stretchings that solve linear
equations. The stretchings that implicitly accompany some
differential equations follow naturally in Section~3, where they
are christened {\sl simple row and column stretchings}. Section~4
shows these stretchings stably solve linear equations when some
parameters are properly chosen. Application to arrow matrices is
made in Section~5, and comparison with deflated block
elimination is made in Section~6. Finally, Section~7 describes the
differential transformations that inspired this work.
Odd-numbered sections are specific and accessible; Section~2 is
more general and introduces notation used throughout the paper;
Sections~4 and~6 are more technical. Applications to parallel
processing and randomly sparse matrices are not developed
beyond the suggestions made in this Introduction. To improve
readability, appendices contain proofs of theorems and
descriptions of numerical experiments. 

\beginsection {2. Stretching Equations}

What is needed to begin is a {\sl stretching\/} process that
associates the matrix $A$ with a larger matrix $A^S$. $$A
\rightarrow A^S$$ The reason for stretching is something about
$A$ makes $A x = y$ difficult to solve and something about $A^S$
makes $A^S z = y^S$ easier. Stretchings and {\sl squeezings\/} are
needed for vectors too. The overall solution process then
consists of first stretching $A \rightarrow A^S$ and $y
\rightarrow y^S$, next solving $A^S z = y^S$, and finally squeezing
$z \rightarrow z_\S = x$. $$\matrix {A^S& z& =& y^S \cr
\uparrow& \downarrow&& \uparrow \cr A& z_\S& =& y \cr }$$ The
superscript ${}^S$ indicates something dimensionally bigger than
the matrix or vector underneath, the subscript ${}_\S$ indicates
something smaller. In this notation the process of solving linear
equations is simply $$x = ((A^S)^{-1} y^S)_\S \; .$$ Little is gained
by greater formalism. Of interest rather are stretchings and
squeezings that work. They can be anything at all provided the
result is something useful like $A x = y$. 

Matrix stretchings with the ancillary vector operations needed to
solve linear equations are difficult to find. The stretchings
illustrated in Figures~2, 3 and~4 fit a common pattern which is
of interest because it may aid the discovery of more. The
pattern springs from a sequence of assumptions which might be
altered to obtain different stretchings. The first assumption is
(1) the stretched matrix be square and nonsingular if the original
is. Alternate courses are possible, for example, stretching might
produce under- or over-determined equations to be solved by
least squares methods.

The next assumption is (2) the vector operations be linear and
more or less independent of $A$ and $A^S$. $$\eqalign {y
\rightarrow y^S& := Y^{-}y\cr z \rightarrow z_\S& := \X z\cr }$$
Alternate courses might employ affine transformations. If the
stretchings and squeezings solve $A x = y$ for all $y$, then
linearity implies $$x = ((A^S)^{-1} y^S)_\S = \X (A^S)^{-1} Y^{-} y$$
and makes the search for stretchings the search for {\sl
oversize\/} factorizations. $$A^{-1} = \X (A^S)^{-1} Y^{-}$$ There
are many of these, but not many whose factors $\X$ and $Y^{-}$ are
independent of $A$ and $A^S$. Lacking some independence the
vector operations could degenerate to applying $A^{-1}$ and
stretching would gain nothing. 

An acceptable situation has the factors depending at most on the
sparsity patterns of $A$ and $A^S$. Factorizations with
restrictions of this kind are unlikely to be found even with
explicit knowledge of $A^{-1}$. Theorem~1 provides a mechanism
to overcome this difficulty by parameterizing the oversize
factorizations of $A^{-1}$.

\Proclaim Theorem~1. If $A$ and $A^S$ are nonsingular and
$$\vcenter {\offinterlineskip \halign {\Hbox {\hfil #\hfil }& \vstrut
\hfil #\hfil& \Hbox {\hfil #\hfil }\cr for some matrix $Y$& \Vstrut
or& for some matrix $X$\cr $\X := A^{-1} Y A^S$& \vrule& $\X :=
{}$ any left inverse of $X$\cr $Y^{-} := {}$ any right inverse of $Y$& 
\vrule& $Y^{-} := A^S X A^{-1}$\cr }}$$ then $A^{-1} = \X (A^S)^{-1}
Y^{-}$ (proof appears in Appendix~1).

The notational symmetry, $X$ and $Y$, $\X$ and $Y^{-}$, is
suggested by the Theorem's corollary.

\Proclaim Corollary to Theorem~1. If in addition $$\vcenter
{\offinterlineskip \halign {\Hbox {\hfil #\hfil }& \vstrut \hfil #\hfil 
& \Hbox {\hfil #\hfil }\cr $X := (A^S)^{-1} Y^{-} A$& \vrule& $Y := A
\, \X (A^S)^{-1}$\cr }}$$ then $\X X = I$, $Y Y^{-} = I$ and $A = Y A^S
X$ (proof appears in Appendix~1).

The third, more restrictive assumption is (3) $\X$ and $Y^{-}$ be
built from one of the Theorem's two sets of formulas. Alternate
courses might seek different expressions for $\X$ and $Y^{-}$, but
the formulas in Theorem~1 allow considerable freedom. A
likely stretching $A \rightarrow A^S$ might have several
matrices $Y$ and $X$ which yield factorizations for $A^{-1}$,
produced by the formulas above, that are appropriate for solving
linear equations. The sole criteria in choosing among them is the
convenience of applying the $\X$ and $Y^{-}$ actually used to
solve equations.

Something concrete begins to appear if the parametric matrix $Y$
or $X$ alone participates in the Corollary's factorization of $A$,
that is, if either $$\vcenter {\offinterlineskip \halign {\Hbox {\hfil 
#\hfil }& \vrule #& \Hbox {\hfil #\hfil }\cr $A^S := \left[ \matrix {B& 
G\cr } \right] P^t$& \vstrut& $A^S := P^t \left[ \matrix {B\cr G\cr }
\right]$\cr and $Y B = A$& \vstrut& and $B X = A$\cr }}$$ in which
$P$ is a permutation matrix. The extra columns and rows, both
denoted $G$ for {\sl glue}, can do more than make the stretched
matrices square. When they lie in the null spaces of $Y$ or $X$,
then Theorem~1's $\X$ or $Y^{-}$ depend only on $P$. $$\vcenter
{\offinterlineskip \halign {\Hbox {\hfil #\hfil }& \vrule #& \Hbox
{\hfil #\hfil }\cr $\X = A^{-1} Y A^S = \left[ \matrix {I& 0\cr } \right]
P^t$& \vstrut& $Y^{-} = A^S X A^{-1} = P^t \left[ \matrix {I\cr 0}
\right]$\cr provided $Y G = 0$& \vstrut& provided $G X =
0$\cr }}$$ The null space condition therefore makes $\X$ or $Y^{-}$
independent of $A^{-1}$, which is assumption (2). Moreover, if
$A^S$ is nonsingular then the extra columns or rows necessarily
are linearly independent, but with the null space condition
conversely, if the extra columns and rows are linearly
independent then $A^S$ is nonsingular, which is assumption (1).
This leads to the fourth and final assumption, which makes the
search for stretchings the search for one-sided factorizations of
$A$. It is embodied in the following Definition. The subsequent
Theorem~2 formalizes the preceding discussion and validates the
use of Definition~1's {\sl row and column stretchings\/} to solve
linear equations.

\Proclaim Definition 1, Row and Column Stretchings. A row or
column stretching $A \rightarrow A^S$ of square matrices has
$$\vcenter {\offinterlineskip \halign {\Hbox {\hfil #\hfil }& \vrule
#& \Hbox {\hfil #\hfil }\cr {\rm row stretching}& \Vstrut& {\rm
column stretching}\cr ${A^S := \left[ \matrix {B& G\cr } \right]
P^t}$& \vstrut& $A^S := P^t \left[ \matrix {B\cr G\cr } \right]$\cr 
with $Y B = A$& \vstrut& with $B X = A$\cr and $Y G = 0$& \vstrut
& and $G X = 0$\cr for some $Y$ of full rank& \vstrut& for some
$X$ of full rank\cr }}$$ in which $G$ has full rank and $P$ is a
permutation matrix, and chooses $$\vcenter {\offinterlineskip
\halign {\Hbox {\hfil #\hfil  }& \vrule #& \Hbox {\hfil #\hfil }\cr $\X
:= \left[ \matrix {I& 0\cr } \right] P^t$& \vstrut& $Y^{-} := P^t
\left[ \matrix {I\cr 0\cr } \right]$\cr $Y^{-} := {}$ any right inverse
of $Y$& \vstrut& $\X := {}$ any left inverse of $X$.\cr }}$$

\Proclaim Theorem~2. If $A \rightarrow A^S$ is a row or column
stretching and $A$ is nonsingular, then $A^S$ is nonsingular and
$A^{-1} = \X (A^S)^{-1} Y^{-}$ (proof appears in Appendix~1).

\Proclaim Corollary to Theorem~2. If $A \rightarrow A^S$ is a row
or column stretching of a nonsingular matrix $A$, if\/ $\X$ and $X$
are the auxiliary matrices in Definition~1, and if $A^S z = y^S$
are the stretched equations corresponding to $A x = y$, then not
only $\X z = x$ but also $z = Xx$ (proof appears in Appendix~1).

In summary, finding stretchings to solve equations involves two
tasks. One is to find better $A^S$, and assuming linear vector
operations, the other is to find workable $\X$ and $Y^{-}$. Row and
column stretchings are valuable because they are a rich class of
matrix stretchings for which acceptable $\X$ and $Y^{-}$ are
readily available.

Row stretchings can be viewed as being built in three stage. The
first, $$_nA_n \rightarrow {_{n+m}}B_n \quad \hbox{with} \quad
{_n}Y_{n+m}B_n = {_n}A_n,$$ increases the row dimension in a way
reversible by multiplication with some matrix $Y$. Whence the
name {\sl row stretching}. The new notation, $_{n+m}B_n$,
indicates a matrix of $n+m$ rows and $n$ columns. The second
stage, $$_{n+m}B_n \rightarrow \left[ \matrix {_{n+m}B_n& 
_{n+m}G_m\cr } \right] \quad \hbox{with} \quad {_n}Y_{n+m}G_n =
0,$$ adds new columns annihilated by $Y$. The third, $$\left[
\matrix {B& G\cr } \right] \rightarrow \left[ \matrix {B& G\cr }
\right] P^t = A^S,$$ scrambles the columns in a way reversible by
a permutation matrix $P$. 

Column stretchings are the transpose of row stretchings. Again
there are three stages. The first, $$_nA_n \rightarrow {_nB_{n+m}}
\quad \hbox{with} \quad {_n}B_{n+m}X_n = {_n}A_n,$$ increases
only the column dimension in a way reversible by multiplication
with some matrix $X$. Whence the name {\sl column stretching}.
The second stage, $$_nB_{n+m} \rightarrow \left[ \matrix
{_nB_{n+m}\cr _mG_{n+m}\cr } \right] \quad \hbox{with} \quad
{_m}G_{n+m}X_n = 0,$$ adds new rows annihilated by $X$. The
third, $$\left[ \matrix {B\cr G\cr } \right] \rightarrow P^t \left[
\matrix {B\cr G\cr } \right] = A^S,$$ scrambles the rows in a way
reversible by a permutation matrix $P$. 

\beginsection {3. Simple Stretchings}

The Introduction's stretchings receive a proper christening here.
Section~7 describes their prior use in the numerical solution of
ordinary differential equations, but now they are seen to be
legitimate offspring of general algebraic methods, and are named
{\sl simple row and column stretchings}. The following derivation
amounts to making specific choices for $B$ and $G$ in
Definition~1. 

A situation in which row stretching may be of use is that of a
single dense row which inhibits row reordering during triangular
factorization. This row represents a linear equation of the form
$$\displaylines {a_{j,1} x_1 + a_{j,2} x_2 + a_{j,3} x_3 + a_{j,4} x_4
+ a_{j,5} x_5 + a_{j,6} x_6 = y_j\cr \noalign {\smallskip}
\Updownarrow\cr \left[ \matrix {a_{j,1}& a_{j,2}& a_{j,3}& a_{j,4}&
a_{j,5}& a_{j,6}\cr } \right]\cr }$$ in which $a_{j,k}$, $x_k$ and $y_j$
are entries of $A$, $x$ and $y$ in $A x = y$. Row stretching might
be used to expand this row into something sparser. $$\left[
\matrix {a_{j,1}& a_{j,2}&&&& \cr &&  a_{j,3}& a_{j,4}&& \cr &&&&
a_{j,5}& a_{j,6}\cr } \right]$$ Section~4 shows this choice can
reduce the computational complexity of triangular factorization.
The first stage of row stretching, $A \rightarrow B$, simply
replaces the $j^{th}$ row of $A$ by the three stretched rows
above and optionally reorders the rows. This stage is undone by a
transformation $Y B = A$ that copies the untouched rows and
sums the three stretched ones. If row $j$ is the last in $A$ and if
the stretched rows replace it at the bottom, then $$Y = \left[
\matrix {1\cr & 1\cr && 1\cr &&& 1\cr &&&&  1\cr &&&&& 1& 1&
1\cr } \right].$$

The second stage, $B \rightarrow \left[ \matrix {B& G\cr } \right]$,
produces a square matrix by appending new columns which, to
make the stretched matrix nonsingular, must span the right null
space of $Y$. A column vector in this null space has zeroes in the
original rows of $A$ and sums to $0$ over the stretched rows.
After the second stage the stretched rows could be $$\left[
\matrix {a_{j,1}& a_{j,2}&&&&& - \sigma_1& \cr && a_{j,3}& 
a_{j,4}&&& + \sigma_1& - \sigma_2\cr &&&& a_{j,5}& a_{j,6}&& +
\sigma_2\cr } \right]$$ for some nonzero $\sigma_1$ and
$\sigma_2$. 

The third and final stage, $\left[ \matrix {B& G\cr } \right]
\rightarrow A^S$, reorders the columns. This is more than a
cosmetic detail because column order affects the complexity of
solving equations. The new columns could become the $3^{rd}$ and
$6^{th}$. Altogether $A^S$ has the following stretched rows.
$$\left[ \matrix {a_{j,1}& a_{j,2}& - \sigma_1&&&& \cr && +
\sigma_1& a_{j,3}& a_{j,4}& - \sigma_2&& \cr &&&&& +
\sigma_2& a_{j,5}& a_{j,6}\cr } \right]$$

$A^S$ can be used to solve $A x = y$ as follows. Step~$1$ forms
$y^S = Y^- y$ where $Y^-$ is any right inverse for $Y$. This
transformation copies entries of unstretched rows from $y$ to
$y^S$ and places numbers that sum to $y_j$ in the three stretched
rows. Step~$2$ solves $A^S z = y^S$. Step~$3$ forms $x = \X z$
where $\X = [ \matrix {I& 0\cr } ] P^t$ and $P$ is the permutation
matrix that reorders the columns in stage 3. This means the old
variables lie among the new in locations corresponding to the
original columns of $A$. The net result is the original equation has
been replaced by $$\vcenter {\openup1\jot \halign {& \hfil $#$\hfil 
\cr a_{j,1} x_1 + a_{j,2} x_2& {} - \sigma_1 s_1&&&& {} = t_1\cr & 
{} + \sigma_1 s_1& {} + a_{j,3} x_3 + a_{j,4} x_4& {}- \sigma_2
s_2&& {} = t_2\cr &&& {}+ \sigma_2 s_2& {}+ a_{j,5} x_5 + a_{j,6}
x_6& {} = t_3\cr }}$$ in which $s_1$ and $s_2$ are the new
variables and any numbers that sum to $y_j$ can appear on the
right. Different choices give different values to the new
variables, but of course the original variables remain unchanged.
\medskip

\begindoublecolumns Although column stretching is the
transpose of row stretching, significant conceptual differences
arise when solving equations. It is best to consider a separate
example---taking care to avoid the page costs of displaying
column vectors. A dense column represents a variable that occurs
in several linear equations of the form $$\vcenter {\openup1\jot
\halign {$\ldots + \hfil a_{#,k} x_k\hfil + \ldots = {}$& $y_#$\hfil 
\cr 1& 1\cr 2& 2\cr 3& 3\cr 4& 4\cr 5& 5\cr 6& 6\cr }} \hquad
\Longleftrightarrow \left[ \matrix {a_{1,k}\cr a_{2,k}\cr a_{3,k}\cr 
a_{4,k}\cr a_{5,k}\cr a_{6,k}} \right]$$ in which $a_{j,k}$, $x_k$ and
$y_j$ are entries of $A$, $x$ and $y$ in $A x = y$. Column
stretching can be used to expand this one column to something
sparser. $$\left[ \matrix {a_{1,k}&& \cr a_{2,k}&& \cr & a_{3,k}& \cr 
& a_{4,k}& \cr && a_{5,k}\cr && a_{6,k}\cr } \right]$$ The first
stage, $A \rightarrow B$, replaces the $k^{th}$ column of $A$ by
the three stretched columns above and optionally reorders the
columns. This stage is undone by a transformation $B X = A$ that
copies the untouched columns and sums the three stretched ones.
If column $k$ is the last in $A$ and if the stretched columns
replace it at the right side, then $$X = \left[ \matrix {1\cr & 1\cr 
&& 1\cr &&& 1\cr &&&& 1\cr &&&&& 1\cr &&&&& 1\cr &&&&& 
1\cr } \right].$$ The second stage, $$B \rightarrow \left[ \matrix
{B\cr G\cr } \right]$$produces a square matrix by appending new
rows that span the left null space of $X$. A row vector in this null
space has zeroes in the original columns of $A$ and sums to $0$
over the stretched columns. After the second stage the stretched
columns could be $$\left[ \matrix {a_{1,k}&& \cr a_{2,k}&& \cr & 
a_{3,k}& \cr & a_{4,k}& \cr && a_{5,k}\cr && a_{6,k}\cr - \sigma_1& 
+ \sigma_1& \cr & - \sigma_2& + \sigma_2\cr } \right]$$ for some
nonzero $\sigma_1$ and $\sigma_2$. The third stage reorders the
rows. If the new rows become the $3^{rd}$ and $6^{th}$, then $A^S$
has stretched columns $$\left[ \matrix {a_{1,k}&& \cr a_{2,k}&& \cr 
- \sigma_1& + \sigma_1& \cr & a_{3,k}& \cr & a_{4,k}& \cr & -
\sigma_2& + \sigma_2\cr && a_{5,k}\cr && a_{6,k}\cr } \right] .$$

Once again $A^S$ can be used to solve $A x = y$, but the steps
differ from the row case in several details. Step~$1$ forms $y^S = 
Y^- y$ $$Y^- = P^t \left[ \matrix {I\cr 0\cr } \right]$$ in which $P$ is
the permutation matrix that reorders the rows in stage $3$. This
transformation copies all entries of $y$ into $y^S$ and places
zeroes in the new rows. Step~$2$ solves $A^S z = y^S$. Step~$3$
forms $x = \X z$ where $\X$ can be any left inverse for $X$. Entries
of $z$ that correspond to original columns copy directly into $x$.
That is, unstretched columns retain their original variables.
Entries of $z$ that correspond to stretched columns coalesce in a
linear combination whose coefficients sum to $1$. That is, the
original variable $x_k$ equals any linear combination, with
coefficients summing to $1$, of the new variables for the
stretched columns. The net result is that the original equations
have been replaced by $$\vcenter {\openup1\jot \halign {#& $#$& 
\hfil $#$& $#$& \hfil $#$& $#$& ${} #$& ${} = #$\hfil \cr \hfil $\ldots +
a_{1,k}$& s_1&&&&& + \ldots& y_1\cr \hfil $\ldots + a_{2,k}$& 
s_1&&&&& + \ldots& y_2\cr \hfil ${} - \sigma_1$& s_1& {} +
\sigma_1& s_2&&&& 0\cr $\ldots + {}$\hfil&& a_{3,k}& s_2&&& +
\ldots& y_3\cr $\ldots + {}$\hfil&& a_{4,k}& s_2&&& + \ldots& 
y_4\cr && {} - \sigma_2& s_2& {} + \sigma_2& s_3&& 0\cr $\ldots
+ {}$\hfil&&&& a_{5,k}& s_3& + \ldots& y_5\cr $\ldots +
{}$\hfil&&&& a_{6,k}& s_3& + \ldots& y_6\cr }}$$ in which $s_1$,
$s_2$, and $s_3$ are the new variables. The equations make the
new variables equal to $x_k$ in principal, but machine
computation makes them different in fact. Section~4 considers
the effect of numerical error. \enddoublecolumns \medskip

Differences between row and column stretching therefore occur
in solving linear equations. Row stretching allows some freedom
in choosing the right side of the stretched equations, but
completely specifies how to recover the solution of the original
equations. The reverse is true for column stretching. Column
stretching completely specifies the right side, but allows some
freedom in recovering the solution. The simple stretchings
described above apply as well to blocks of rows and columns.

\Proclaim Definition~2, Simple Row and Column Stretchings. For a
system of linear equations $A x = y$ whose coefficient matrix
$$A = \left[ \matrix{A_1\cr A_2\cr } \right]$$ has a block of dense
rows $A_2$, simple row stretching partitions the columns into
$m$ blocks $$\left[ \matrix {A_{1,1}& A_{1,2}& A_{1,3}& \ldots& 
A_{1,m}\cr A_{2,1}& A_{2,2}& A_{2,3}& \ldots& A_{2,m}\cr } \right]
\left[ \matrix {x_1\cr x_2\cr x_3\cr \Vdots\cr x_m\cr } \right] = 
\left[ \matrix {y_1\cr y_2\cr } \right]$$ and replaces the equations
by $$\left[ \matrix {A_{1,1}& A_{1,2}& A_{1,3}& \ldots& A_{1,m}& 
0& 0& \ldots& 0\cr A_{2,1}&&&&& -D_1\cr & A_{2,2}&&&& +D_1& 
-D_2\cr && A_{2,3}&&&& +D_2& \Ddots\cr &&& \Ddots&&&& 
\Ddots& -D_{m-1}\cr &&&& A_{2,m}&&&& +D_{m-1}\cr } \right]
\left[ \matrix {x_1\cr x_2\cr x_3\cr \Vdots\cr x_m\cr s_1\cr 
s_2\cr \Vdots\cr s_{m-1}\cr } \right] = \left[ \matrix {y_1\cr 
t_1\cr t_2\cr t_3\cr \Vdots\cr t_m\cr } \right]$$ in which $D_1$,
$D_2$, $\ldots$, $D_{m-1}$ are nonsingular (presumably diagonal)
matrices and $y_2 = t_1 + t_2 + \ldots + t_m$. Alternatively, for
a system of linear equations $A x = y$ whose coefficient matrix
$$A = \left[ \matrix {A_1& A_2\cr } \right]$$ has a block of dense
columns $A_2$, simple column stretching partitions the rows into
$m$ blocks $$\left[ \matrix {A_{1,1}& A_{1,2}\cr A_{2,1}& 
A_{2,2}\cr A_{3,1}& A_{3,2}\cr \vdots& \vdots\cr A_{m,1}& 
A_{m,2}\cr } \right] \left[ \matrix {x_1\cr x_2\cr } \right] = \left[
\matrix {y_1\cr y_2\cr y_3\cr \vdots\cr y_m\cr } \right]$$ and
replaces the equations by $$\left[ \matrix {A_{1,1}& A_{1,2}\cr 
A_{2,1}&& A_{2,2}\cr A_{3,1}&&& A_{3,2}\cr \vdots&&&& 
\Ddots\cr A_{m,1}&&&&& A_{m,2}\cr 0& -D_1& +D_1\cr 0&& 
-D_2& +D_2\cr \vdots&&& \Ddots& \Ddots\cr 0&&&& -D_{m-1}& 
+D_{m-1}\cr } \right] \left[ \matrix {x_1\cr s_1\cr s_2\cr s_3\cr 
\vdots\cr s_m\cr } \right] = \left[ \matrix {y_1\cr y_2\cr y_3\cr 
\vdots\cr y_m\cr 0\cr 0\cr \vdots\cr 0\cr } \right]$$ in which
$D_1$, $D_2$, $\ldots$, $D_{m-1}$ are nonsingular (presumably
diagonal) matrices and $x_2 = s_1 = s_2 = \ldots = s_m$. The
matrices can be reordered both before and after the stretchings
and may assume a final appearance quite different from the
templates above. 

\Proclaim Theorem~3. A simple row stretching in the sense of
Definition~2 is a row stretching in the sense of Definition~1, and
similarly for column stretchings (proof appears in Appendix~1).

The general row and column stretchings of Definition~1 are
parameterized by auxiliary matrices $X$ and $Y$, respectively.
Ignoring reorderings, the simple row stretching of Definition~2
has $$Y = \left[ \matrix {I_1\cr & I_2& I_2& \ldots& I_2\cr }
\right]$$ in which $I_1$ and $I_2$ are identity matrices whose
orders match the row orders of $A_1$ and $A_2$. Simple column
stretching has $$X = \left[ \matrix {I_1\cr & I_2\cr & I_2\cr &
\vdots\cr & I_2\cr } \right]$$ where the orders of $I_1$ and $I_2$
match the column orders of $A_1$ and $A_2$. Theorem~2 can be
invoked with Theorem~3 and these $X$ and $Y$ to confirm the
nonsingularity of the stretched matrices and the procedure for
solving linear equations. Yet in this simple case these conclusions
can be obtained more directly. Theorem~4 shows the stretched
matrices are nonsingular.

\Proclaim Theorem~4. If $A \rightarrow A^S$ is a simple row or
column stretching as in Definition~2, then $$\det A^S = \det A
\prod_{j = 1}^{m-1} \det D_j$$ with perhaps a sign change when the
rows and columns are reordered as the definition allows (proof
appears in Appendix~1).

\beginsection {4. Numerical Stability}

Bounds on the rounding error for solving equations with and
without stretching compare favorably because prop\-er\-ly
formed stretchings increase matrix condition numbers at worst
moderately. This is the paper's major analytic result. 

Analyses of stretching's errors must consider more than matrix
condition numbers. The overall process for $A x = y$
first stretches $A \rightarrow A^S$ and $y \rightarrow y^S$,
then solves $A^S z = y^S$, and finally squeezes $z \rightarrow
z_\S = x$. $$\vcenter {\openup0.5\jot \halign {\hfil $\,#\,$\hfil&& 
\hfil $\,#\,$\hfil \cr A&&& y\cr \downarrow&&& \downarrow \cr 
A^S& z& {} = {}& y^S \cr & \downarrow \cr & z_\S\cr }}$$ The
manipulative steps introduce errors beyond those of solving $A^S
z = y^S$. Nevertheless, if the embedded solution process is
stable in the customary sense, and if $A$ is well conditioned, then
the overall process accurately solves $A x = y$.

The present analysis takes the standard approach toward
understanding finite precision computation. The errors are
interpreted as being governed by both the equations and the
algorithm. {\sl Error analyses\/} follow individual arithmetic
errors through an algorithm and are technically demanding, but in
stretching's case most errors arise in the solution of $A^S z = 
y^S$ and can be assumed accounted for by other analyses. These
accumulated errors are viewed as perturbing the equations rather
than the solution, and the solution's accuracy is assessed by two
inequalities. 

The {\sl vector perturbation inequality\/} emphasizes the role of
the equations. Any approximate solution $\barx$ for $A x = y$
exactly solves $A \barx = y - r$ in which $r$ is the {\sl residual\/}
$y - A \barx$. $$\vcenter {\openup1\jot {\halign {\hfil #\cr if $A x
=  y$ and\cr $\barx$ is any vector\cr }}} \quad \Longrightarrow
\quad {\| x - \barx \| \over \| x \|} \le \kappa (A) \, {\| r \| \over \| y
\|}$$ An algorithm might produce an $\barx$ with a small relative
residual (bars traditionally denote computed quantities), but no
matter how small, the {\sl condition number} $$\kappa (A) = \| A \|
\, \| A^{-1} \| \ge 1$$ scales the bound and perhaps the error. The
bound may not be sharp because the error varies with $A$ and $y$
(not merely linearly with the condition number as the bound
suggests). Yet the bound is valuable because the residual is
directly observable and because the condition number is intrinsic
to the matrix (and to the measurement of errors by norms---this
and the next inequality are valid for any consistent matrix-vector
norm). 

The {\sl matrix perturbation inequality\/} relies on details of the
solution method, with the approximate solution expected to be
the exact solution of an approximate problem derived from error
analysis. $$\vcenter {\openup1\jot {\halign {\hfil $#$\cr A x = y\cr 
(A + E) \barx = y\cr }}} \quad \Longrightarrow \quad {\| x - \barx \|
\over \| x \|} \le {\| A^{-1} E \| \over 1 - \| A^{- 1} E \|}$$ The solution
algorithm determines $\barx$ from $A$ and $y$, but the
perturbation $E$ may be any for which $E \barx = r$. The bound
additionally requires $\| A^{-1} E \| < 1$, implying $\| E \| < \| A \|$,
and a {\sl stable\/} algorithm has some $E$ provably small
relative to $A$, so the bound is usefully weakened as follows. $$ {\|
x - \barx \| \over \| x \|} \le {\| A^{-1} E \| \over 1 - \| A^{- 1} E \|}
\le \kappa (A) \, {\| E \| \over \| A \| - \| E \|}$$ Bounds upon some $\|
E \|$ that depend on $A$ but not on $y$ represent the error as
being independent of $y$ and have been derived for several
algorithms. They are primarily the work of J.~H.~Wilkinson and are
beyond the scope of this discussion. They can be found in many
texts including [\Duff ] [\Golub ] and references therein.

The inequalities above are too flexible to be of predictive value
when the error is small, but they diagnose the cause when the
error is large. They prove stable algorithms applied to
well-conditioned matrices yield accurate solutions. Section~6
illustrates the risk of calculating without a performance
guaranty. There, a plausible but imperfect method is found to
produce unexpectedly large errors. 

\beginsubsection {4a. Condition Numbers}

The condition numbers of stretched matrices vary with the newly
introduced nonzeroes which in some sense bind the stretched
matrices together. The glue lies in the submatrices $G$ and $D_1$,
$D_2$, $\ldots\,$, $D_{m-1}$ of Definitions~1 and~2. This section
finds glue that favorably bounds the condition of matrices
stretched by Definition~2.

The bounds are stated not for a single stretching, $A \rightarrow
A^S$, but rather for a sequence of stretchings. $$A \rightarrow
A^S \rightarrow A^{SS} \rightarrow \cdots \rightarrow A^{SS
\cdots S}$$ This generality anticipates sparse factorization
software that might stretch many times. For example, a row $$
\lineskip = \abovedisplayskip \tabskip=\centering \halign
to\displaywidth{\hfil $\displaystyle #$\tabskip=0pt&
$\displaystyle {}#$\hfil \tabskip=\centering \cr r &= \,\,
\left[\matrix {a_1& a_2& a_3& a_4& a_5& a_6\cr }\right]\cr
\noalign {\vskip \belowdisplayskip \hbox {could stretch once}}
r^{S} &= \left[\matrix {a_1& \cdot& \cdot& \cdot& \cdot& \cdot&
{-}\cr \cdot& a_2& a_3& a_4& \cdot& \cdot& {+}& {-}\cr \cdot&
\cdot& \cdot& \cdot& a_5& \cdot&& {+}& {-}\cr \cdot& \cdot&
\cdot& \cdot& \cdot& a_6&&& {+}\cr }\right]\cr \noalign {\hbox
{and then a descendent could stretch again.}} r^{SS} &=
\left[\matrix {a_1&&&&&& {-}\cr \cdot& a_2& \cdot& \cdot&
\cdot& \cdot& \cdot& \cdot& \cdot& {-}\cr \cdot& \cdot& a_3&
\cdot& \cdot& \cdot& \cdot& \cdot& \cdot& {+}& {-}\cr \cdot&
\cdot& \cdot& a_4& \cdot& \cdot& {+}& {-}& \cdot&& {+}\cr
&&&& a_5&&& {+}& {-}\cr &&&&& a_6&&& {+}\cr }\right]\cr
\noalign {\vbox {\noindent \normalbaselines The $\pm$'s are the
glue. Rows that contain mostly glue could stretch to rows that
contain only glue, with little apparent order.}} r^{SSS} &=
\left[\matrix {a_1&&&&&& {-}\cr & a_2&&&&&&&& {-}\cr \cdot&
\cdot& a_3& \cdot& \cdot& \cdot& \cdot& \cdot& \cdot& \cdot&
\cdot& {-}\cr \cdot& \cdot& \cdot& \cdot& \cdot& \cdot& \cdot&
\cdot& \cdot& {+}& \cdot& {+}& {-}\cr \cdot& \cdot& \cdot&
\cdot& \cdot& \cdot& \cdot& \cdot& \cdot& \cdot& {-}&& {+}&
{-}\cr \cdot& \cdot& \cdot& \cdot& \cdot& \cdot& \cdot& \cdot&
\cdot& \cdot& \cdot&&& {+}\cr &&& a_4&&& {+}& {-}&&& {+}\cr
&&&& a_5&&& {+}& {-}\cr &&&&& a_6&&& {+}\cr }\right]\cr}$$ Yet
glue links the stretched rows or columns in a tree structure
described by the following Definition.

\Proclaim Definition 3, (Weighted) Row and Column Graphs. The\/
{\rm row graph} of a matrix takes rows for vertices and connects
two by an edge of\/ {\rm weight} $w$ if they have nonzeroes in the
same $w$ columns. The row graph of matrix $G$ is $\row (G)$. Row
graphs of submatrices are subgraphs of row graphs, and so on.
The\/ {\rm column graph} is similar.\footnote{$^3$}{\rm This
concept isn't in [\Duff ] [\George ] and may be new. Other weights
can be used, for example, the inner product.}

A simple row stretching links its stretched rows by linear trees
within the row graph of its glue columns; compounded stretchings
build more elaborate trees.\footnote{$^4$}{A {\sl tree\/} is a
connected graph that breaks in two with the loss of any non-loop
edge. Equivalently, a tree has exactly one non-repeating path
between every two vertices.} Figure~11 displays the successive
trees for the compound stretching of a single row $r \rightarrow
r^S \rightarrow r^{SS} \rightarrow r^{SSS}$ in the text above.
When many rows in a matrix stretch, then the descendents of each
become separate maximally connected subgraphs in the glue's
row graph, and those subgraphs are trees. The rows in each tree
contain all the scattered pieces of some original row and can be
summed to recover that row. For column stretchings, exchange
{\sl row\/} and {\sl column\/} in this discussion.


\bigskip \vbox {\Picture 11. (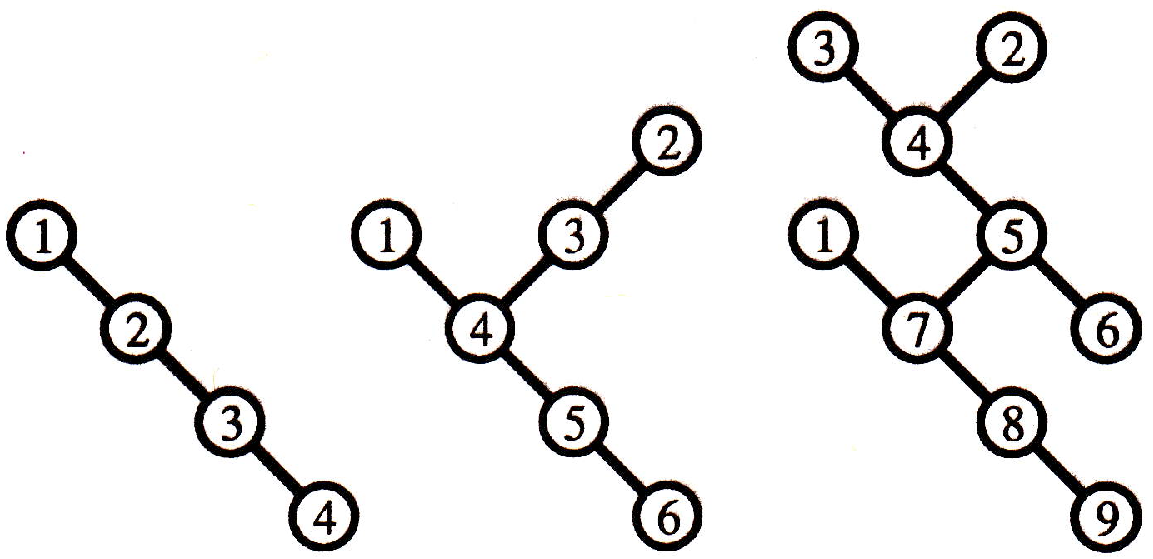 4.56 2.18) Row graphs of the
matrices stretched from a single row in the text of Section~4a.
Loop edges have been omitted. \par }\bigbreak

The trees in the row and column graphs of the glue enable the
following bound on condition numbers. The bound is particularly
pleasing because it depends on the maximum descendents of any
one row or column---the size of the largest tree---but it does
not depend on the total descendents of all rows and columns. For
example, if many rows stretch in the block fashion of
Definition~2, then the bound varies with the number of
descendents for any one row, as though only one row stretched. In
this way the bound is independent of the total growth in the size
of the matrix. 

\Proclaim Theorem~5. If $$A \rightarrow A^S \rightarrow A^{SS}
\cdots \rightarrow \A$$ is a sequence of simple row or column
stretchings\/ {\rm but not both}, and if each row or column of $A$
stretches to at most $m$ rows or columns of $A^{SS \cdots S}$,
and if Definition~2's matrices $D_i$ have the form $\sigma I$ for
the same $\sigma$, then the following choices for $\sigma$
$$\vbox {\offinterlineskip \halign {\hfil $\;$#& \quad\vrule\quad#& 
\hfil $#$\hfil \quad& \hfil $#\;$\hfil \cr $\sigma$\vrule depth6pt
height9pt width0pt&& p =1& p =\infty\cr \noalign {\hrule}
row\vrule depth3pt height15pt width0pt&& \| A \|_p / 2& \| A
\|_p\cr column\vrule depth3pt height15pt width0pt&& \| A \|_p& 
\| A \|_p / 2\cr }}$$ yield a final stretched matrix $\A$ with
bounded condition number $$\kappa_p (\A) \le c \, \kappa_p (A)$$
in which the multiplier $c$ is given below. $$\vbox
{\offinterlineskip \halign {\hfil $\;$#& \quad\vrule\quad#& 
\hfil #\hfil \quad& \hfil #$\;$\hfil \cr $c$\vrule depth6pt height9pt
width0pt&& $p =1$& $p =\infty$\cr \noalign {\hrule} row\vrule
depth3pt height15pt width0pt&& $2m-1$& $m^2$\cr column\vrule
depth3pt height15pt width0pt&& $m^2$& $2m-1$\cr }}$$ When the
sequence of stretchings is disjoint in the sense that later
stretchings do not alter the rows or columns of earlier
stretchings, then $3m$ can replace $m^2$ in this table. All these
bounds are sharp for some matrices (proof appears in
Appendix~1).

Figure~12 illustrates Theorem~5. The matrices of Figure~6 are
stretched in each of the four ways indicated by the Theorem's
tables. Either row or column stretching is performed, and the glue
is chosen to bound either the $1$-norm or the $\infty$-norm
condition numbers. Figure~12 plots the condition numbers before
and after stretching. In all cases the condition numbers increase
less than the moderate bounds allow.

\midinsert {\Picture 12. (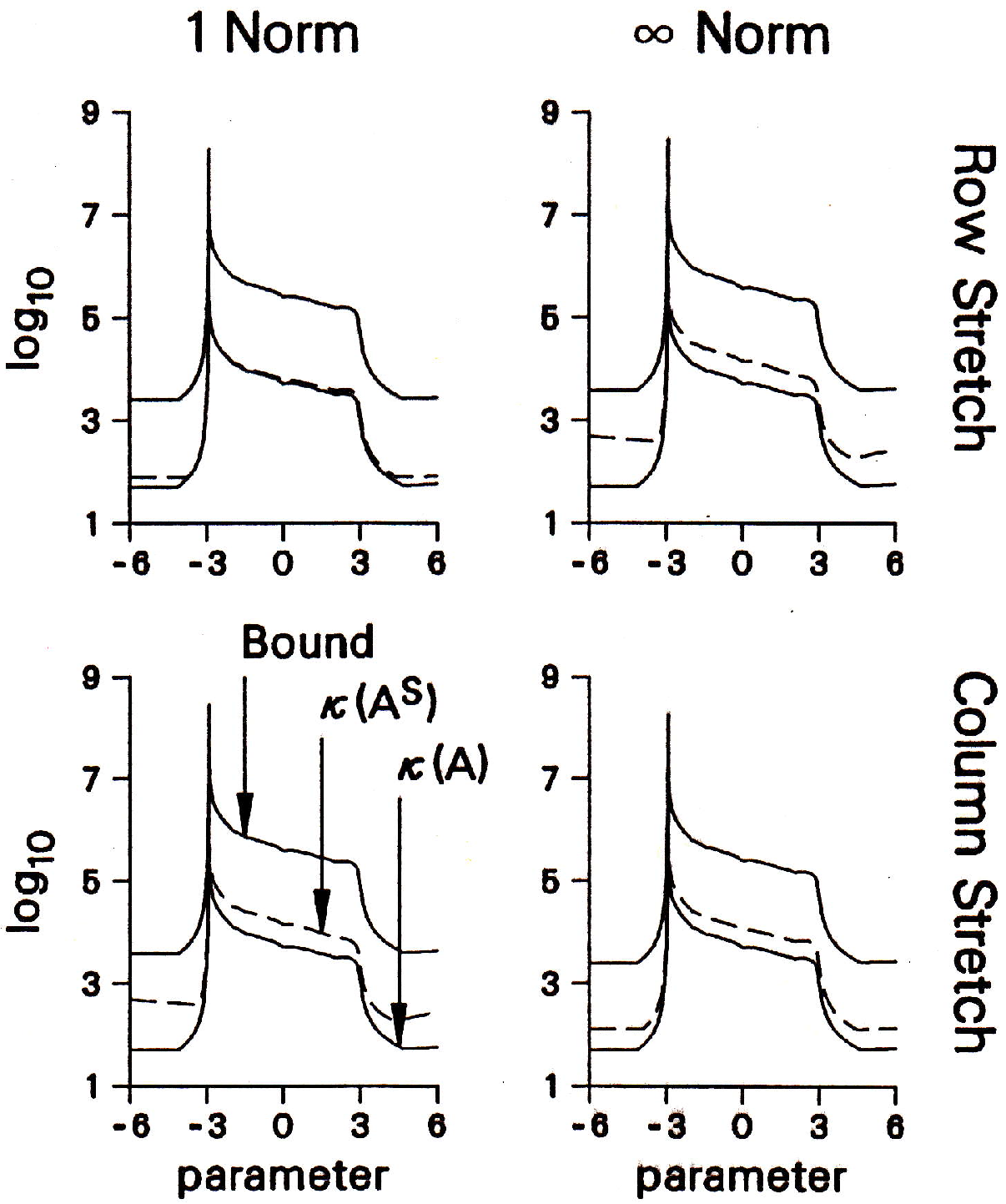 xxx 6.0) Condition numbers of the matrices in
Figure~6 (lower solid lines) and condition numbers after
stretching (dashed lines) to remove either bordering rows or
columns. Theorem~5 specifies glue that bounds (upper solid lines)
either the $1$- or $\infty$-norm condition numbers. Appendix~2
and Section~4a explain the calculations. \par }\endinsert

\beginsubsection {4b. A Priori Accuracy}

J.~H.~Wilkinson called error bounds {\sl a priori\/} when they
guaranty accuracy without measuring residuals and the like. He
obtained a~priori bounds for solving $A x = y$ under the two
conditions discussed in Section~4. When simple stretchings are
used, Theorem~5 assures $A^S$ is well conditioned if $A$ is, and
the algorithm that solves $A^S z = y^S$ must be stable. These are
Wilkinson's conditions. Thus, the requirements for a~priori
accuracy are no more stringent with stretching than without. 

Theorem~6 presents a~priori bounds for solving $A x = y$ by
iterated simple row or column stretching with glue chosen by
Theorem~5. As a practical matter, this glue sometimes may be
unnecessary. Replacing Theorem~5's $\sigma$ with $1$ rescales
the columns in simple row stretching and has little effect on
triangular factorization with row reordering, a popular and often
stable algorithm. Nevertheless, the a priori bounds require
Theorem~5's glue. The computed stretched matrix thus differs
from the ideal matrix because $\| A \|$ must be computed with
imprecise machine arithmetic. Theorem~6 accounts for this
discrepancy. 

Elaborate vector manipulations $y \rightarrow y^S$ and $z
\rightarrow z_\S$ also generate errors, but the simplest
conveniently eliminate the need for additional error analysis. As
explained in Section~3, a simple row stretching admits several
vector stretchings, but $z \rightarrow z_\S$ must copy entries out
of $z$, that is, must {\sl gather}. Conversely, simple column
stretching admits several vector squeezings, but $y \rightarrow
y^S$ must copy entries into $y^S$, that is, must {\sl scatter}. Both
simple stretchings can accept both simple vector operations, and
when they do, then forming $A^S$ and solving $A^S z = y^S$ are
the only sources of machine arithmetic error. In this case,
relative errors in $z_\S$ bound relative errors in $z$, and
Appendix~1 combines these bounds with Theorem~5 to obtain the
following result.

\Proclaim Theorem~6. If $A \rightarrow A^S$ is stretching of a
nonsingular matrix obtained from a sequence of simple row or
column stretchings\/ {\rm but not both}, and if glue is chosen by
Theorem~5, and if the stretched matrix $\AS$ is computed in
finite precision arithmetic with unit roundoff $\epsilon$, and if the
vector operations used to solve linear equations are error-free
scatter $y \rightarrow y^S$ and gather $z \rightarrow z_\S$
operations, and if the approximate solution $\barz$ to the
computed stretched equations $\AS z = y^S$ exactly satisfies
some perturbed equations $(\AS + E) \barz = y^S$, then
$$\delta_1 := c_1 \, [(1 + \epsilon)^{n}-1] < 1 \qquad \hbox {and}
\qquad \delta_2 := { \| E \|_p \over \| \AS \|_p} < {1 - \delta_1
\over 1 + \delta_1 \vphantom {\AS}}$$ imply $$\eqalign {{\| x -
\barz_\S \|_p \over \| x \|_p} &< c_2 \, \kappa (A) \, {(\delta_1 +
\delta_2 + \delta_1 \delta_2) \over 1 - (\delta_1 + \delta_2 +
\delta_1 \delta_2)}\cr \noalign {\smallskip} &\approx c_2 \,
\kappa (A) \left( c_1 n \epsilon + \| E \|_p / \|
\AS \|_p \right) \cr}$$ \smallskip \noindent in which $c_1$ and
$c_2$ are given by the tables $$\vcenter {\offinterlineskip \halign
{\hfil $\;$#& \quad \vrule \quad#& \hfil #\hfil \quad& \hfil
#$\;$\hfil \cr  $c_1$\vrule depth6pt height9pt width0pt&& $p =
1$& $p = \infty$\cr \noalign {\hrule} row\vrule depth3pt
height15pt width0pt&& $2$& $m$\cr column\vrule depth3pt
height15pt width0pt&& $m$& $2$\cr }} \qquad \vcenter
{\offinterlineskip \halign {\hfil $\;$#& \quad\vrule\quad#& \hfil
#\hfil \quad&  \hfil #$\;$\hfil \cr $c_2$\vrule depth6pt height9pt
width0pt&&  $p =1$& $p =\infty$\cr \noalign {\hrule} row\vrule
depth3pt height15pt width0pt&& $(2m-1)^2$& $m^2$\cr
column\vrule depth3pt height15pt width0pt&& $m^3$& $2m-1$\cr
}}$$ where $n$ is the order of $A$ and each row of $A$ stretches to
at most $m$ rows of $A^S$. When the sequence of stretchings is
disjoint in that later stretchings do not alter the rows or columns
of earlier stretchings, then the tables can be replaced by the ones
below. $$\vcenter {\offinterlineskip \halign {\hfil $\;$#& \quad
\vrule \quad#& \hfil #\hfil \quad& \hfil #$\;$\hfil \cr  $c_1$\vrule
depth6pt height9pt width0pt&& $p =1$&  $p =\infty$\cr \noalign
{\hrule} row\vrule depth3pt height15pt width0pt&& $2$& $2$\cr
column\vrule depth3pt height15pt width0pt&& $2$& $2$\cr }}
\qquad \vcenter {\offinterlineskip \halign {\hfil $\;$#& \quad
\vrule \quad#& \hfil #\hfil \quad& \hfil #$\;$\hfil \cr $c_2$\vrule
depth6pt height9pt width0pt&&  $p =1$& $p =\infty$\cr \noalign
{\hrule} row\vrule depth3pt height15pt width0pt&& $(2m-1)^2$&
$3m$\cr column\vrule depth3pt height15pt width0pt&& $3m^2$&
$2m-1$\cr }}$$ Thus, if $A$ is well-conditioned, if $\epsilon$ and $\|
E \|_p / \| \AS \|_p$ are very small, and if $m$ and $n$ are not
excessively large, then $\barz_\S$ is a good approximate solution
to $A x = y$ (proof appears in Appendix~1).

\beginsection {5. Arrow Matrices}

The important class of bordered, banded matrices demonstrates
stretching's utility. These matrices have the shape of the matrix
in Figure~2, namely $$\left[ \matrix {B& C\cr R& E\cr } \right]$$ in
which $B$ is banded and the bordering rows and columns $R$ and
$C$ are dense. Bordered matrices with general, sparse $B$ occur
frequently. The banded kind of interest here sometimes are
called {\sl arrow matrices}. For them, stretching significantly
improves the solving equations by triangular factorization. In the
next section, stretching compares favorably even with algorithms
designed specifically for bordered systems.

Stretching has applications beyond arrow matrices, but more
general sparse matrices pose questions that cannot be settled by
mathematical proof. Investigation of these, like other issues
involving randomly sparse matrices, requires extensive
comparison of examples that is beyond the scope of this paper.
Arrow equations are considered because they allow precise
quantification of stretching's economies. In general, only
experience proves the effectiveness of sparse matrix methods.

Bordered matrices pose a significant dilemma in the use of
triangular factorization methods. The row or column order usually
must change to insure numerical accuracy, but when a row or
column moves out of the dense border, then the factors can
become completely dense. The computational complexity of the
dense case is an upper bound for all matrices but a severe
overestimate for many sparse ones [\Duff ] [\George ]. Special
reordering strategies that avoid creating new nonzeroes and
special data structures that manipulate only the nonzeroes yield
significant economies that can be precisely quantified for banded
matrices and some others. Theorem~7 shows banded matrices
reduce factorization complexity from $2n^3/3$ operations to $2
\ell (\ell + u + 1) n$ in which $\ell$ and $u$ are the strict lower and
upper bandwidths. When the matrix is bordered, however, then the
pessimistic dense case cannot be ruled out and in some cases is
even likely.

\Proclaim Theorem~7. An $n \times n$, dense system of linear
equations can be solved by triangular factorization with row
reordering for stability using $$\vcenter {\openup1\jot
\halign{\hfil #\quad& #\hfil \cr $2 n^3 / 3 - 2n/3$& arithmetic
operations for the factorization and\cr $2n^2-n$& operations for
the solution phase.\cr }}$$ However, if the matrix is banded with
strict lower and upper bandwidths $\ell$ and $u$, and if $\ell + u <
n$, then the operations reduce to $$\vcenter {\openup1\jot
\halign{\hfil #\quad& #\hfil \cr $2 \ell (\ell + u + 1) n - \ell (4 \ell^2
+ 6 \ell u + 3 u^2 + 6 \ell + 3 u + 2)/3$& for the factorization
and\cr $(4 \ell + 2 u + 1) n - (2 \ell^2 + 2 \ell u + u^2 + 2 \ell + u)$& 
for the solution phase\cr }}$$ (proof appears in Appendix~1). 

Stretching eliminates the possibility of catastrophe for
bordered, banded matrices by eliminating the border. With the
customary row reordering it is sufficient to remove only the
dense rows. This can be done by the simple row stretching of
Definition~2. A row and column reordering then gives the
stretched matrix a banded structure for which factorization with
row reordering is clearly efficient. Both the stretching and the
reordering depend on the following blocking of the rows and
columns.

\Proclaim Theorem~8. This row and column partitioning makes a
banded matrix into a block-bidiagonal one. For a matrix of order
$n$ with strict lower and upper bandwidths $\ell$ and $u$, and with
$0 < \ell + u < n$, the columns and rows partition into blocks of the
following size. {\def \Ls {,{\kern 0.75em}} $$\vcenter
{\openup1\jot \halign {\hfil #\hquad \quad& $#$\hfil \cr columns& 
\listfive {a + u} {\ell + u} {\ldots} {\ell + u} {\ell + c} \endlist\cr 
rows& \listsix {a} {u + \ell} {u + \ell} {\ldots} {u + \ell} {c} \endlist
\cr \noalign {\medskip} \multispan2 {\hfil $0 \le a \le \ell \qquad 0
\le c \le u \qquad 0 < a + c$ \hfil }\cr }}$$}The block-column
di\-men\-sion is $m = \lceil n / (\ell + u) \rceil$, and the
block-row dimension is $m + 1$ or $m$ (since one of $a$ or $c$ may
be zero). Moreover, the upper diagonal blocks are lower
triangular and the lower diagonal blocks are upper triangular
(proof appears in Appendix~1). 

The partitioning of Theorem~8 applied to the banded portion of an
arrow matrix results in a blocked matrix. $$\left[ \matrix {B& C\cr 
R& E\cr } \right] = \left[ \matrix {L_1&&&&& C_1\cr U_1& 
L_2&&&& C_2\cr & U_2& L_3&&& C_3\cr && U_3& \Ddots&& 
\Vdots\cr &&& \Ddots& L_m& C_m\cr &&&& U_m& C_{m+1}\cr 
R_1& R_2& R_3& \cdots& R_m& E\cr } \right]$$ The row blocks
containing $L_1$ and $U_m$ have row orders $a$ and $c$ which the
Theorem allows to be zero, but no harm results from including null
blocks in the display. Both the picture and the Theorem assume
$B$ is square, in other words, the bordering rows number the
same as the bordering columns, and the banded portion ends as
shown where the bordering rows and columns intersect. Simple
row stretching replaces this matrix by $$\left[ \matrix
{L_1&&&&& C_1\cr U_1& L_2&&&& C_2\cr & U_2&  L_3&&&
C_3\cr && U_3& \Ddots&& \Vdots\cr &&& \Ddots& L_m&  C_m\cr
&&&& U_m& C_{m+1}\cr R_1&&&&&& -D_1\cr &  R_2&&&&&
+D_1& -D_2\cr && R_3&&&&& +D_2& \Ddots\cr &&& 
\Ddots&&&&& \Ddots& -D_{m-1}\cr &&&& R_m& E&&&& 
+D_{m-1}\cr } \right]$$ and then applies a perfect shuffle to the
blocks of rows and columns. The $j^{th}$ block of new rows goes
after the $j^{th}$ block of old rows, and similarly for columns. The
result is arresting, even beautiful. $$\left[ \matrix
{L_1&&&&&&&&& C_1\cr R_1& -D_1&&&&&&&& 0\cr U_1&& 
L_2&&&&&&& C_2\cr & +D_1& R_2& -D_2&&&&&& 0\cr && 
U_2&& L_3&&&&& C_3\cr &&& +D_2& R_3& \Ddots&&&& 
\Vdots\cr &&&& U_3&& \Ddots&&& \cr &&&&& \Ddots& \Ddots& 
-D_{m-1}&& 0\cr &&&&&& \Ddots&& L_m& C_m\cr &&&&&&& 
+D_{m-1}& R_m& E\cr &&&&&&&& U_m& C_{m+1}\cr } \right]$$
Theorem~8 goes to some trouble to insure this matrix is as good
as it looks. The banded portion is seamless with uniform strict
lower and upper bandwidths $d + \ell$ and $u$, in which $d$ is the
depth of the border and $\ell$ and $u$ are the bandwidths in the
original arrow matrix.

When the stretched matrix is used to solve equations, then the
computational complexity varies linearly with the size of the
banded portion of the original matrix, as in the purely banded
case. Theorems~7 and~9 supply the following operation counts
for the factorization and solution phases, respectively. $$\vcenter
{\offinterlineskip \halign {\hfil $\displaystyle #$& \hfil #\hfil&
$\displaystyle #$\hfil \cr \hbox {banded}\hfil&  \vrule depth9pt
height9pt& \hfil \hbox {additional operations when bordered and
stretched}\cr \noalign {\hrule}\cr & \vrule depth0pt height6pt\cr
\noalign {\smallskip} 2 \ell (\ell + u + 1) n&  $\hquad + \hquad$&
2d(2 d + 3 \ell + u + 1) n \; + \; 2d (d + \ell) (2d + \ell + u + 1) \left
\lceil {n \over \ell + u} \right \rceil\cr \noalign {\medskip} (4\ell +
2u + 1) n& $\hquad + \hquad$& 4dn \; + \; d(4d + 4\ell + 2u + 1 )
\left \lceil {n \over \ell + u} \right \rceil\cr }}$$ In these formulas,
$n + d$ is the size of the unstretched arrow matrix and $n$ is the
size of its banded portion.

\Proclaim Theorem~9. An order $n + d$, bordered, banded system
of linear equations, whose coefficient matrix has $d$ dense rows
and columns in the bordering portion and has strict lower and
upper bandwidths $\ell$ and $u$ in the $n \times n$ banded portion,
where $0 < \ell + u < n$, can be solved by simple row stretching
and triangular factorization with row reordering for stability in
$$\displaylines {(4 d^2 + 6 d \ell + 2 d u + 2 \ell^2 + 2 \ell u + 2 d +
2 \ell) N \cr {}- (d + \ell) (13 d^2 + 14 d \ell + 12 d u + 4 \ell^2 + 6
\ell u + 3 u^2 + 9 d + 6 \ell + 3 u + 2) / 3\cr }$$ arithmetic
operations for the factorization and $$\displaylines {(4 d + 4 \ell +
2 u + 1) N \cr {}- (2 d^2 + 4 d \ell + 2 d u + 2 \ell^2 + 2 \ell u + u^2 +
2 d + 2 \ell + u)\cr }$$ operations for the solution phase, in which
$$N = n + d \left \lceil {n \over \ell + u} \right \rceil$$ is the
size of the stretched matrix (proof appears in Appendix~1). 

\removelastskip \beginsection {6. Deflated Block Elimination}

Bordered matrices occur with sufficient frequency to receive
special treatment. Although stretching is a general sparse matrix
method, it compares favorably with specialized algorithms for
bordered equations. Chan's deflated block elimination [\Chan ]
solves bordered, banded systems with computational complexity
near stretching's, but its numerical accuracy can be much worse.
Stretching therefore is more reliable for arrow matrices, and also
more versatile. Deflated block elimination may be suited for other
contexts, however, and should not be judged solely by this
comparison.

Pure block elimination employs a representation of the inverse
$${\left[ \matrix {B& C \cr R& E} \right]}^{- 1} = \left[ \matrix {I& 
- B^{-1} C \cr 0& I} \right] \left[ \matrix {B^{- 1}& 0\cr 0& (E - R
B^{-1} C)^{- 1}} \right] \left[ \matrix {I& 0\cr - R B^{-1}& I}
\right]$$ obtained by inverting the block factorization $$\left[
\matrix {B& C \cr R& E} \right] = \left[ \matrix {I& 0\cr R B^{-1}& 
I} \right] \left[ \matrix {B& 0\cr 0& E - R B^{-1} C} \right] \left[
\matrix {I& B^{-1} C \cr 0& I} \right] .$$ The factorization is a
recipe for applying the inverse without evaluating the factors.
There are two phases. Steps in the first phase depend on the
matrix alone, those in the second repeat for each right side.
Table~1 lists the steps and counts the arithmetic operations to
solve $$\left[ \matrix {B& c\cr r^t& \alpha_0\cr } \right] \left[
\matrix {v_*\cr \beta_*\cr } \right] = \left[ \matrix {v_0\cr 
\beta_0\cr } \right]$$ in which $r^t$ and $c$ replace $R$ and $C$ to
indicate a single bordering row and column. The weakness of block
elimination is the need to solve equations with a coefficient
matrix, $B$, that can be badly conditioned or singular even when
the larger matrix is neither. 

 \midinsert{
 \table 1. {Factorization and solution phases of block elimination
with operation counts. $B$ has order $n$ and strict lower and
upper bandwidths $\ell$ and $u$, $\ell + u < n$. The costs of
factoring $B$ and applying $B^{-1}$ are from Theorem~7. Terms
independent of $n$ are omitted. Section~6 provides further
explanation.} {\vbox {\openup0.5\jot \halign {\quad#\hfil& \qquad
\hfil #\hfil \quad \cr \omit \hfil step\hfil& operations\cr \noalign
{\medskip \hrule \medskip} construct a triangular& $2 \ell (\ell +
u + 1) n$\cr \quad factorization of $B$\cr $u_1 := B^{- 1} c$& $(4
\ell + 2 u + 1) n$\cr $\alpha_1 := \alpha_0 - r^t u_1$& $2n$\cr 
\noalign {\medskip \hrule \medskip} $v_1 := B^{- 1} v_0$& $(4 \ell
+ 2 u + 1) n$\cr $\beta_1 := \beta_0 - r^t v_1$& $2n$\cr $\beta_*
:= \beta_1 / \alpha_1$\cr $v_* := v_1 - u_1 \beta_*$& $2 n$\cr 
\noalign {\medskip \hrule }}}}\par }\endinsert

Deflated block elimination [\Chan ] attempts to correct the
deficiency of the original method. The rationale for the deflated
algorithm appears to be the following. It has been developed for
matrices with one dense row and column, that is, $B$ must be a
maximal submatrix. Maximal submatrices of nonsingular matrices
have null spaces dimensioned at most $1$, so if the entire matrix
is well-conditioned then it is inferred either $B$ is
well-conditioned too or has a well-separated, smallest singular
value. Whenever $B^{- 1}$ must be applied to a vector, a
component that lies in or near the space corresponding to the
smallest singular value might be separated from the product. A
modification of the block factorization recipe manipulates these
decomposed vectors without loss of accuracy. There results a
more complicated algorithm that gives accurate results even
when $B$ has a small singular value. When not, then there is no
harm in evaluating the more elaborate formulas. The weakness of
the algorithm is the need to estimate the smallest singular value
of $B$ and the associated singular vectors. 

Table~2 lists the steps and counts the arithmetic operations
performed by the deflated algorithm to solve the same equations
as Table~1. The Table makes the following choices among the
algorithm's many variations. First, several approximations might
be made to the the smallest singular value and its singular
vectors. Table~2 uses the ``orthogonal projector'' estimates in
the final row of Chan's Table~3.1 [\Chan,~p.~125]. This choice
appears to be the most economical. Second, there is some
apparent variation in computing the vector decompositions.
Table~2 uses ``Algorithm NIA'' [\Chan,~p.~126] without the final
step. Chan does not recognize NIA's final step unnecessarily
applies a linear transformation to a vector invariant for the
transformation. He omits the step for other reasons
[\Chan,~p.~130]. Third, $B$ might be factored and $B^{- 1}$
applied by several means. Table~2 employs the Doolittle
triangular decomposition with row reordering for stability, as
does Chan [\Chan,~p.~130].

 \midinsert{
 \table 2. {Factorization and solution phases of deflated block
elimination, with cross-references to the original notation\/ {\rm
[\Chan ]}, and with operation counts. $B$ has order $n$ and strict
lower and upper bandwidths $\ell$ and $u$, $\ell + u < n$. The costs
of factoring $B$ and applying $B^{-1}$ are from Theorem~7. Terms
independent of $n$ are omitted. Notation $e_k$ is column $k$ of an
identity matrix. Section~6 provides further explanation.} {\vbox
{\openup0.5\jot \halign {\quad\hfil #\hfil& \quad#\hfil&
\qquad\hfil #\hfil \quad\cr \noalign {\medskip} original& \omit
\hfil step\hfil& operations\cr \noalign {\medskip \hrule
\medskip} $A$& construct a triangular& $2 \ell (\ell + u + 1) n$\cr
& \quad factorization of $B$\cr $k$& choose $k$, the index of&
$n$\cr & \quad the smallest pivot\cr $\psi$& $u_1 := B^{- t} e_k /
\| B^{- t} e_k \|_2$& $(4 \ell + 2 u + 4) n$\cr & $u_2 := B^{- 1} u_1$&
$(4 \ell + 2 u + 1) n$\cr $\delta$& $\alpha_1 := 1 / \| u_2 \|_2$& 
$2n$\cr $\phi$& $u_3 := u_2 \alpha_1$& $n$\cr $c_b$& $\alpha_2 :=
u_1^t c$& $2n$\cr $v_D$& $u_4 := B^{- 1} (c - u_1 \alpha_2)$& $(4
\ell + 2 u + 3) n$\cr $h_2$& $\alpha_3 := \alpha_0 - r^t u_4$& 
$2n$\cr & $\alpha_4 := r^t u_3$& $2n$\cr $D$& $\alpha_5 :=
\alpha_2 \alpha_4 - \alpha_1 \alpha_3$\cr \noalign {\medskip
\hrule \medskip} $c_f$& $\beta_1 := u_1^t v_0$& $2 n$\cr $w_D$& 
$v_1 := B^{- 1} (v_0 - u_1 \beta_1)$& $(4 \ell + 2 u + 3) n$\cr 
$h_1$& $\beta_2 := \beta_0 - r^t v_1$& $2n$\cr $h_3$& $\beta_3 :=
\alpha_2 \beta_2 - \alpha_3 \beta_1$\cr $h_4$& $\beta_4 :=
\alpha_4 \beta_1 - \alpha_1 \beta_2$\cr $x$& \omit \quad \rlap
{$v_* := v_1 + (u_3 \beta_3 - u_4 \beta_4) / \alpha_5$}& $4 n$\cr 
$y$& $\beta_* := \beta_4 / \alpha_5$\cr \noalign {\medskip \hrule
}}}}\par }\endinsert 

The numerical performance of deflated block elimination varies
with the quality of approximations to small singular values and
their vectors. The approximations made by Table~2 depend on
``the smallest pivot having the magnitude of the smallest singular
value, which is definitely not valid in general, but which is shown
empirically and theoretically to be valid in practice''
(paraphrasing [\Chan,~p.~124]). It is well known ``there is no
correlation between small pivots and ill-conditioning''
[\Golub,~p.~63], but in light of Chan's remarks it is surprising
Figure~13 shows his approximations are invalid for the
commonplace matrices of Figure~6. In this case deflated block
elimination degenerates to pure block elimination, and Figure~14
shows the two methods have nearly identical, spectacularly large
errors. Better approximations to the singular values and vectors
could remedy this, but their costs would favor stretching even
more.

 \pageinsert {\vfil
 {\Picture 13. (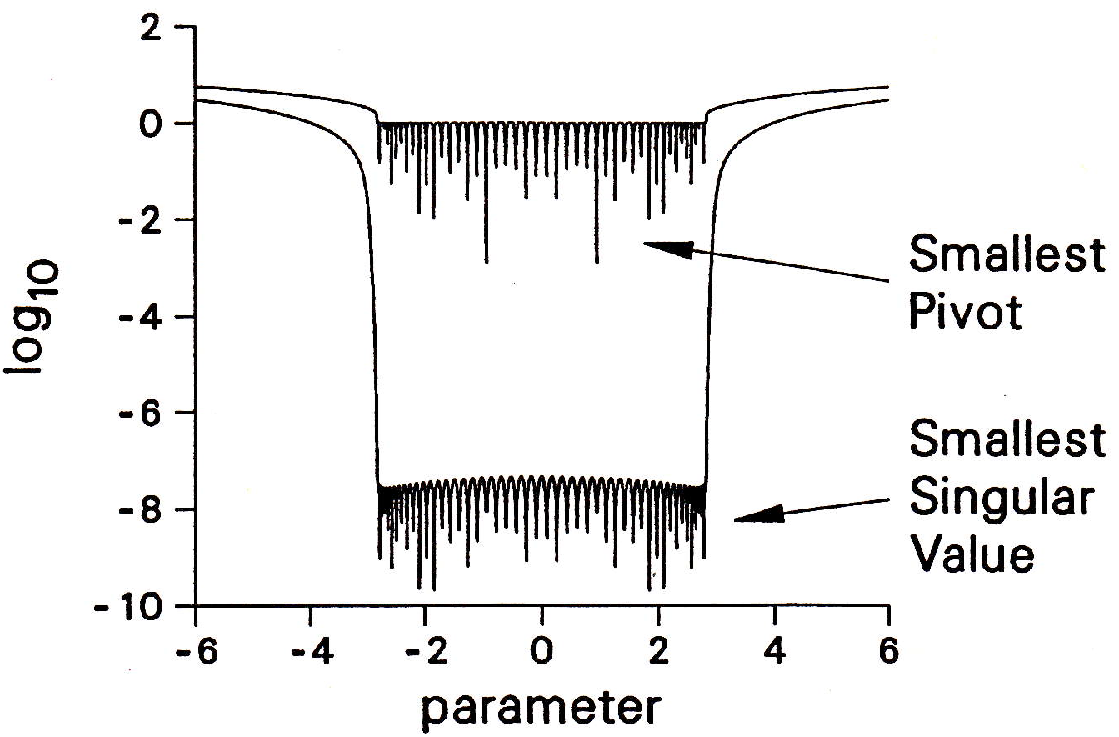 xxx 3.5) Smallest pivot and singular value for the banded
portion of the matrices in Figure~6. Table~2's version of deflated
block elimination assumes the pivot and singular value have the
same magnitude ``which is definitely not valid in general, but
which is shown empirically and theoretically to be valid in
practice'' {\rm [\Chan,~p.~124].} \par }
 \vfil
 {\Picture 14. (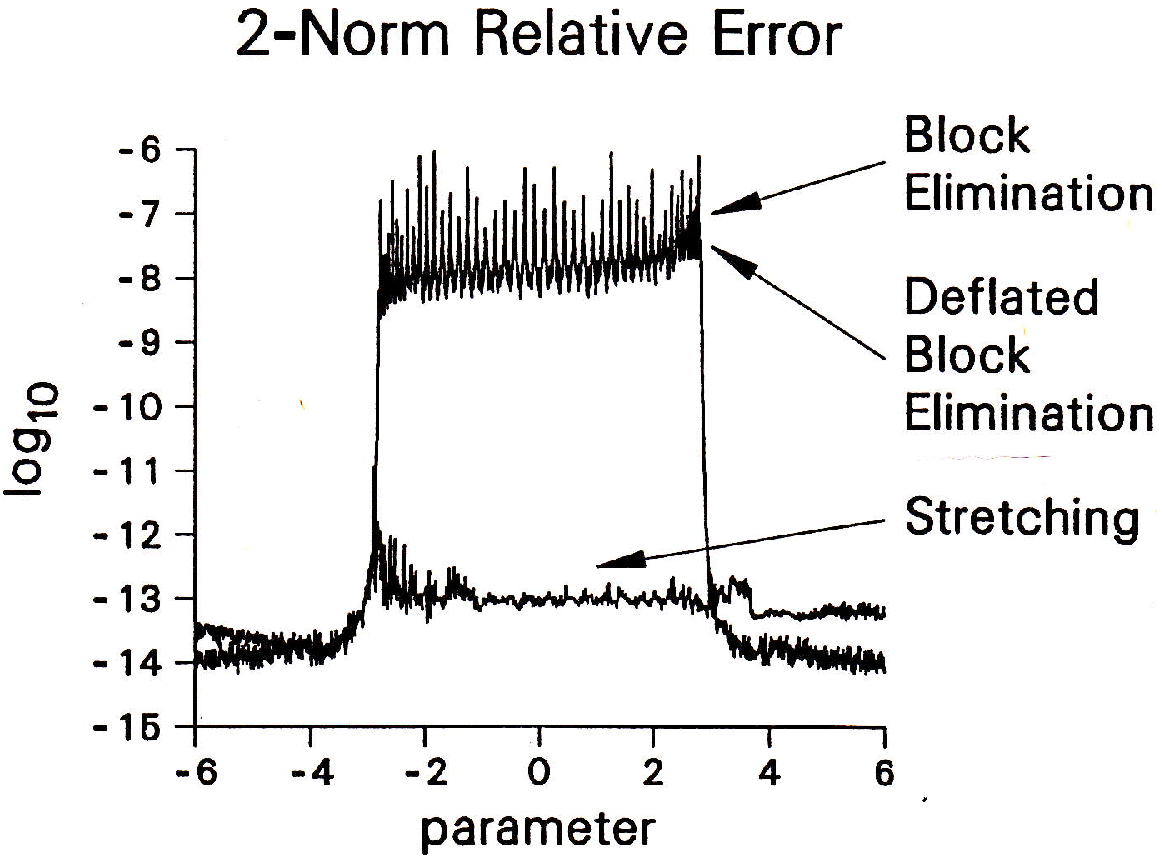 xxx 3.5) Maximum $2$-norm relative errors for equations
$A x = y$ with $20$ different $y$'s solved by block elimination,
deflated block elimination, and simple row stretching. The
parameterized coefficient matrices $A$ are those of Figure~6.
The elimination methods cannot be distinguished at this
plotting resolution. The stretching data also appears in Figure~9.
Appendix~2 and Section~6 explain the calculations. \par }
 \vfil}\endinsert 

The arithmetic costs for Table~2's version of deflated block
elimination and for simple row stretching are nearly the same, but
stretching's are generally smaller. Operations in the
``factorization'' phase are $$\vcenter {\openup1\jot \halign {\hfil
#& \hquad \quad $n \, [ \, \hfil #$& $\displaystyle {{} + #}$\hfil \cr 
Block Elimination& 2 \ell^2 +2 \ell u& \phantom{1}6 \ell + 2 u +
\phantom{1}3 \, ]\cr Simple Row Stretching& \ldots \;& 10 \ell + 2
u + \phantom{1}8 + {6 \ell + 6 \over \ell + u} \, ]\cr Deflated Block
Elimination& \ldots \;& 14 \ell + 6 u + 18 \, ] \, \cr }}$$ and in the
``solution'' phase they are $$\vcenter {\openup1\jot \halign {\hfil
#& \hquad \quad $n \, [ \, \hfil #$& $\displaystyle {{} + #}$\hfil \cr 
Block Elimination& 4 \ell +2 u& 5 \, ]\cr Simple Row Stretching& 
\ldots \;& 9 + {2 \ell + 7 \over \ell + u} \, ]\cr Deflated Block
Elimination& \ldots \;& 11 \, ] \, \cr }}$$ The operation counts for
the block elimination algorithms are from Tables~1 and~2. Those
for simple row stretching are from Theorem~9 with the size, $N$,
of the stretched matrix estimated as follows. $$N \; = \; n + \left
\lceil {n \over \ell + u} \right \rceil \; \approx \; n \, {\ell + u + 1
\over \ell + u}$$ By this estimate stretching has the advantage in
the factorization phase whenever $$10 \ell + 2u + 8 + {6 \ell + 6
\over \ell + u} < 14 \ell + 6 u + 18 \quad \Longleftrightarrow \quad
0 < \ell + u$$ and in the solution phase whenever $$9 + {2 \ell + 7
\over \ell + u} < 11 \quad \Longleftrightarrow \quad 3 < u.$$
Stretching is therefore more economical than deflated block
elimination for all but the smallest bandwidths.

A more significant advantage of stretching is the ability to treat
arbitrarily many bordering rows and columns. Theorem~9 already
reports stretching's arithmetic costs for bordered, banded
matrices with many borders. In contrast, deflated block
elimination has been developed for only one bordering row and
column. It cannot easily be applied recursively to a succession of
single borders, or be extended to multiple borders in some more
direct way, because in both cases it encounters the difficult
problem of submatrices with multiple small singular values.

\beginsection {7. Antecedents}

Some analytic methods make some ordinary differential
equations more amenable to numerical solution, and inspire
matrix stretching. The following description il\-lustrates the
intellectual leap to the algebraic process and may suggest
additional stretchings. This section views differential
equations from the standpoint of software engineering. 

For ease of notation the equations are assumed to include only
1st-order differentials and only nondifferential boundary
conditions. Thus, the equations are $$\openup2\jot \displaylines
{F(u, u') = 0\cr G(u(a)) = 0\qquad H(u(b)) = 0\cr F \,\colon\, \R^m
\times \R^m \rightarrow \R^m\qquad (G \times H) \,\colon\, \R^m
\rightarrow \R^m \qquad u \,\colon [a,b] \rightarrow \R^m \cr }$$ in
which $F$ defines the system of $m$ differential equations while
$G$ and $H$ enforce the boundary conditions. $\R$ is the set of
real numbers.

A discrete approximate solution, $u_k \approx u(x_k)$, is
determined at $n+1$ points $$a = x_0 < x_1 < x_2 < \ldots <x_n =
b$$ by algebraic (or rather, nondifferential) equations. $$\vcenter
{\openup2\jot \halign {\hfil $\displaystyle #$\hfil& 
\qquad$\displaystyle #$& \qquad $\displaystyle #$\cr F_k := F
\left( {u_k + u_{k-1} \over 2}, {u_{k} - u_{k-1} \over x_{k} -
x_{k-1}} \right) = 0& G(u_0) = 0& H(u_n) = 0\cr k = 1, 2, \ldots,
n\cr }}$$ When these discrete equations are solved by methods
like Newton's, then matrix equations must be solved to obtain the
Newton corrections. The matrices are Jacobian matrices for the
ensemble of functions above. Ordering the variables and
equations in the natural way $$\listfive u_0 u_1 u_2 {\ldots} u_n
\endlist \qquad \listsix G F_1 F_2 {\ldots} F_n H \endlist$$ gives
the matrices a banded structure $$J = \left[ \matrix {A_0&
\phantom {B_0}& \phantom {B_0}&  \phantom {B_0}& \phantom
{B_0}\cr B_0& A_1\cr \phantom {B_0}&  B_1& A_2\cr \phantom
{B_0}&& B_2& \Ddots\cr \phantom {B_0}&&& \Ddots& A_n\cr
\phantom {B_0}&&&& B_n\cr } \right]$$ in which $A_0$, $A_{j>0}$,
$B_{k < n}$ and $B_n$ are Jacobian matrices for $G$, $F_j$,
$F_{k+1}$ and $H$ with respect to $u_0$, $u_j$, $u_k$ and $u_n$,
respectively. The matrix $J$ is square because $A_0$ and $B_n$
may not be. For example, if all the boundary conditions are
applied at the left endpoint then $H$ and $B_n$ are vacuous. The
banded structure allows the linear equations to be solved by the
efficient matrix factorization process of Theorem~7. 

Arrow matrices occur when parameters and constraints
accompany the differential equations $$\openup2\jot
\displaylines {F(u, \lambda, u') = 0 \qquad 1 - \int_a^b (w^tu)^2 =
0\cr G(u(a)) = 0 \qquad H(u(b)) = 0\cr }$$ in which $w$ is a vector
of coefficients that define the integral constraint. With the
dependent parameter $\lambda$ placed after the discrete
variables, and with a discrete analogue of the constraint $$E (u_0,
u_1, u_2, \ldots, u_n) \; := \; 1 - \sum_{k=1}^n (x_k - x_{k-1}) \,
{(w^t u_k)^2 + (w^t u_{k-1})^2 \over 2} \; = \; 0$$ placed after the
other equations, the matrices acquire borders $$J = \left[ \matrix
{A_0& \phantom {B_0}& \phantom {B_0}& \phantom {B_0}&
\phantom {B_0}& 0\cr  B_0& A_1&&&& c_1\cr \phantom {B_0}&
B_1& A_2&&& c_2\cr  \phantom {B_0}&& B_2& \Ddots&&
\Vdots\cr \phantom {B_0}&&&  \Ddots& A_n& c_n\cr \phantom
{B_0}&&&& B_n& 0\cr r_0^t&  r_1^t& r_2^t& \ldots& r_n^t& 0\cr }
\right]$$ in which $c_k$ and $r_k^t$ are the Jacobian matrices of
$F_k$ and $E$ with respect to $\lambda$ and $u_k$, respectively.
Section~5 shows arrow matrix equations can be solved efficiently
after stretching to remove the border.

However, the banded structure is traditionally recovered by
transforming the differential system rather the algebraic
equations, as follows. Differentiating the parameter and the
constraint $$\openup2\jot \displaylines {F(u, \lambda, u') = 0
\qquad \lambda' = 0\qquad v' - \| u \|_2^2 = 0\cr G(u(a)) = 0\quad
v(a) = 0\qquad H(u(b)) = 0\quad v(b) = 1\cr }$$ produces an
un-parameterized, un-constrained differential system whose
discretization again results in banded Jacobian matrices. These
differential transformations are familiar simplifying devices in
the theory of differential systems. They find widespread
application in numerical software [\Kee ] and special mention in
numerical textbooks [\Kellera,~p.~47]. Simple row and column
stretchings correspond to differentiating the constraint and the
parameter, respectively. 

The correspondence between differentiating and stretching is not
precise because the discretization intervenes. The analytic
transformation (differentiating) precedes the discretization,
while the algebraic one (stretching) follows. $$\matrix {\vcenter
{\halign {\hfil #\hfil \cr constrained, parameterized\cr differential
system\cr }}& \buildrel \hbox {\sl differentiate} \over {\hbox
to1.25in {\rightarrowfill}}& \vcenter {\halign {\hfil #\hfil \cr 
differential system\cr }}\cr \llap {\hbox {\sl discretize}\quad}
\Bigg\downarrow& \phantom {\vrule height24pt depth18pt}& 
\Bigg\downarrow \rlap {\quad \hbox {\sl discretize}}\cr \vcenter
{\halign {\hfil #\hfil \cr bordered\cr banded matrix\cr }}& \buildrel
\hbox {\sl stretch} \over {\hbox to1.25in {\rightarrowfill}}& 
\vcenter {\halign {\hfil #\hfil \cr banded matrix\cr }} \cr }$$ The
two paths through the diagram may not lead to the same Jacobian
matrix. If the discretization employs a high order approximation
to $\lambda'$, for example, then the matrix rows for the equation
$\lambda' = 0$ unnecessarily contain more than the two nonzeroes
inserted by simple column stretching. The differential
transformation therefore results in a Jacobian matrix indirectly
chosen and likely suboptimal. Moreover, any change to the
differential system can subtly perturb the discretization. The
integral constraint may yield slightly different discrete
approximations in its differentiated form. In contrast, stretching
guarantees efficient solution of the matrix equations---however
derived---and thus allows the most appropriate discretization of
the differential system. 

The originality of the algebraic approach can be appreciated by
considering related problems for which there is no convenient
analytic interpretation, and consequently, to which nothing
approximating stretching has been applied. These problems are
parameterized equations whose solution is desired as a function
of the parameter. $$\openup2\jot \displaylines {F(u,  \lambda,
u')\cr G(u(a)) = 0 \qquad H(u(b)) = 0\cr F \,\colon\, \R^m \times
\R^m \rightarrow \R^m\qquad (G \times H) \,\colon\, \R^m
\rightarrow \R^m \qquad u \,\colon [a,b] \rightarrow \R^m \cr }$$ In
this case the discrete equations are underdetermined $$\vcenter
{\openup2\jot \halign {\hfil $\displaystyle #$\hfil& 
\qquad$\displaystyle #$& \qquad $\displaystyle #$\cr F_k := F
\left( {u_k + u_{k-1} \over 2}, \lambda, {u_{k} - u_{k-1} \over x_{k}
- x_{k-1}} \right) = 0& G(u_0) = 0& H(u_n) = 0\cr k = 1, 2, \ldots,
n\cr }}$$ and the extra variable, $\lambda$, gives the Jacobian
matrices more columns than rows. $$J = \left[ \matrix {A_0&
\phantom {B_0}& \phantom {B_0}& \phantom {B_0}&  \phantom
{B_0}& 0\cr B_0& A_1&&&& c_1\cr \phantom {B_0}&  B_1&
A_2&&& c_2\cr \phantom {B_0}&& B_2& \Ddots&&  \Vdots\cr
\phantom {B_0}&&& \Ddots& A_n& c_n\cr \phantom {B_0}&&&&
B_n& 0\cr } \right] \qquad \vec u = \left[ \matrix {u_0\cr u_1\cr
u_2\cr \vdots\cr u_n\cr \lambda\cr } \right]$$ Mild assumptions
guarantee a differentiable curve of discrete solutions, $\vec u$,
which can be traced in a variety of ways. 

Keller [\Kellerb ] appears to be the first to consider the numerical
problem of following the solution curve through $\R^{n+2}$. He
locally parameterizes the curve as $\vec u (\sigma)$ in which
$\sigma$ is a local approximation to arclength. The points on the
curve near a known point $\vec u (0)$ are specified by
augmenting the discrete equations with the projection equation
$$\vec \tau \cdot \left[ \vec u (\sigma) - \vec u (0) \right] -
\sigma = 0 \qquad \vec \tau = \left[ \matrix {\tau_0\cr \tau_1\cr
\tau_2\cr \vdots\cr \tau_n\cr \mu\cr } \right]$$ in which $\cdot$
is the vector dot product and the ancillary vector $\vec \tau$
approximates a tangent to the curve (equivalently, a right null
vector of $J$) at the known solution. Given the known solution at
$\sigma = 0$, another is found by solving the augmented
equations for some $\sigma > 0$, then the curve is
re-parameterized, and the process is repeated. Only sufficiently
small $\sigma$ can be expected to uniquely determine $\vec u
(\sigma)$, and the tangent must be approximated in some way,
but the details of Keller's {\sl pseudo-arclength continuation
method\/} are of no concern here.

If variants of Newton's method are used to solve the augmented
equations, then the correction equations again feature arrow
matrices. $$J = \left[ \matrix {A_0& \phantom {B_0}& \phantom
{B_0}& \phantom {B_0}& \phantom {B_0}& 0\cr  B_0& A_1&&&&
c_1\cr \phantom {B_0}& B_1& A_2&&& c_2\cr  \phantom {B_0}&&
B_2& \Ddots&& \Vdots\cr \phantom {B_0}&&&  \Ddots& A_n&
c_n\cr \phantom {B_0}&&&& B_n& 0\cr \tau_0^t&  \tau_1^t&
\tau_2^t& \ldots& \tau_n^t& \mu \cr } \right]$$ Keller [\Kellerb ]
suggests pure block elimination for these matrices, while Chan
[\Chan ] appears to have developed deflated block elimination
with this problem in mind. These methods are discussed in
Section~6. A stretching process would be preferable but isn't
employed---apparently because the bordered matrices are
conceived in a discrete rather than an analytic context. 

No doubt many devices have been used over the years to trade
inconvenient matrices for more convenient, larger matrices. Until
this paper, however, only the {\sl capacitance matrix meth\-od\/}
has received much attention. It obtains larger matrices by domain
embedding of partial differential equations, and views the
enlarged problems as perturbed ones susceptible to the
Sherman-Morrsion-Woodbury formula [\Golub ] [\Householder ].
Buzbee, Dorr, George and Golub [\Buzbee ] explain this approach
and provide references to earlier work. The latitude in applying
the Woodbury formula produces variations in numerical accuracy
which continue to be of research interest [\Dryja ], but the matrix
enlargement process has not been generalized. The capacitance
matrix method thus retains the flavor and terminology of domain
embedding. Presumably, it can be recast as a matrix {\sevenrm
\one(S) \two(T) \one(R) \two(E) \one(T) \two(C) \one(H) \two(I)
\one(N) \two(G) }.

\vfil \eject \beginsection {References}

\item{[\Buzbee ]} B.~L.~Buzbee, F.~W.~Dorr, J.~A.~George and
G.~H.~Golub, {\sl The direct solution of the discrete Poisson
equation on irregular regions}, SIAM Journal on Numerical
Analysis, {\bf 8}:722--736 (1971).

\item{[\Chan ]} T.~F.~Chan, {\sl Deflation techniques and
block-elimination algorithms for solving bordered singular
systems}, SIAM Journal on Scientific and Statistical Computing,
{\bf 5}:121--134 (1984).

\item{[\Dongarra ]} J.~Dongarra, J.~R.~Bunch, C.~B.~Moler and
G.~W.~Stewart, {\sl LINPACK Users Guide}, Society for Industrial
and Applied Mathematics, Philadelphia, 1978.

\item{[\Dryja ]} M.~Dryja, {\sl A finite element-capacitance
matrix method for the elliptic problem}, SIAM Journal on
Numerical Analysis, {\bf 4}:671--680 (1983).

\item{[\Duff ]} I.~S.~Duff, A.~M.~Erisman and J.~K.~Reid, {\sl Direct
Methods for Sparse Matrices}, Oxford University Press, Oxford,
1986.

\item{[\George ]} A.~George and J.~W.~H.~Liu, {\sl Computer
Solution of Large Sparse Positive-Definite Systems},
Prentice-Hall, Englewood Cliffs, 1981.

\item {[\Golub ]} G.~H.~Golub and C.~F.~Van Loan, {\sl Matrix
Computations}, Johns Hopkins University Press, Baltimore, 1983.
Page references are to the paperback form of this first edition.
The second edition is much revised.

\item {[\Kee ]} R.~J.~Kee, L.~R.~Petzold, M.~D.~Smooke and
J.~F.~Grcar, {\sl Implicit Methods in Combustion and Chemical
Kinetics Modeling}, in Multiple Time Scales, edited by
J.~U.~Brackbill and B.~I.~Cohen, Academic Press, Orlando, 1985.

\item {[\Kellera ]} H.~B.~Keller, {\sl Numerical Solution of Two
Point Boundary Value Problems}, Society for Industrial and
Applied Mathematics, Philadelphia, 1976.

\item {[\Kellerb ]} H.~B.~Keller, {\sl Numerical Solution of
Bifurcation and Nonlinear Eigenvalue Problems}, in Applications of
Bifurcation Theory, edited by P.~H.~Rabinowitz, Academic Press,
1977, 359--384.

\item {[\Householder ]} A.~S.~Householder, {\sl The Theory of
Matrices in Numerical Analysis}, Blaisdell, New York, 1964.
Reprinted by Dover, New York, 1975.

\item{[\Parter ]} S.~Parter, {\sl The use of linear graphs in Gauss
elimination}, SIAM Review, {\bf 8}:119--130 (1961).

\vfil \eject \beginsection {Appendix 1. Proofs}

This appendix proves the theorems cited in the text.

\Proclaim Theorem~1. If $A$ and $A^S$ are nonsingular and
$$\vcenter {\offinterlineskip \halign {\Hbox {\hfil #\hfil }& \vstrut
\hfil #\hfil& \Hbox {\hfil #\hfil }\cr for some matrix $Y$& \Vstrut
or& for some matrix $X$\cr $\X := A^{-1} Y A^S$& \vrule& $\X :=
{}$ any left inverse of $X$\cr $Y^{-} := {}$ any right inverse of $Y$& 
\vrule& $Y^{-} := A^S X A^{-1}$\cr }}$$ then $A^{-1} = \X (A^S)^{-1}
Y^{-}$.

\Proclaim Corollary to Theorem~1. If in addition $$\vcenter
{\offinterlineskip \halign {\Hbox {\hfil #\hfil }& \vstrut \hfil #\hfil 
& \Hbox {\hfil #\hfil }\cr $X := (A^S)^{-1} Y^{-} A$& \vrule& $Y := A
\, \X (A^S)^{-1}$\cr }}$$ then $\X X = I$, $Y Y^{-} = I$ and $A = Y A^S
X$.

{\sl Proof by substitution and multiplication.} In the case
parameterized by $Y$ $$\eqalign {A \, \big[ \X (A^S)^{-1} Y^{-}
\big]& = A \, \big[ \, ( \, A^{-1} Y A^S \, ) \, (A^S)^{-1} Y^{-} \big] =
Y Y^{-} = I\cr \noalign {\smallskip} \X X& = (A^{-1} Y A^S) \big[
(A^S)^{-1} Y^{-} A \big] = I\cr \noalign {\smallskip} Y A^S X& = Y
A^S \big[ \, (A^S)^{-1} Y^{-} A \, \big] = A\cr }$$ and similarly in the
case parameterized by $X$. {\sl End of proof.}

\Proclaim Theorem~2. If $A \rightarrow A^S$ is a row or column
stretching and if $A$ is nonsingular, then $A^S$ is nonsingular and
$A^{-1} = \X (A^S)^{-1} Y^{-}$.

{\sl Proof.} In the row stretching case there are matrices $B$, $G$
of full rank, $P$ a permutation matrix, $Y$ and $Y^{-}$ with the
following relationships. $$A^S = \left[ \matrix {B& G} \right] P^t
\qquad \vcenter {\openup 1\jot \halign {\hfil $#$& $#$\hfil& \qquad
\hfil $#$\hfil& $#$\hfil \cr Y B& {} = A& \X& {} = \left[ \matrix {I& 
0\cr } \right] P^t\cr Y G& {} = 0& Y^{-}& {} = \hbox {any right
inverse of $Y$}\cr }}$$ Suppose $A^S u = 0$, and let $$P \left[
\matrix {v_1 \cr v_2} \right] = u$$ in which the order of $v_1$ is
the column order of $B$, and the order of $v_2$ is the column
order of $G$. With this notation, $B v_1 + G v_2 = A^S u = 0$ so
$$A v_1 = (Y B) v_1 = Y B v_1 - Y \left( {B v_1 + G v_2} \right) = -
Y G v_2 = 0$$ and thus $v_1 = 0$ because $A$ is nonsingular. This
and $B v_1 + G v_2 = 0$ imply $v_2 = 0$ because $G$ has full rank.
Altogether $u = 0$ hence $A^S$ is nonsingular. Moreover, $$A^{-1}
Y A^S = A^{-1} \left( Y A^S P \right) P^t = A^{-1} \left[ \matrix
{A& 0\cr } \right] P^t = \left[ \matrix {I& 0\cr } \right] P^t = \X$$
from which Theorem~1 asserts $A^{-1} = \X (A^S)^{-1} Y^{-}$. The
column stretching case is similar. {\sl End of Proof.}

\Proclaim Corollary to Theorem~2. If $A \rightarrow A^S$ is a row
or column stretching of a nonsingular matrix $A$, if\/ $\X$ and $X$
are the auxiliary matrices in Definition~1, and if $A^S z = y^S$
are the stretched equations corresponding to $A x = y$, then not
only $\X z = x$ but also $z = Xx$.

{\sl Proof.} In the row stretched case, $X x = [(A^S)^{-1} Y^{-} A] x 
= (A^S)^{-1} y^S = z$. The column case is more interesting. The
Theorem says the stretched equations can be used to solve the
unstretched equations by way of the identity $\X z = x$ in which
$\X$ can be any left inverse for the $X$ which parameterizes the
stretching. Since $\X (X x) = x$, so $z - Xx$ lies in the right null
space of every left inverse for $X$. Those left inverses are
precisely $(X^t X)^{-1} X^t + N$ in which $N$ is any left annihilator
of $X$ with the proper row dimension. Let $z - Xx = X v_1 + v_2$
decompose $z - Xx$ over the column space of $X$ and its
orthogonal complement. $$0 = [(X^t X)^{-1} X^t + N] \, (X v_1 +
v_2) = v_1 + N v_2$$ The choice $N = 0$ shows $v_1 = 0$, thus $0
= N v_2$ for every $N$. The rows of $N$ can be any vectors in the
left null space of $X$, which contains ${v_2}^t$. Thus $v_2 = 0$ and
altogether $z - Xx = 0$. {\sl End of Proof.}

\Proclaim Theorem~3. A simple row stretching in the sense of
Definition~2 is a row stretching in the sense of Definition~1, and
similarly for column stretchings. 

{\sl Proof for row stretching.} A simple row stretching has 
$$\vcenter {\openup2\jot \halign {$#$&  \hfil $\displaystyle
#$\hfil& #\cr Q_1 A Q_2 = {}& \left[ \matrix {A_{1,1}& A_{1,2}&
A_{1,3}& \ldots& A_{1,m}\cr A_{2,1}& A_{2,2}&  A_{2,3}& \ldots&
A_{2,m}\cr } \right]\cr & \Big\downarrow \cr Q_3 A^S Q_4 = {}&
\multispan2$\displaystyle \left[ \matrix {A_{1,1}&  A_{1,2}&
A_{1,3}& \ldots& A_{1,m}\cr A_{2,1}\cr & A_{2,2}\cr && 
A_{2,3}\cr &&& \Ddots\cr &&&& A_{2,m}\cr } \hquad \matrix {0& 
0& \ldots& 0\cr -D_1\cr +D_1& -D_2\cr & +D_2& \Ddots\cr && 
\Ddots& -D_{m-1}\cr &&& +D_{m-1}\cr } \right] $\hfil \cr & 
\multispan2$ \displaystyle \hquad \, \underbrace {\hphantom
{\matrix {A_{1,1}& A_{1,2}& A_{1,3}& \ldots& A_{1,m}\cr 
A_{2,1}\cr }}}_{\displaystyle \widetilde B} \hquad \underbrace
{\hphantom {\matrix {-D_1& -D_2& \Ddots& -D_{m-1}\cr
}}}_{\displaystyle \widetilde G} $\hfil \cr }}$$ where the $Q_i$ are
permutation matrices and $\widetilde Y Q_3 A^S Q_4 = \left[
\matrix {Q_1 A Q_2& 0\cr } \right]$ with $$\widetilde Y = \left[
\matrix {I_1\cr & I_2& I_2& \ldots& I_2\cr } \right]$$ in which
$I_1$ and $I_2$ are identity matrices. The permutation matrices
perform the reorderings that Definition~2 allows before and after
the stretching. The matrices $$B = {Q_3}^t \widetilde B {Q_2}^t
\qquad G = {Q_3}^t \widetilde G \qquad P = Q_4  \left[ \matrix
{{Q_2}^t\cr & I\cr } \right] \qquad Y = {Q_1}^{t\,} \widetilde Y
Q_3$$ are needed to give the simple row stretching the
appearance of a general row stretching. $G$ has full rank because
the staircase sparsity pattern of the glue columns makes them
linearly independent (the $D_j$ are non\-sin\-gu\-lar by
definition). $P$ is a permutation matrix. $Y$ has full rank.
Identities like those in Definition~1 follow by marshalling the
permutations. $$\eqalign {A^S P &= \left( {Q_3}^t \left[ \matrix
{\widetilde B& \widetilde G\cr} \right] {Q_4}^t \right) \left( Q_4
\left[ \matrix {{Q_2}^t\cr & I\cr } \right] \right) = \left[ \matrix
{{Q_3}^t \widetilde B {Q_2}^t& {Q_3}^t \widetilde G\cr } \right] =
\left[ \matrix {B& G\cr} \right]\cr Y B &= \left( {Q_1}^{t\,}
\widetilde Y Q_3 \right) \left( {Q_3}^t \widetilde B {Q_2}^t \right)
= {Q_1}^{t\,} \widetilde Y \widetilde B {Q_2}^t = A\cr Y G &= \left(
{Q_1}^{t\,} \widetilde Y Q_3 \right)  \left( {Q_3}^t \widetilde G
\right) = {Q_1}^{t\,} \widetilde Y \widetilde G = 0\cr }$$ {\sl End of
proof.}

\Proclaim Theorem~4. If $A \rightarrow A^S$ is a simple row or
column stretching as in Definition~2, then $$\det A^S = \det A
\prod_{j =1}^{m-1} \det D_j$$ with perhaps a sign change when the
rows and columns are reordered as the definition allows.

{\sl Proof.} With no reordering, simple row stretching has $$A =
\left[ \matrix{A_1\cr  A_2\cr } \right] = \left[ \matrix {A_{1,1}&
A_{1,2}& A_{1,3}&  \ldots& A_{1,m}\cr A_{2,1}& A_{2,2}& A_{2,3}&
\ldots&  A_{2,m}\cr } \right]$$ and $$A^S = \left[ \matrix
{A_{1,1}& A_{1,2}&  A_{1,3}& \ldots& A_{1,m}& 0& 0& \ldots& 0\cr
A_{2,1}&&&&&  -D_1\cr & A_{2,2}&&&& +D_1& -D_2\cr &&
A_{2,3}&&&& +D_2&  \Ddots\cr &&& \Ddots&&&& \Ddots&
-D_{m-1}\cr &&&&  A_{2,m}&&&& +D_{m-1}\cr } \right] \; .$$ Let
$$S = \left[ \matrix{I_1\cr & I_2& I_2& I_2& \ldots& I_2\cr &&
I_2\cr &&&  I_2\cr &&&& \Ddots\cr &&&&& I_2\cr } \right]$$ in
which $I_1$ and $I_2$ are identity matrices with orders equal to
the row orders of $A_1$ and $A_2$. The product $$S A^S = \left[
\matrix {A_{1,1}&  A_{1,2}& A_{1,3}& \ldots& A_{1,m}& 0& 0&
\ldots& 0\cr A_{2,1}&  A_{2,2}& A_{2,3}& \ldots& A_{2,m}& 0& 0&
\ldots& 0\cr &  A_{2,2}&&&& +D_1& -D_2\cr && A_{2,3}&&&&
+D_2& \Ddots\cr  &&& \Ddots&&&& \Ddots& -D_{m-1}\cr &&&&
A_{2,m}&&&&  +D_{m-1}\cr } \right]$$ is $2 \times 2$ block lower
triangular. The first diagonal block is $A$ and the second is $(m-1)
\times (m-1)$ block upper triangular with diagonal blocks $D_1$,
$D_2$, $\ldots \,$, $D_{m-1}$. Thus $$\det A^S = (\det S) (\det A^S)
= \det (S A^S) = \det A \prod_{j =1}^{m-1} \det D_j \; .$$ {\sl End
of proof.}

\Proclaim Lemma~1 to Theorem~5. A sequence of row stretchings
is a row stretching, and similarly for column stretchings. 

{\sl Proof for row stretching by induction on the sequence length.}
Each of two row stretchings $A \rightarrow A^S \rightarrow
A^{SS}$ has matrices $B_i$, $G_i$ of full column rank, $P_i$ a
permutation matrix, and $Y_i$ of full row rank with the following
relationships. $$\vcenter {\openup 2\jot \halign {& \hfil
$\displaystyle #$\tabskip = 0pt& ${}\displaystyle #$\tabskip =
2em\hfil \cr A^S& = \left[ \matrix {B_1& G_1} \right] {P_1}^t& Y_1
B_1& = A& Y_1 G_1& = 0\cr  A^{SS}&= \left[ \matrix {B_2& G_2}
\right] {P_2}^t& Y_2 B_2& = A^S& Y_2 G_2& = 0\cr }}$$ Thus $B_2
= \left[ \matrix {{B_1}^S& {G_1}^S\cr} \right] {P_1}^t$ in which $Y_2
{B_1}^S = B_1$ and $Y_2 {G_1}^S = G_1$. $$A^{SS} = [ \; \overbrace
{{B_1}^S}^{\displaystyle B} \quad \overbrace { {G_1}^S  \quad
G_2}^{\displaystyle G} \; ] \overbrace {\left[ \matrix {{P_1}^t\cr &
I\cr} \right] P_2}^{\displaystyle P} \qquad \overbrace {Y_1 \,
Y_2}^{\displaystyle Y} B = A \qquad Y G = 0$$ $P$ is a permutation
matrix and $Y$ has full row rank. Suppose $G u = 0$. With the rows
of $u$ blocked to match the columns of $G$, then $$0 = Y_2 G u =
Y_2 \left[ \matrix{ {G_1}^S& G_2\cr} \right] \left[ \matrix {u_1\cr
u_2\cr} \right] = G_1 u_1$$ so $u_1 = 0$ because $G_1$ has full
rank, leaving $G_2 u_2 = 0$, so $u_2 = 0$ because $G_2$ has full
rank. Altogether $u = 0$, therefore $G$ has full rank. {\sl End of
proof.}

\Proclaim Lemma~2 to Theorem~5. A sequence of simple row
stretchings produces a row stretching $A \rightarrow A^S$  $$A^S
P = \left[ \matrix {B&  G\cr } \right] \qquad Y B = A \qquad Y G =
0$$ in which $B$, $G$, $P$ and $Y$ are as in Definition~1, and
additionally (1) the entries of $A$ scatter into $B$ in a way that
preserves columns and segregates entries from different rows,
and (2) $\,Y$ is a matrix of\/ $0$'s and\/ $1$'s with exactly one\/
$1$ per column. If the matrices $D_j$ in the simple stretchings are
diagonal, then (3) the columns of $G$ have exactly two nonzeroes
and those nonzeroes have equal magnitudes and opposite signs,
(4) non-loop edges in\/ $\row (G)$ have weight\/~$1$, (5) the
maximally connected subgraphs of\/ $\row (G)$ are trees, and (6)
the nonzeroes in each row of $\,Y$ pick out one maximally
connected subgraph of\/ $\row (G)$. For column stretchings,
replace $Y$ by $X$, replace the equations by $$P A^S = \left[
\matrix {B\cr G\cr } \right] \qquad B X = A \qquad G X = 0$$ and
exchange\/ {\rm row} and\/ {\rm column} in the text.

{\sl Proof for row stretching.} Some parts of this omnibus Lemma
mightn't need proof, but the whole is more easily argued together.
Simple row stretchings are row stretchings (Theorem~3), and
sequences of row stretchings are row stretchings (Lemma~1), so
the matrices $B$, $G$, $P$ and $Y$ exist as required by
Definition~1. 

(1) A simple row stretching $A \rightarrow A^S$ copies entries of
$A$ to $A^S$. Entries from different columns (or rows) go to
different columns (respectively, rows), and those from the same
column go to the same column. The stretching therefore
preserves columns but only segregates rows. Moreover, it places 
glue in separate new columns. An iterated stretching $$A = A_{0}
\rightarrow {A_{0}}^S = A_1 \rightarrow {A_1}^S = A_2
\rightarrow {A_2}^S = A_3 \; \cdots \; A_q = A^S$$ merely copies
the entries of $A$ and the accumulating glue several times. 

(2) For a simple row stretching, $Y$ is a matrix of $0$'s and $1$'s
with exactly one $1$ per column (see the proof of Theorem~3).
This distinctive sparsity pattern is inherited by the product of
such matrices. For a sequence of row stretchings, the overall $Y$
is the product of the $Y_i$ for each stretching (see the proof of
Lemma~1).

(3) When the $D_j$ are diagonal, a simple row stretching $A
\rightarrow A^S$ places exactly two pieces of glue with identical
magnitudes and opposite signs in new columns of $A^S$.
Subsequent stretchings preserve columns, and while they may
rearrange old columns of glue, they add nothing to them. 

(4) Non-loop edges in $\row (G)$ have weight $1$ because each
glue column has exactly two nonzeroes, and no two columns have
the same sparsity pattern, so each column is the unique edge
between two rows. If two columns were to have the same sparsity
pattern, then they'd be linearly dependent, and $G$ couldn't have
full rank.

(5) Discarding reorderings, the glue columns of a simple row
stretching $A \rightarrow A^S$ are those new columns containing
the $D_j$. $$\vcenter {\openup2\jot \halign {$\displaystyle #$\hfil&
#\cr \left[ \matrix {A_{1,1}& A_{1,2}& \ldots& A_{1,m}\cr
A_{2,1}& A_{2,2}& \ldots& A_{2,m}\cr } \right]\cr \hfil
\Big\downarrow \cr \multispan2{$\displaystyle \left[ \matrix {
A_{1,1}& A_{1,2}& \ldots& A_{1,m}\cr A_{2,1}\cr & A_{2,2}\cr &&
\Ddots\cr &&&  A_{2,m}\cr } \hquad \matrix {0& \ldots& 0\cr
-D_1\cr +D_1& \Ddots\cr & \Ddots& -D_{m-1}\cr && +D_{m-1}\cr
} \right]$ \hfil} \cr }}$$ When the $D_j$ are diagonal the glue is
better viewed after shuffling rows and columns to group entries
from the same diagonal position $$\left[ \matrix {0\cr \phantom
{T_2} & T_1\cr&& T_2\cr &&& \Ddots \cr &&&& T_k} \right] \quad
\hbox {in which} \hquad T_i = \left[ \matrix {-d_{i,1}\cr +d_{i,1}&
-d_{i,2}\cr & +d_{i,2}& \Ddots\cr && \Ddots& -d_{i,m-1}\cr &&&
+d_{i,m-1}\cr } \right]$$ where $k$ is the order of the $D_j$ and
$d_{i,j}$ is the $i$-th diagonal entry of $D_j$. The row graph of each
$T_i$ is a tree (a linear tree with two leaves and no branches). No
two $T_i$ overlap in columns of $G$, so the rows containing each
$T_i$ are a maximally connected subgraph of $\row (G)$. Any other
row of $G$ is zero, hence connected to nothing, hence a maximally
connected subgraph and a trivial tree. 

Moreover, each row of $A$ that stretches has its own $T_i$ in $G$,
and each row of $A$ that doesn't stretch has its own zero row in
$G$, so altogether, the rows of $A$ number the same as the
maximally connected subgraphs in $\row (G)$. 

For a sequence of simple row stretchings, all but the last can be
collapsed to one stretching, leaving $A \rightarrow A^S
\rightarrow A^{SS}$ in which, by induction hypothesis, the
maximally connected subgraphs in the row graph of the glue
columns of $A^S$ are trees. With a suitable ordering, the last
stretching $A^S \rightarrow A^{SS}$ changes the glue columns as
follows. $$\vcenter {\openup2\jot \halign {$\displaystyle #$\hfil&
#\cr \left[ \matrix {G_{1,1}& G_{1,2}& \ldots& G_{1,m}\cr G_{2,1}&
G_{2,2}& \ldots& G_{2,m}\cr } \right]\cr \hfil \Big\downarrow \cr
\multispan2{$\displaystyle \left[ \matrix { G_{1,1}& G_{1,2}&
\ldots& G_{1,m}\cr G_{2,1}\cr & G_{2,2}\cr && \Ddots\cr &&& 
G_{2,m}\cr } \hquad \matrix {0& \ldots& 0\cr -D_1\cr +D_1&
\Ddots\cr & \Ddots& -D_{m-1}\cr && +D_{m-1}\cr } \right]$ \hfil}
\cr }}$$ Some of the glue blocks $G_{i,j}$ may have zero column
dimension because the stretching $A^S \rightarrow A^{SS}$
needn't copy old glue to every newly stretched row. For example,
if no row of $A^S$ containing old glue stretches, then the glue
columns of $A^{SS}$ would be $$\left[ \matrix { G_{1,1}& 0&
\ldots& 0\cr & -D_1\cr & +D_1& \Ddots\cr && \Ddots&
-D_{m-1}\cr &&& +D_{m-1}\cr } \right]$$ but there is no harm in
allowing null blocks into the original  display. As before, the glue
is better viewed after shuffling rows and columns $$\left[ \matrix
{B_0\cr B_1& T_1\cr B_2&& T_2\cr \Vdots&&& \Ddots \cr B_k
&&&& T_k} \right] \quad \hbox {in which} \left\{ \eqalign {B_0& =
\; \left[ \matrix {G_{1,1}& G_{1,2}& \ldots& G_{1,m}\cr } \,
\right]\cr B_i& = \left[ \matrix {g_{i,1}& \hphantom {G_{1,2}}&
\hphantom {\ldots}& \hphantom {G_{1,m}}\cr \hphantom {G_{1,1}}&
g_{i,2}\cr && \Ddots\cr &&& g_{i,m}\cr } \right]} \right.$$ where
now $g_{i,j}$ is the $i$-th row of $G_{2,j}$, and the $T_i$ are as
pictured earlier. When a glue-bearing row of $A^S$ stretches
$$\left[ \matrix {g_{i,1}& g_{i,2}& \ldots& g_{i,m}\cr } \, \right] \;
\longrightarrow \; B_i = \left[ \matrix {g_{i,1}& \hphantom
{G_{1,2}}& \hphantom {\ldots}& \hphantom {G_{1,m}}\cr \hphantom
{G_{1,1}}& g_{i,2}\cr && \Ddots\cr &&& g_{i,m}\cr } \right]$$ the
effect on the row graph is to replace one vertex by several,
dividing the former's edges among the latter, thus breaking one
maximally connected subgraph into smaller trees (the original
maximally connected subgraphs are trees by the induction
hypothesis, and branches of trees are trees), and finally grafting
all the branches onto the new linear tree created by $T_i$.  Thus,
the old trees merely rearrange and grow to form new trees---of
the same number. That the rows of $A$ number the same as the
maximally connected subgraphs of $\row (G)$ is needed below.

(6) In the multiplication $Y G = 0$, only the two rows that contain a
glue column's two nonzeroes can be summed to annihilate that
column's glue. The multiplication therefore sums entire maximally
connected subgraphs in $\row (G)$, that is, the nonzeroes in each
row of $Y$ pick out whole maximally connected subgraphs in $\row
(G)$. Each row of $Y$ must pick out at least one maximally
connected subgraph because $Y$ has no zero rows ($Y$ has full row
rank). Each maximally connected subgraph must be chosen by
some row of $Y$ because $Y$ has no zero columns (columns of $Y$
have one nonzero apiece). The rows of $Y$ number the same as the
rows of $A$, and as explained in the proof of (5) above, the rows
of $A$ number the same as the maximally connected subgraphs in
$\row (G)$. Therefore, each row of $Y$  picks out exactly one
maximally connected subgraph. {\sl End of proof.}

\Proclaim Lemma~3 to Theorem~5. If the row graph of a matrix is
a tree whose non-loop edges have weight $1$, and if each column
has exactly $2$ nonzeroes, then the removal of any row leaves a
nonsingular matrix. If the nonzeroes in the matrix are $\pm 1$,
then the nonzeroes in the inverse are $\pm 1$. The column
dimension therefore bounds the $1$-norm and $\infty$-norm of
the inverse.

{\sl Proof by induction on the number of columns.} If there is one
column, then the column's two nonzeroes must connect all rows
(since the row graph is a tree), so there are just two rows.
Removing any row leaves a non\-singular matrix and so on.

The inductive step uses the following observations. (1) Discarding
any leaf-row and its column-edges makes a smaller matrix that
satisfies the Lemma's hypotheses (the graph remains a tree
because only a leaf is lost; each remaining column has exactly two
nonzeroes because columns that lost nonzeroes with the removed
row are removed too). (2) A leaf-row has exactly one nonzero
(columns have two nonzeroes so the nonzeroes in a row must
connect to other rows; edges have weight one so the nonzeroes in
a row connect to separate neighbors; a leaf has just one neighbor
hence no more than one nonzero). These facts supports many
inductive proofs, for example (3) the matrices described by the
Lemma are square but for one extra row (if the matrix has one
column then it has been observed to have two rows; if it has more
than one column then removing a leaf row and its single column
edge leaves a smaller matrix like the original).

The proof's inductive step has two cases. First, the removed row
neighbors every other. In this case the columns of the surviving
matrix have one nonzero apiece, and these lie in distinct rows
because all non-loop edges in the original row graph have weight
$1$. This means the new, square matrix has the sparsity pattern of
a permutation matrix. Its inverse is formed by transposing and
reciprocating.

Second, the removed row doesn't neighbor every other. Some path
from the row therefore ends at a non-neighboring leaf. With the
row to be removed placed last and with the non-neighboring
leaf-row and the leaf's single column-edge placed first, the matrix
has the form $$\left[ \matrix {\sigma& 0\cr c& B\cr 0& r^t\cr }
\right]$$ in which the number $\sigma$ is not zero and the column
vector $c$ has one nonzero entry. The induction hypothesis
applied to $$\left [ \matrix {B\cr r^t} \right]$$ makes $B$
nonsingular so $$\left[ \matrix {\sigma& 0\cr c& B\cr } \right]$$ is
nonsingular too. If the nonzeroes are $\pm 1$ then $$\left[ \matrix
{\pm 1& 0\cr c& B\cr } \right]^{-1} = \left[ \matrix {\pm 1\cr \mp
B^{-1} c& B^{-1}\cr } \right]$$ in which $\mp B^{-1} c$ is ${\pm}$ a
column of $B^{-1}$. Again by the induction hypothesis, the
nonzeroes in the inverse must be $\pm 1$. {\sl End of proof.}

\Proclaim Lemma~4 to Theorem~5. If\/ $A \rightarrow A^S$ is a
row stretching produced by a sequence of simple row stretchings,
and if $$A^S P = \left[ \matrix {B&  G\cr } \right] \qquad Y B = A
\qquad Y G = 0$$ are as in Definition~1, and if $m$ is the size of
the largest maximally connected subgraph in\/ $\row (G)$, and if\/
$Q$ is a permutation matrix that places any member of the glue's
maximally connected subgraph for row $j$ of $A$ into row $j$ of
$A^S$ (resulting in the following blockings) $$Q B = \left[ \matrix
{B_1\cr B_2\cr} \right] \qquad Q G = \left[ \matrix {G_1\cr
G_2\cr} \right] \qquad Y Q^t = \left[ \matrix {I& Y_2\cr} \right]$$
then $G_2$ is invertible. If $A$ is nonsingular, then there is an
explicit representation for $(A^S)^{-1}$. $$P^t (A^S)^{-1} Q^t =
\left[ \matrix {A^{-1}& A^{-1} Y_2\cr - {G_2}^{-1} B_2 A^{-1}&
{G_2}^{-1} (I - B_2 A^{-1} Y_2)\cr } \right]$$ Moreover, if all the
nonzero entries of\/ $G$ are\/ $\pm \sigma$ for the same
$\sigma$, then $\| {G_2}^{-1} B_2 \|$ has the following bounds. $$\|
{G_2}^{-1} B_2 \|_1 \le (m-1) \| A \|_1 / | \sigma | \quad  \hbox {and}
\quad \| {G_2}^{-1} B_2 \|_\infty \le \| A \|_\infty / | \sigma | $$

{\sl Proof.} $Y$ is a matrix of $0$'s and $1$'s with exactly one $1$
per column  (part~2 of Lemma~2), and the nonzeroes in each row
of $Y$ pick out a maximally connected subgraph of $\row (G)$
(part~6), so the chosen row ordering places an identity matrix in
the first block of $Y Q^t$ as shown.

The columns of $G$ have exactly two nonzeroes (part~3 of
Lemma~2), non-loop edges in $\row (G)$ have weight $1$ (part~4),
and the maximally connected subgraphs of\/ $\row (G)$ are trees
(part~5). Each maximally connected component of $\row (G)$
therefore satisfies Lemma~3's hypotheses. Any reordering of $G$
that groups the maximally connected subgraphs of $\row (G)$ and
their column edges makes $G$ block diagonal. $G_2$ is obtained
by removing one row from each block, so $G_2$ is nonsingular
(Lemma~3), and in the reordering described, $G_2$ is block
diagonal with blocks no larger than $m-1$. 

If the nonzero entries of $G$ are $\pm \sigma$, then those of
${G_2}^{-1}$ are $\pm 1 / \sigma$ (Lemma~3), and in view of the
blocking described above, no more than $m-1$ nonzeroes
populate any row or column of $G_2$. The $1$-norm bound follows
from this and the fact each column of $B_2$ has entries copied
from only one column of $A$ (part~1 of Lemma~4). $$\| {G_2}^{-1}
B_2 \|_1 \le \| {G_2}^{-1} \|_1 \, \| B_2 \|_1 \le | 1 / \sigma | (m-1) \|
A \|_1$$ The rows of $B$ segregate entries from different rows of
$A$ (part~1 of Lemma~4), so the product ${G_2}^{-1} B_2$ only
sums rows stretched from the same row of $A$. Thus, each row of
the product partly reassembles some row of $A$, perhaps with
different signs, and has $\infty$-norm bounded by $\| A \|_{\infty}
/ | \sigma |$. 

Finally, the identities $$\eqalign {A &= Y B = (Y Q^t) (Q B) = B_1 +
Y_2 B_2 \cr 0 &= Y G = (Y Q^t) (Q G) = G_1 + Y_2 G_2\cr}$$ allow
$Q A^S P$ to be written succinctly $$Q A P = \left[ \matrix {A - Y_2
B_2& - Y_2 G_2\cr B_2& G_2\cr } \right]$$ so the formula for the
inverse can be verified by multiplication. {\sl End of proof.}

\Proclaim Theorem~5. If $$A \rightarrow A^S \rightarrow A^{SS}
\cdots \rightarrow \A$$ is a sequence of simple row or column
stretchings\/ {\rm but not both}, and if each row or column of $A$
stretches to at most $m$ rows or columns of $A^{SS \cdots S}$,
and if Definition~2's matrices $D_i$ have the form $\sigma I$ for
the same $\sigma$, then the following choices for $\sigma$
$$\vbox {\offinterlineskip \halign {\hfil $\;$#& \quad\vrule\quad#& 
\hfil $#$\hfil \quad& \hfil $#\;$\hfil \cr $\sigma$\vrule depth6pt
height9pt width0pt&& p =1& p =\infty\cr \noalign {\hrule}
row\vrule depth3pt height15pt width0pt&& \| A \|_p / 2& \| A
\|_p\cr column\vrule depth3pt height15pt width0pt&& \| A \|_p& 
\| A \|_p / 2\cr }}$$ yield a final stretched matrix $\A$ with
bounded condition number $$\kappa_p (\A) \le c \, \kappa_p (A)$$
in which the multiplier $c$ is given below. $$\vbox
{\offinterlineskip \halign {\hfil $\;$#& \quad\vrule\quad#& 
\hfil #\hfil \quad& \hfil #$\;$\hfil \cr $c$\vrule depth6pt height9pt
width0pt&& $p =1$& $p =\infty$\cr \noalign {\hrule} row\vrule
depth3pt height15pt width0pt&& $2m-1$& $m^2$\cr column\vrule
depth3pt height15pt width0pt&& $m^2$& $2m-1$\cr }}$$ When the
sequence of stretchings is disjoint in the sense that later
stretchings do not alter the rows or columns of earlier
stretchings, then $3m$ can replace $m^2$ in this table. All these
bounds are sharp for some matrices.

{\sl Proof.} Only row stretchings need be considered because
column stretchings are the transpose. Moreover, multiple
stretchings can be treated simultaneously because the Lemmas
have done the dirtiest work. A straightforward proof remains. It
blocks the columns of $\A$ to bound $\| \A \|$, it blocks the rows
of $\A$ to bound $\| ({\A})^{-1} \|$, and then it minimizes the
product of the bounds.

Order the columns of $\A$ as in Definition~1 so the entries of $A$
lie in their original columns and $\pm \sigma$ lie in the others.
$$\A P = \left[ \matrix {B& \sigma G\cr } \right]$$ The nonzeroes
in each column of $B$ are exactly the nonzeroes in the same
column of $A$. The rows have been stretched however, so the
nonzeroes in each row only lie among the nonzeroes in some row
of $A$. Thus $$\| B \|_1 = \| A \|_1 \quad \hbox {and} \quad \| B
\|_{\infty} \le \| A \|_{\infty}.$$ Simple row stretchings build the
glue columns, $\sigma G$, by placing one $\pm \sigma$ pair per
column. If the stretchings in the sequence do not alter rows
stretched earlier, then each row of $G$ acquires at most $2$
nonzero entries. Otherwise, the nonzeroes in a row of $G$ could
number as large as $m-1$. Thus $$\| G \|_1 = 2 \quad \hbox {and}
\quad \| G \|_{\infty} \le \hbox {$2$ or $m-1$}$$ and the following
bounds on $\| \A \|$ have been easily obtained. $$\| {\A} \|_1 \le
\max \, \{ \| A \|_1 ,\, 2 \sigma \} \qquad \| {\A} \|_{\infty} \le \| A
\|_{\infty} + (\hbox {$2$ or $m-1$}) \sigma$$ The bounds on $\| ( \A
)^{-1} \|$ are more subtle.

Retain the column blocking, above, and additionally order the
rows as in Lemma~4 to place any member of the glue's maximally
connected subgraph for row $j$ of $A$ into row $j$ of $\A$. Place
the other stretched rows in a second block. $$\openup 2\jot
\displaylines {Q \A P = \left[ \matrix {B_1& \sigma G_1\cr B_2&
\sigma G_2\cr } \right]\cr Y Q^t = \left[ \matrix {I& Y_2\cr} \right]
\quad \quad Y \A P = \left[ \matrix {A& 0\cr} \right]\cr}$$ $Y$ is
the matrix of Definition~1 and Lemma~2 that reverses the
stretching. $Y_2$ contains only $0$'s and $1$'s, and row $j$ of $Y_2$
picks out the other members of the subgraph for row $j$ of $A$,
hence $\| Y_2 \|_\infty \le m - 1$. Lemma~4 proves the following.
$$\openup 3\jot \displaylines {P^t (\A)^{-1} Q^t = \left[ \matrix
{A^{-1}& A^{-1} Y_2\cr - \sigma^{-1} {G_2}^{-1} B_2 A^{-1}&
\sigma^{-1} {G_2}^{-1} (I - B_2 A^{-1} Y_2)\cr } \right]\cr
\displaystyle \| {G_2}^{-1} B_2 \|_1 \le (m-1) \| A \|_1 \quad  \hbox
{and} \quad \| {G_2}^{-1} B_2 \|_\infty \le \| A \|_\infty\cr }$$ The
present notation differs from the Lemma's because $\sigma$ is
implicit there ($G_2$ in Lemma~4 is $\sigma G_2$ here).

At this point the proof divides into separate cases for each norm.
The $1$-norm of the first block-column of $(\A)^{-1}$ is bounded
by $$\eqalign {\| A^{-1} \|_1 + {}& \| \sigma^{-1} {G_2}^{-1} B_2
A^{-1} \|_1\cr \noalign {\medskip} {} \le \; \| A^{-1} \|_1 + {}& 
\sigma^{-1} (m-1) \| A \|_1 \| A^{-1} \|_1.\cr }$$ The more
complicated second block-column doesn't need separate
attention because the row ordering for $\A$ can place any row in
the first block of $\A$, and thus can place any column in the first
block of $(\A)^{-1}$. The bound above therefore applies to all
columns of $(\A)^{-1}$ and so to $\| (\A)^{-1} \|_1$. 

The product of the bounds on $\| \A \|_1$ and $\| (\A)^{-1} \|_1$ is a
maximum of two functions, one increasing with $\sigma$ and the
other decreasing. $$\eqalign {& \| {\A} \|_1 \| (\A)^{-1} \|_1\cr 
\noalign {\medskip}& \quad \le \max \, \{ \| A \|_1 ,\, 2 \sigma \}
\times \Big[ \, \| A^{-1} \|_1 + \sigma^{-1} (m-1) \| A \|_1 \| A^{-1}
\|_1\Big]\cr \noalign {\medskip}& \quad = \max \; \left\{ \vcenter
{\openup1\jot \halign {$#$\hfil \cr \| A \|_1 \| A^{-1} \|_1 +
\sigma^{-1} (m-1) {\| A \|_1}^2 \| A^{-1} \|_1\cr 2 \sigma \| A^{-1}
\|_1 + 2 (m-1) \| A \|_1 \| A^{-1} \|_1\cr }} \right. \cr }$$ The minimax
with respect to $\sigma$ occurs where the two functions match.
The bound on $\| \A \|_1$ indicates this point is $\sigma = \| A \|_1
/ 2$ where the bound on $\kappa_1 (\A)$ is $(2m - 1) \kappa_1
(A)$. This bound is sharp for the following simple stretching of the
$1 \times 1$ identity matrix to an $m \times m$ matrix. $$[1]^S =
\left[ \matrix {1& -0.5\cr & +0.5& \Ddots\cr && \Ddots& -0.5\cr 
&&& +0.5\cr } \right] \quad \hbox {and} \quad ([1]^S)^{-1} = \left[
\matrix {1& 1& \cdots& 1\cr & 2& \cdots& 2\cr && \Ddots& 
\Vdots\cr &&& 2\cr } \right].$$ 

The $\infty$-norm case requires a closer look at $(\A)^{-1}$. The
norm of the first block-row is easily bounded by $$\| A^{-1}
\|_{\infty} + \| A^{-1} Y_2 \|_{\infty} \le m \| A^{-1} \|_{\infty}$$
because $\| Y_2 \|_{\infty} \le m-1$. The norm of the lower left
block is cleverly bounded by $$\| - \sigma^{-1} {G_2}^{-1} B_2
A^{-1} \|_{\infty} \le \sigma^{-1} \| A \|_{\infty} \| A^{-1}
\|_{\infty}$$ using an inequality supplied by Lemma~4. As in the
$1$-norm case, the columns in the left block of $({\A})^{-1}$
correspond to a specific choice of representatives for connected
components in the row graph of $G$. Each column has at most $m$
choices, so the $\infty$-norm of the entire lower block-row of
$(\A)^{-1}$ can't exceed $m \sigma^{-1} \| A \|_{\infty} \| A^{-1}
\|_{\infty}$. $$\| (\A)^{-1} \|_{\infty} \le \max \, \{ m \| A^{-1}
\|_{\infty}, m \sigma^{-1} \| A \|_{\infty} \| A^{-1} \|_{\infty} \}$$

Repeating the earlier argument, the product of the bounds on $\|
\A \|_{\infty}$ and $\| (\A)^{-1} \|_{\infty}$ is a maximum of two
functions, one increasing with $\sigma$ and the other decreasing.
$$\eqalign {& \| {\A} \|_{\infty} \| (\A)^{-1} \|_{\infty}\cr \noalign
{\medskip}& \quad \le \left( \| A \|_{\infty} + \alpha \sigma \right)
\times \max \, \{ m \| A^{-1} \|_{\infty}, m \sigma^{-1} \| A
\|_{\infty} \| A^{-1} \|_{\infty} \}\cr \noalign {\medskip}& \quad =
\max \left\{ \; \vcenter {\openup1\jot \halign {$#$\hfil \cr m \| A
\|_{\infty} \| A^{-1} \|_{\infty} + \alpha \sigma m \| A^{-1}
\|_{\infty}\cr m \sigma^{-1} \| A \|_{\infty}^{\; 2} \| A^{-1}
\|_{\infty} + \alpha m \| A \|_{\infty} \| A^{-1} \|_{\infty} \cr }} \right.
\cr }$$ If the stretchings in the sequence are disjoint, then the
$\alpha$ in this formula is $2$, but if they operate on
each others' rows, then $\alpha$ is $m-1$. In any case, the minimax
with respect to $\sigma$ occurs where the two functions match.
The bound on $\| (\A)^{-1} \|_{\infty}$ indicates this is $\sigma = \|
A \|_{\infty}$ where the bound on $\kappa_{\infty} (\A)$ is $(\alpha
+ 1) m \kappa_{\infty} (A)$. When $\alpha = m - 1$ this bound is
sharp for the following iterated stretching of the $1 \times 1$
negative identity matrix to an $m \times m$ matrix. $$[-1]^{SS
\cdots S} = \left[ \matrix { -1& -1& -1& \cdots& -1\cr & +1\cr &&
+1\cr &&&  \Ddots\cr &&&& +1\cr } \right] \quad \hbox {and}
\quad ([-1]^{SS \cdots S})^{-1} = [-1]^{SS \cdots S}$$ When $\alpha
= 2$ the bound is sharp for the following simple stretching of the
$1 \times 1$ negative identity matrix to an $m \times m$ matrix.
$$[-1]^S = \left[ \matrix { 0& -1\cr -1& +1& -1\cr && +1&
\Ddots\cr &&&  \Ddots& -1\cr &&&& +1\cr } \right] \quad \hbox
{and} \quad ([-1]^S)^{-1} = \left[ \matrix { 1& 1& 1& \cdots& 1\cr
1& 0& 0&  \cdots& 0\cr && 1& \cdots& 1\cr &&& \Ddots&
\Vdots\cr &&&&  1\cr } \right]$$ {\sl End of proof.}

\Proclaim Lemma~1 to Theorem~6. If $A \rightarrow A^S$ is a row
or column stretching of a nonsingular matrix obtained from a
sequence of simple row or column stretchings\/ {\rm but not
both}, and if the glue is chosen by Theorem~5, and if the stretched
matrix is computed in finite precision arithmetic with unit
roundoff $\epsilon$, then the computed $\AS$ differs from the
ideal $A^S$ as follows $$\AS = A^S + E \qquad \| E \|_p \le c_{5.1} \,
[(1+ \epsilon)^{n}-1] \, \| A \|_p$$ where $n$ is the order of $A$ and
$c_{5.1}$ is given by the table $$\vbox {\offinterlineskip \halign
{\hfil $\;$#& \quad\vrule\quad#& \hfil #\hfil \quad&  \hfil #$\;$\hfil
\cr $c_{5.1}$\vrule depth6pt height9pt width0pt&&  $p =1$& $p
=\infty$\cr \noalign {\hrule} row\vrule depth3pt height15pt
width0pt&& $2$& $m$\cr column\vrule depth3pt height15pt
width0pt&& $m$& $2$\cr }}$$ in which  each row of $A$ stretches
to at most $m$ rows of $A^S$. When the sequence of stretchings is
disjoint in the sense that later stretchings do not alter the rows
or columns of earlier stretchings, then $2$ can replace $m$ in this
table.

{\sl Proof for $p = \infty$.\/} In the column case, Definition~2's
stretchings copy the entries of $A$ to $A^S$ and additionally place
two pieces of Theorem~5's glue, $\pm \| A \|_{\infty} / 2$, in new
rows. The computed and ideal stretched matrices thus differ by a
matrix $E$ whose nonzero entries equal the error in computing
$\pm \| A \|_{\infty} / 2$. If $\delta_c$ bounds this error then $\| E
\|_{\infty} \le 2 \delta_c$.

In the row case, Definition~2's stretchings again copy the entries
of $A$ to $A^S$ and place one or two pieces of Theorem~5's glue,
$\pm \| A \|_{\infty}$, in the stretched rows. Should the stretchings
in the sequence not alter rows stretched earlier, then each row
acquires at most two pieces of glue, otherwise, a row could
acquire as many as $m-1$. If $\delta_r$ bounds the error in
computing $\| A \|_{\infty}$ then $\| E \|_{\infty} \le (\hbox{$2$ or
$m-1$}) \delta_r$. 

The standard model of machine computation gives to each
arithmetic operation a small relative perturbation bounded by
the optimistically named {\sl unit roundoff\/}. The $j$-th absolute
row sum of $A$ may be computed as follows $$\vcenter
{\openup1\jot \halign {\hfil $#$& ${}#$\hfil \qquad& \hfil $#$& 
${}#$\hfil \cr s_{j,1}& = | a_{j,1} |& \bars_{j,1}& = | a_{j,1} |\cr 
s_{j,2}& = s_{j,1} + | a_{j,2} |& \bars_{j,2}& = (\bars_{j,1} + | a_{j,2} |)
(1 + \epsilon_2)\cr & \;\, \vdots&& \;\, \vdots \cr s_{j,n} & =
s_{j,n-1} + | a_{j,n} |& \bars_{j,n}& = (\bars_{j,n-1} + | a_{j,n} |) (1 +
\epsilon_n)\cr }}$$ in which bars denote computed quantities, sign
changes are errorless, the $a_{j,k}$ are entries of $A$, and $|
\epsilon_i | \le \epsilon$. The ideal and computed sums therefore
are $$s_{j,n} = \sum_{k =1}^n |a_{j,k}| \qquad \bars_{j,n} = \sum_{k
=1}^n \left[ |a_{j,k}| \prod_{i =k\wedge2}^n (1 + \epsilon_i)
\right]$$ so the computed sum lies between $(1 \pm
\epsilon)^{n-1}$ of the ideal. The ideal and computed $\| A
\|_{\infty}$ are the largest $s_{i,n}$ and $\bars_{j,n}$, respectively.
No computed sum exceeds $(1 + \epsilon)^{n-1}$ of the largest
ideal sum, so the computed $\| A \|_{\infty}$ lies between $(1 \pm
\epsilon)^{n-1} \| A \|_{\infty}$. Thus $\delta_r \le [(1 \pm
\epsilon)^{n-1}-1] \, \| A \|_{\infty}$ whence $$\| E \|_\infty \le
(\hbox{$2$ or $m-1$}) [(1 \pm \epsilon)^{n-1}-1] \, \| A \|_\infty$$
for row stretching. The weaker bound  $(\hbox{$2$ or $m$}) [(1 \pm
\epsilon)^{n}-1] \, \| A \|_{\infty}$ is more concise.  As for
$\delta_c$, halving the computed $\| A \|_{\infty}$ introduces one
more relative pertubation with the result that $\delta_c \le [(1
\pm \epsilon)^{n}-1] \, \| A \|_\infty / 2$ hence $$\| E \|_\infty \le
[(1 \pm \epsilon)^{n}-1] \, \| A \|_\infty$$ for column stretching.
{\sl End of proof.}

\Proclaim Lemma~2 to Theorem~6. If $A \rightarrow A^S$ is a row
or column stretching of a nonsingular matrix obtained from a
sequence of simple row or column stretchings\/ {\rm but not
both}, and if in the row case the glue is chosen by Theorem~5 (the
column case may choose any glue), and if the vector operations
used to solve linear equations are scatter $y \rightarrow y^S$ and
gather $z \rightarrow z_\S$ operations, and if $\barz$ is an
approximate solution to the stretched equations $A^S z = y^S$,
then $\barx := {\barz }_\S$ can be regarded as an approximate
solution to the unstretched equations $A x = y$. The error in
$\barz$ bounds the error in $\barx$ $${\| x - \barx \|_p \over \| x
\|_p} \le c_{5.2} \, {\| z - \barz \|_p \over \| z \|_p}$$ where $c_{5.2}$
is given by the table $$\vbox {\offinterlineskip \halign {\hfil
$\;$#& \quad\vrule\quad#&  \hfil #\hfil \quad& \hfil #$\;$\hfil \cr
$c_{5.2}$\vrule depth6pt height9pt width0pt&& $p =1$& $p
=\infty$\cr \noalign {\hrule} row\vrule depth3pt height15pt
width0pt&& $2m-1$& $1$\cr  column\vrule depth3pt height15pt
width0pt&& $m$& $1$\cr }}$$ in which $n$ is the order of $A$ and
each row of $A$ stretches to at most $m$ rows of $A^S$. 

{\sl Proof.} When the squeezing $z \rightarrow z_\S := \X z$ is a
gather operation then $\X$ is a matrix of $0$'s and $1$'s with one
$1$ per row and no more than one per column. $$\| x - \barx \| = \|
\X (z - \barz) \| \le \| \X \| \, \| z - \barz \| = \| z - \barz \|$$ This and
$\| z \| \le c_{5.2} \| x \|$ will imply the Lemma's bound on the
relative error.

For a column stretching, $z = Xx$ (by the Corollary to
Theorem~2), in which $X$ is a matrix of $0$'s and $1$'s with exactly
one $1$ per row and from one to $m$ per column (parts~2 and~6
for the column case of Lemma~1 to Theorem~5). Hence $\| X \| =
c_{5.2}$ and $\| z \| \le c_{5.2} \| x \|$.

For a row stretching, $Y^{-}$ can be any right inverse of the matrix
$Y$ that parameterizes the stretching in Definition~1 $$A^S =
\left[ \matrix {B& G\cr } \right] P^t \qquad Y B = A \qquad Y G =
0$$ but when $y \rightarrow y^S = Y^{-} y$ is a scatter operation
then $$Y^{-} = Q \left[ \matrix {I\cr 0\cr } \right]$$ in which $I$ is an
identity matrix and $Q$ is a permutation matrix. Let $$Y = \left[
\matrix {I& Y_2\cr } \right] Q^t \qquad B = Q \left[ \matrix {B_1\cr 
B_2\cr } \right] \qquad G = Q \left[ \matrix {G_1\cr G_2\cr }
\right]$$ match the blocking of $Y^{-}$. Lemma~4 to Theorem~5
shows $$(A^S)^{-1} = P \left[ \matrix {A^{-1}& A^{-1} Y_2\cr -
{G_2}^{-1} B_2 A^{-1}&  {G_2}^{-1} (I - B_2 A^{-1} Y_2)\cr } \right]
Q^t$$ which provides a formula for $z$. $$z = (A^S)^{-1} y^S =
(A^S)^{-1} Y^{-} y = (A^S)^{-1} Y^{-} A x = P \left[ \matrix {I\cr -
{G_2}^{-1} B_2\cr } \right] x$$ With Theorem~5's choice of glue the
Lemma also provides the following bounds. $$\openup 2\jot
\eqalign {\sigma = \| A \|_1 / 2 &\; \Longrightarrow \; \| {G_2}^{-1}
B_2 \|_1 \le \| A \|_1 (m - 1) / \sigma = 2 (m - 1)\cr \sigma = \| A
\|_\infty & \; \Longrightarrow \; \| {G_2}^{-1} B_2 \|_\infty \le \| A
\|_\infty / \sigma = 1\cr}$$ Thus, $c_{5.2}$ is $2m-1$ for the
$1$-norm, and $1$ for the $\infty$-norm. {\sl End of proof.}

\Proclaim Theorem~6. If $A \rightarrow A^S$ is stretching of a
nonsingular matrix obtained from a sequence of simple row or
column stretchings\/ {\rm but not both}, and if glue is chosen by
Theorem~5, and if the stretched matrix $\AS$ is computed in
finite precision arithmetic with unit roundoff $\epsilon$, and if the
vector operations used to solve linear equations are error-free
scatter $y \rightarrow y^S$ and gather $z \rightarrow z_\S$
operations, and if the approximate solution $\barz$ to the
computed stretched equations $\AS z = y^S$ exactly satisfies
some perturbed equations $(\AS + E) \barz = y^S$, then
$$\delta_1 := c_1 \, [(1 + \epsilon)^{n}-1] < 1 \qquad \hbox {and}
\qquad \delta_2 := { \| E \|_p \over \| \AS \|_p} < {1 - \delta_1
\over 1 + \delta_1 \vphantom {\AS}}$$ imply $$\eqalign {{\| x -
\barz_\S \|_p \over \| x \|_p} &< c_2 \, \kappa (A) \, {(\delta_1 +
\delta_2 + \delta_1 \delta_2) \over 1 - (\delta_1 + \delta_2 +
\delta_1 \delta_2)}\cr \noalign {\smallskip} &\approx c_2 \,
\kappa (A) \left( c_1 n \epsilon + \| E \|_p / \|
\AS \|_p \right) \cr}$$ \smallskip \noindent in which $c_1$ and
$c_2$ are given by the tables $$\vcenter {\offinterlineskip \halign
{\hfil $\;$#& \quad \vrule \quad#& \hfil #\hfil \quad& \hfil
#$\;$\hfil \cr  $c_1$\vrule depth6pt height9pt width0pt&& $p =
1$& $p = \infty$\cr \noalign {\hrule} row\vrule depth3pt
height15pt width0pt&& $2$& $m$\cr column\vrule depth3pt
height15pt width0pt&& $m$& $2$\cr }} \qquad \vcenter
{\offinterlineskip \halign {\hfil $\;$#& \quad\vrule\quad#& \hfil
#\hfil \quad&  \hfil #$\;$\hfil \cr $c_2$\vrule depth6pt height9pt
width0pt&&  $p =1$& $p =\infty$\cr \noalign {\hrule} row\vrule
depth3pt height15pt width0pt&& $(2m-1)^2$& $m^2$\cr
column\vrule depth3pt height15pt width0pt&& $m^3$& $2m-1$\cr
}}$$ where $n$ is the order of $A$ and each row of $A$ stretches to
at most $m$ rows of $A^S$. When the sequence of stretchings is
disjoint in that later stretchings do not alter the rows or columns
of earlier stretchings, then the tables can be replaced by the ones
below. $$\vcenter {\offinterlineskip \halign {\hfil $\;$#& \quad
\vrule \quad#& \hfil #\hfil \quad& \hfil #$\;$\hfil \cr  $c_1$\vrule
depth6pt height9pt width0pt&& $p =1$&  $p =\infty$\cr \noalign
{\hrule} row\vrule depth3pt height15pt width0pt&& $2$& $2$\cr
column\vrule depth3pt height15pt width0pt&& $2$& $2$\cr }}
\qquad \vcenter {\offinterlineskip \halign {\hfil $\;$#& \quad
\vrule \quad#& \hfil #\hfil \quad& \hfil #$\;$\hfil \cr $c_2$\vrule
depth6pt height9pt width0pt&&  $p =1$& $p =\infty$\cr \noalign
{\hrule} row\vrule depth3pt height15pt width0pt&& $(2m-1)^2$&
$3m$\cr column\vrule depth3pt height15pt width0pt&& $3m^2$&
$2m-1$\cr }}$$ Thus, if $A$ is well-conditioned, if $\epsilon$ and $\|
E \|_p / \| \AS \|_p$ are very small, and if $m$ and $n$ are not
excessively large, then $\barz_\S$ is a good approximate solution
to $A x = y$.

{\sl Proof.} Lemma~1 shows $\AS = A^S + E_{5.1}$ with $\|
E_{5.1} \| \le \delta_1 \| A \|$ hence $\| E_{5.1} \| \le \delta_1 \| A^S
\|$ because Theorem~5's glue makes $\| A \| \le \| A^S \|$. Therefore
$$\eqalign {\| E_{5.1} + E \|& \le \| E_{5.1} \| + \| E \|\cr & = \| E_{5.1}
\| + \delta_2 \| \AS \|\cr & \le \delta_1 \| A^S \| + \delta_2 (1 +
\delta_1) \| A^S \|\cr & = (\delta_1 + \delta_2 + \delta_1
\delta_2) \| A^S \|\cr }$$ in which the hypotheses for $\delta_1$
and $\delta_2$ assure $\delta_1 + \delta_2 + \delta_1 \delta_2 <
1$.\footnote{$^5$}{The bound on $\delta_1$ also implies $\AS$ is
nonsingular, but this fact is not needed.} From $$(A^S + E_{5.1} + E)
\barz = (\AS + E) \barz = y^S$$ the matrix perturbation inequality
of Section~4 now implies $${\| z - \barz \| \over \| z \|} \le \kappa
(A^S) \, {\| E_{5.1} + E \| \over \| A^S \| - \| E_{5.1} + E \|} < \kappa
(A^S) \, {(\delta_1 + \delta_2 + \delta_1 \delta_2) \over 1 -
(\delta_1 + \delta_2 + \delta_1 \delta_2)}$$ while Lemma~2 and
Theorem~5 $${\| x - \barx \| \over \| x \|} \le c_{5.2} \, {\| z - \barz \|
\over \| z \|} \qquad \hbox {and} \qquad \kappa_p (A^S) \le c_4 \,
\kappa_p (A)$$ complete the chain of inequalities with $c_2 =
c_{5.2} c_4$. {\sl End of proof.}

\Proclaim Theorem~7. An $n \times n$, dense system of linear
equations can be solved by triangular factorization with row
reordering for stability using $$\vcenter {\openup1\jot
\halign{\hfil #\quad& #\hfil \cr $2 n^3 / 3 - 2n/3$& arithmetic
operations for the factorization and\cr $2n^2-n$& operations for
the solution phase.\cr }}$$ However, if the matrix is banded with
strict lower and upper bandwidths $\ell$ and $u$, and if $\ell + u <
n$, then the operations reduce to $$\vcenter {\openup1\jot
\halign{\hfil #\quad& #\hfil \cr $2 \ell (\ell + u + 1) n - \ell (4 \ell^2
+ 6 \ell u + 3 u^2 + 6 \ell + 3 u + 2)/3$& for the factorization
and\cr $(4 \ell + 2 u + 1) n - (2 \ell^2 + 2 \ell u + u^2 + 2 \ell + u)$& 
for the solution phase.\cr }}$$

{\sl Proof.} The factorization algorithm proceeds down the main
diagonal of the coefficient matrix $A$ by subtracting multiples of
the row in which the diagonal entry lies from lower rows to place
zeroes in the column beneath the diagonal entry. This column is
first inspected and the rows reordered to insure that the diagonal
entry has larger magnitude than any below. At the $k^{th}$
diagonal entry of the dense matrix, $n-k$ comparisons select the
largest entry, $n-k$ divisions form the multipliers, and $(n-k)^2$
each of multiplications and subtractions perform the row
operations. The total of these over all diagonal entries is
$$\sum_{k =1}^n \left [2(n-k) + 2(n-k)^2 \right] = 2 n^3 / 3 - 2 n /
3.$$ There results a factorization of the reordered matrix, $PA =
LU$, in which $P$ is the permutation matrix of the row reordering,
$L$ is the unit lower triangular matrix of multipliers, and $U$ is
the upper triangular matrix that contains what remains of the
original matrix. A particular system of linear equations $A x = y$
can then be written as $P^t L U x = y$ and solved by substitution.
Substitution with $L$ performs one multiplication and one
subtraction for each entry below the main diagonal. Substitution
with $U$ additionally divides by each diagonal entry. Altogether
$$(n-1) n + n + (n-1) n = 2 n^2 - n$$ operations are needed to
obtain the solution. 

The lower and upper triangular factors of banded matrices inherit
the lower and upper bandwidths of the original. The effect of row
reordering is merely to increase the upper bandwidth by the
lower bandwidth [\Golub ]. The reordering performs one
comparison, and the preparation of multipliers performs one
division, for each lower diagonal entry. These account for
$$\sum_{j =1}^\ell 2 (n - j) = 2 \ell n - (\ell^2 + \ell).$$ operations.
The $j$ in this summation indexes the lower diagonals. The row
operations subtract multiples of each strictly upper triangular
entry from the $\ell$ entries immediately below. With $\ell + u <
n$ upper diagonals, the first $n - \ell$ rows account for $$2 \ell
\big[ (\ell + u) (n - \ell) - u (u + 1) / 2 \big] = 2 \ell (\ell + u) n -
\ell (2 \ell^2 + 2 \ell u + u^2 + u)$$ operations, and the final $\ell$
rows account for $$\sum_{k =1}^\ell 2 (\ell-k)^2 = \ell (2 \ell^2 - 3
\ell + 1) / 3 .$$ The sum of the three expressions above simplifies
to the Theorem's formula for the factorization phase. By
reasoning similar to the dense case, the substitution phase of the
banded case performs $$2 \left[ \ell n - {\ell (\ell + 1) \over 2}
\right] + n + 2 \left[ (\ell + u) n - {(\ell + u) (\ell + u + 1) \over 2}
\right]$$ operations. {\sl End of proof.} 

\Proclaim Theorem~8. This row and column partitioning makes a
banded matrix into a block-bidiagonal one. For a matrix of order
$n$ with strict lower and upper bandwidths $\ell$ and $u$, and with
$0 < \ell + u < n$, the columns and rows partition into blocks of the
following size. {\def \Ls {,{\kern 0.75em}} $$\vcenter
{\openup1\jot \halign {\hfil #\hquad \quad& $#$\hfil \cr columns& 
\listfive {a + u} {\ell + u} {\ldots} {\ell + u} {\ell + c} \endlist\cr 
rows& \listsix {a} {u + \ell} {u + \ell} {\ldots} {u + \ell} {c} \endlist
\cr \noalign {\medskip} \multispan2 {\hfil $0 \le a \le \ell \qquad 0
\le c \le u \qquad 0 < a + c$ \hfil }\cr }}$$}The block-column
di\-men\-sion is $m = \lceil n / (\ell + u) \rceil$, and the
block-row dimension is $m + 1$ or $m$ (since one of $a$ or $c$ may
be zero). Moreover, the upper diagonal blocks are lower
triangular and the lower diagonal blocks are upper triangular.

{\sl Proof.} Since $n$ is strictly greater than $\ell + u$ it can be
decomposed as $$n = \left \lceil {n - \ell - u \over \ell + u} \right
\rceil (\ell + u) + r$$ where $0 < r \le \ell + u$. Thus, any of the sums
$r = a + c$ with $0 \le a \le \ell$ and $0 \le c \le u$ produce the
claimed partitioning of the columns by way of the decomposition
$$n = (a + u) + (m - 2) (\ell + u) + (\ell + c)$$ in which $$m = \left
\lceil {n - \ell - u \over \ell + u} \right \rceil + 1 = \left \lceil {n
\over \ell + u} \right \rceil .$$ There are exactly $m$ blocks of
columns because neither $a + u$ nor $\ell + c$ can be zero, for
example, $0 < r = a + c \le a + u$. The row partitioning stems
similarly from the decomposition $$n = a + (m - 1) (\ell + u) + c$$
with the precise number of blocks varying from $m + 1$ to $m$
because one of $a$ and $c$ can be zero.

The matrix can be enlarged by placing $\ell + u - a$ rows of zeroes
at the top and $\ell + u - c$ at the bottom, $\ell - a$ columns of
zeroes at the left and $u - c$ at the right. The augmented matrix
has block dimension $(m + 1) \times m$ in which every block is
square of order $\ell + u$. Entry $(i, j\,)$ of the original matrix
becomes entry $(\newi, \newj\,) = (\ell + u - a + i, \ell - a + j\,)$ of
the larger matrix. If the entry is nonzero, then $- \ell \le j - i \le u$
because the original matrix has lower and upper bandwidths $\ell$
and $u$, and then by a little arithmetic $- (\ell + u) \le \newj -
\newi \le 0$. Decomposing $$\eqalign {\newi& = \lfloor \newi /
(\ell + u) \rfloor (\ell + u) + b\cr \newj& = \lfloor \newj / (\ell + u)
\rfloor (\ell + u) + d\cr }$$ in which $0 \le (\hbox {$b$ and $d$}) <
\ell + u$, it then follows that $$- 2(\ell + u) < - (\ell + u) + d - b \le
\newj - \newi + d - b = \left( \lfloor \newj / (\ell + u) \rfloor +
\lfloor \newi / (\ell + u) \rfloor \right) (\ell + u) .$$ Thus $- 1 \le
\lfloor \newj / (\ell + u) \rfloor - \lfloor \newi / (\ell + u) \rfloor
\le 0$, and since $\lfloor \newi / (\ell + u) \rfloor$ and $\lfloor
\newj / (\ell + u) \rfloor$ are the block indices of the nonzero
entry, the inequality above means the entry lies in either a
diagonal block or a block immediately below a diagonal block. The
matrix therefore is block bidiagonal. The earlier inequality $-(\ell
+ u) \le \newj - \newi \le 0$ means the nonzero region extends
from the main diagonal of the subdiagonal blocks up to the main
diagonal of the main diagonal blocks. The subdiagonal blocks
therefore are upper triangular, and the main diagonal blocks are
lower triangular. {\sl End of proof.}

\Proclaim Theorem~9. An order $n + d$, bordered, banded system
of linear equations, whose coefficient matrix has $d$ dense rows
and columns in the bordering portion and has strict lower and
upper bandwidths $\ell$ and $u$ in the $n \times n$ banded portion,
where $0 < \ell + u < n$, can be solved by simple row stretching
and triangular factorization with row reordering for stability in
$$\displaylines {(4 d^2 + 6 d \ell + 2 d u + 2 \ell^2 + 2 \ell u + 2 d +
2 \ell) N \cr {}- (d + \ell) (13 d^2 + 14 d \ell + 12 d u + 4 \ell^2 + 6
\ell u + 3 u^2 + 9 d + 6 \ell + 3 u + 2) / 3\cr }$$ arithmetic
operations for the factorization and $$\displaylines {(4 d + 4 \ell +
2 u + 1) N \cr {}- (2 d^2 + 4 d \ell + 2 d u + 2 \ell^2 + 2 \ell u + u^2 +
2 d + 2 \ell + u)\cr }$$ operations for the solution phase, in which
$$N = n + d \left \lceil {n \over \ell + u} \right \rceil$$ is the
size of the stretched matrix.

{\sl Proof.} The simple row stretching of Definition 2, when based
upon the partitioning of Theorem~8 and when accompanied by the
reordering described in the text, results in a coefficient matrix
with order $N$, with strict lower and upper bandwidths $d + \ell$
and $u$, and with $d$ dense columns but no dense rows. Without
reordering, the triangular factors inherit this nonzero pattern. 
With row reordering, the upper bandwidth increases by the lower
bandwidth as in the purely banded case because there are no
dense rows. The results of Theorem~7 therefore account for
everything but the portion of the dense columns outside the
increased band. Ignoring these, Theorem~7 reports the
factorization phase performs $$\displaylines {2 (d + \ell) (d + \ell
+ u + 1) N\cr {} - (d + \ell) \left[ 4 (d + \ell)^2 + 6 (d + \ell) u + 3
u^2 + 6 (d + \ell) + 3 u + 2 \right] / \, 3}$$ operations, and the
solution phase performs $$(4 \left[ d + \ell) + 2 u + 1\right] N -
\left[ 2 (d + \ell)^2 + 2 (d + \ell) u + u^2 + 2 (d + \ell) + u \right] .$$
The dense columns outside the increased band contain $d \, (2N -
3d - 2\ell - 2u - 1) / 2$ entries. The factorization phase subtracts
multiples of these from the $d + \ell$ entries below, for $$d (d +
\ell) (2N - 3d - 2\ell - 2u - 1)$$ more operations. The solution
phase performs two operations for each strictly upper triangular
entry, for $$d (2N - 3d - 2\ell - 2u - 1)$$ more operations. The
sums of the expressions above simplify to the formulas in the
statement of the Theorem. {\sl End of proof.}

\vfil \eject \beginsection {Appendix 2. Figure Explanations}

This appendix explains the numerical experiments reported in the
Figures. All calculations are performed by a Cray {\sevenrm
Y-MP}8/264 with unit roundoff $3{.}5 \times 10^{-15}$. Matrix
factorizations, solutions of linear equations, and singular values
are computed using Linpack's {\sevenrm SGEFA}, {\sevenrm
SGESL} and {\sevenrm SSVDC} [\Dongarra ].

\figure 6. $2$-norm condition numbers for parameterized
matrices of order $51$ with sparsity patterns like the matrix in
Figure~2. 

The parameterized matrices are tridiagonal except in their final
rows and columns where all entries equal $1$. The tridiagonal
portion is a Toeplitz matrix with $-1$ and $-2$ on the lower and
upper diagonals, and with the parameter on the main diagonal.
There are $1201$ matrices for parameter values uniformly
distributed from $-6$ to $6$. The condition numbers for
parameter values between $\pm 3$ exceed surrounding condition
numbers by two orders of magnitude, but are still no worse than
$10^3$. A parameter value near $-3$ evidently produces a singular
matrix, and nearby values produce matrices with high condition
numbers---but not so high to trouble a machine with a $3{.}5
\times 10^{-15}$ unit roundoff. The Figure's vertical axis has been
lengthened to ease comparison with later figures. 

\figure 7. Maximum $2$-norm relative errors for equations $A x =
y$, with $20$ different $y$'s and the parameterized matrices $A$
of Figure~6, solved by triangular factorization. The lower curve
allows full row reordering. The upper curve restricts row
reordering to the tridiagonal band. 

\figure 8. Percent of non-zeroes in the triangular factors of the
matrices of Figure~6. The upper curve allows full row reordering.
The lower curve restricts row reordering to the tridiagonal band.

Many right hand sides reduce the possibility of serendipity and
smooth the curves by removing occasional outliers. In a sense,
Figure~7 consists only of outliers because it reports the
maximum error for any of the vectors, which repeat for each
matrix. The vectors have entries uniformly and randomly
distributed between $\pm 1$. With limited row reordering, the
errors for parameter values between $\pm 3$ exceed surrounding
errors by a greater margin than do the condition numbers in
Figure~6.

Restricted row reordering is performed by a modified {\sevenrm
SGEFA} which limits its search for elimination rows as though the
matrices were tridiagonal. Figure~8's percentages omit the main
diagonals of the lower triangular factors, which are identically $1$
and not stored. The downward spikes in Figure~8 indicate a few
matrices have abnormally sparse factors. 

\figure 9. $2$-norm relative errors for the equations of Figure~7
solved by triangular factorization with full row reordering after
stretching in the manner of Figure~3.

\figure 10. Percent of non-zeroes in the triangular factors of the
stretched matrices of Figure~9. The percentages are relative to
the size of the unstretched matrices.

The Introduction cites Figures~9 and~10 but later sections more
precisely explain the stretching. Section~3 defines simple row
stretching. Section~5 applies the stretching to bordered, banded
matrices. Theorem~9 shows the stretched matrices have size $$N
= n + d \left \lceil {n \over \ell + u} \right \rceil = 50 + 1 \left
\lceil {50 \over 1 + 1} \right \rceil = 75$$ which is relatively much
larger than the original matrices only because the bandwidth is
very small in this example.

Figure~9 resembles Figure~7. The equations' right hand sides are
those of the earlier Figure stretched by in\-serting zeroes as
explained in Section~3. The relative errors are for the original
variables, that is, they exclude the extraneous new variables also
inserted by the stretching.

Figure~10 resembles Figure~8. Nonzeroes of the $75 \times 75$
factors are reported as percentages of the $51^2$ entries of the
unstretched matrices. The percentages again omit the identically
$1$ main diagonal of one factor. 

\figure 12. Condition numbers of the matrices in Figure~6 (lower
solid lines) and condition numbers after stretching (dashed lines)
to remove either bordering rows or columns. Theorem~5
specifies glue that bounds (upper solid lines) either the $1$- or
$\infty$-norm condition numbers.

The inverse matrices are explicitly formed to evaluate
the $1$-norm and $\infty$-norm condition numbers. The multiplier
$c$ in Theorem~5 is $2m-1$ or $3m$. $$m = \left \lceil {n
\over \ell + u} \right \rceil = \left \lceil {50 \over 1 + 1} \right
\rceil = 25$$

\figure 13. Smallest pivot and singular value for the banded
portion of the matrices in Figure~6. Table~2's version of deflated
block elimination assumes the pivot and singular value have the
same magnitude ``which is definitely not valid in general, but
which is shown empirically and theoretically to be valid in
practice'' {\rm [\Chan,~p.~124].}

\figure 14. Maximum $2$-norm relative errors for equations $A x
= y$ with $20$ different $y$'s solved by block elimination,
deflated block elimination, and simple row stretching. The
parameterized coefficient matrices $A$ are those of Figure~6.
The elimination methods cannot be distinguished at this plotting
resolution. The stretching data also appears in Figure~9.

{\sevenrm SGEFA}'s factorization is $B = P (I + L) U$ in which $P$ is
a permutation matrix, $L$ is strictly lower triangular and $U$ is
upper triangular. The smallest pivot is the entry $u_{k,k}$ of
smallest magnitude on the main diagonal of $U$. Chan [\Chan ]
offers no guidance on the proper choice of the pivot's index. It
appears to be $k$---disregarding the row permutation---because
the deflated algorithm applies $B^{-t} = P^t (I + L)^{-t} U^{-t}$ to
column $k$ of an identity matrix. 

Figure~14 uses the same vectors $y$ and reports the solution
errors in the same manner as Figure~7. Block elimination and
deflated block elimination are performed as shown in Tables~1
and~2. 

 \vfill \eject

{
 \headline={\hfil}

 \def\footremark{{\ }}
 \def\pagenumber{{\tenrm \the \pageno}}
 \footline={\ifnum \pageno=1{\hfil}\else {\ifnum
 \pageno=2{\hfil}\else {\ifnum \pageno=5{\noindent \hfil
 \footremark \hfil {\tenrm 5/6}}\else {\ifnum
 \pageno=6{\hfil}\else {\ifodd \pageno {\noindent \hfil
 \footremark \hfil \pagenumber}\else {\noindent \pagenumber
 \hfil \footremark \hfil}\fi}\fi}\fi}\fi}\fi}

 \baselineskip = 12pt
 \def \inc #1{\global \advance \number by#1}
 \hsize = 2.45in
 \newbox \leftcolumn
 \newcount \number \number = 0
 \newdimen \fullhsize \fullhsize = 5.0in
 \pageno = 56
 \parindent = 0pt
 \parskip = 6pt plus 9pt minus 3pt

 \def \columnbox{\leftline{\pagebody}}
 \let \lr = L 
 \def \doubleformat{\shipout \vbox {\makeheadline
      \hbox to\fullhsize{\box \leftcolumn \hss \columnbox}
   \makefootline}
   \advancepageno}
 \output={\if L\lr
      \global \setbox \leftcolumn = \columnbox \global \let \lr = R
   \else \doubleformat \global \let \lr = L\fi
   \ifnum \outputpenalty>-20000\else \dosupereject \fi}

 \def \eol {\hfil \break}
 \def \Eblock#1\par{\inc 1\noindent#1\filbreak}
 \def \Lblock#1\par{\noindent#1}
 
UNLIMITED RELEASE\par INITIAL DISTRIBUTION 

\frenchspacing \helvetica 
{\pretolerance = 10000

\Eblock{I. K. Abu-Shumays\eol
Bettis Atomic Power Lab.\eol
Box 79\eol
West Mifflin, PA 15122}

\Eblock{Loyce M. Adams\eol
ETH-Zentrum\eol
Inst. Wiss. Rechnen\eol
CH-8092 Zurich\eol
Switzerland}

\Eblock{Peter W. Aitchison\eol
Univ. of Manitoba\eol
Applied Math. Dept.\eol
Winnipeg, Manitoba R3T 2N2\eol
Canada}

\Eblock{Leena Aittoniemi\eol
Tech. Univ. Berlin\eol
Comp. Sci. Dept.\eol
Franklinstrasse 28-29\eol
D-1000 Berlin 10\eol
Germany}

\Eblock{Fernando L. Alvarado\eol
Univ. of Wisconsin\eol
Electrical and Comp. Eng.\eol
1425 Johnson Drive\eol
Madison, WI 53706}

\Eblock{Patrick Amestoy\eol
CERFACS\eol
42 ave g Coriolis\eol
31057 Toulouse\eol
France}

\Eblock{Ed Anderson\eol
Argonne National Lab.\eol
Math. and Comp. Sci. Div.\eol
9700 S Cass Av.\eol
Argonne, IL 60439}

\Eblock{Johannes Anderson\eol
Brunel Univ.\eol
Dept. of Math. and Statistics\eol
Uxbridge Middlesex UB8 3PH\eol
United Kingdom}

\Eblock{Peter Arbenz\eol
ETH-Zentrum\eol
Inst. fur Informatik\eol
CH-8092 Zurich\eol
Switzerland}

\Eblock{Mario Arioli\eol
CERFACS\eol
42 ave g Coriolis\eol
31057 Toulouse\eol
France}

\Eblock{Cleve Ashcraft\eol
Boeing Comp. Services\eol
Mail Stop 7L-21\eol
P. O. Box 24346\eol
Seattle, WA 98124-0346}

\Eblock{Owe Axelsson\eol
Katholieke Univ.\eol
Dept. of Math.\eol
Toernooiveld\eol
6525 ED Nijmegen\eol
The Netherlands}

\Eblock{Zhaojun Bai\eol
New York Univ.\eol
Courant Institute\eol
251 Mercer Street\eol
New York, NY 10012}

\Eblock{R. E. Bank\eol
Univ. of California, San Diego\eol
Dept. of Math.\eol
La Jolla, CA  92093}

\Eblock{Jesse L. Barlow\eol
Pennsylvania State Univ.\eol
Dept. of Comp. Sci.\eol
Univ. Park, PA 16802}

\Eblock{Chris Bischof\eol
Argonne National Lab.\eol
Math. and Comp. Sci. Div.\eol
9700 S Cass Av.\eol
Argonne, IL 60439}

\Eblock{R. H. Bisseling\eol
Koninklijke/Shell\eol
Laboratorium Amsterdam\eol
P. O. Box 3003\eol
1003 AA Amsterdam\eol
The Netherlands}

\Eblock{Ake Bjorck\eol
Linkoping Univ.\eol
Dept. of Math.\eol
S-581 83 Linkoping\eol
Sweden}

\Eblock{Petter E. Bjorstad\eol
Univ. of Bergen\eol
Thormonhlensgt. 55\eol
N-5006 Bergen\eol
Norway}

\Eblock{Daniel Boley\eol
Stanford Univ.\eol
Dept. of Comp. Sci.\eol
Stanford, CA 94305}

\Eblock{Randall Bramley\eol
Univ. of Illinois\eol
305 Talbot Lab.\eol
104 S Wright Street\eol
Urbana, IL 61801}

\Eblock{Richard A. Brualdi\eol
Univ. of Wisconsin\eol
Dept. of Math.\eol
480 Lincoln Drive\eol
Madison, WI 53706}

\Eblock{James R. Bunch\eol
Univ. of California, San Diego\eol
Dept. of Math.\eol
La Jolla, CA  92093}

\Eblock{Angelica Bunse-Gerstner\eol
Univ. Bielefeld\eol
Fakultat fur Math.\eol
postfach 8640\eol
D-4800 Bielfeld 1\eol
Germany}

\Eblock{Ralph Byers\eol
Univ. of Kansas\eol
Dept. of Math.\eol
Lawrence, KA 66045-2142}

\Eblock{Jose Castillo\eol
San Diego State Univ.\eol
Dept. of Math. Sciences\eol
San Diego, CA  92182}

\Eblock{Tony F. Chan\eol
Univ. of California, Los Angles\eol
Dept. of Math.\eol
Los Angeles, CA 90024}

\Eblock{S. S. Chow\eol
Univ. of Wyoming\eol
Dept. of Math.\eol
Laramie, WY 82071}

\Eblock{Anthony T. Chronopoulos\eol
Univ. of Minnesota\eol
Dept. of Comp. Sci.\eol
200 Union Street SE\eol
Minneapolis, MN 55455-0159}

\Eblock{Eleanor C. H. Chu\eol
Univ. of Waterloo\eol
Dept. of Comp. Sci.\eol
Waterloo, Ontario N2L 3G1\eol
Canada}

\Eblock{Len Colgan\eol
South Australian Inst. of Tech.\eol
Math. Dept.\eol
The Levels, 5095, South Australia\eol
Australia}

\Eblock{Paul Concus\eol
Univ. of California, Berkeley\eol
Lawrence Berkeley Lab.\eol
50A-2129\eol
Berkeley, CA 94720}

\Eblock{W. M. Coughran\eol
AT\&T Bell Laboratories\eol
600 Mountain Av.\eol
Murray Hill, NJ 07974-2070}

\Eblock{Jane K. Cullum\eol
IBM T. J. Watson Res. Center\eol
P. O. Box 218\eol
Yorktown Heights, NY 10598}

\Eblock{Carl De Boor\eol
Univ. of Wisconsin\eol
Center for Math. Sciences\eol
610 Walnut Street\eol
Madison, WI 53706}

\Eblock{John E. De Pillis\eol
Univ. of California, Riverside\eol
Dept. of Math. and Comp. Sci.\eol
Riverside, CA 92521}

\Eblock{James W. Demmel\eol
New York Univ.\eol
Courant Institute\eol
251 Mercer Street\eol
New York, NY 10012}

\Eblock{Julio Cesar Diaz\eol
Univ. of Tulsa\eol
Dept. of Math. and Comp. Sci.\eol
600 S College Av.\eol
Tulsa, OK 74104-3189}

\Eblock{David S. Dodson\eol
Convex Comp. Corp.\eol
701 N Plano Rood\eol
Richardson, TX 75081}

\Eblock{Jack Dongarra\eol
Univ. of Tennessee\eol
Comp. Sci. Dept.\eol
Knoxville, TN 37996-1300}

\Eblock{Ian S. Duff\eol
Rutherford Appleton Lab.\eol
Numerical Analysis Group\eol
CC Dept.\eol
OXON OX11 OQX\eol
United Kingdom}

\Eblock{Pat Eberlein\eol
State Univ. of New York\eol
Dept. of Comp. Sci.\eol
Buffalo, NY 14260}

\Eblock{W. Stuart Edwards\eol
Univ. of Texas\eol
Center for Nonlinear Dynamics\eol
Austin, TX 78712}

\Eblock{Louis W. Ehrlich\eol
Johns Hopkins Univ.\eol
Applied Physics Lab.\eol
Johns Hopkins Road\eol
Laurel, MD 20707}

\Eblock{Michael Eiermann\eol
Univ. Karlsruhe\eol
Inst. fur Prakticshe Math.\eol
Postfach 6980\eol
D-6980 Karlsruhe 1\eol
Germany}

\Eblock{Stanley C. Eisenstat\eol
Yale University\eol
Dept. of Comp. Sci.\eol
P. O. Box 2158, Yale Station\eol
New Haven, CT 06520}

\Eblock{Lars Elden\eol
Linkoping Univ.\eol
Dept. of Math.\eol
S-581 83 Linkoping\eol
Sweden}

\Eblock{Howard C. Elman\eol
Univ. of Maryland\eol
Dept. of Comp. Sci.\eol
College Park, MD 20742}

\Eblock{Albert M. Erisman\eol
Boeing Comp. Services\eol
565 Andover Park West\eol
Mail Stop 9C-01\eol
Tukwila, WA 98188}

\Eblock{D. J. Evans\eol
Univ. of Technology\eol
Dept. of Computer Studies\eol
Leicestershire LE11 3TU\eol
United Kingdom}

\Eblock{Vance Faber\eol
Los Alamos National Lab.\eol
Group C-3\eol
Mail Stop B265\eol
Los Alamos, NM 87545}

\Eblock{Bernd Fischer\eol
Univ. Hamburg\eol
Inst. fur Angew. Math.\eol
D-2000 Hamburg 13\eol
Germany}

\Eblock{Geoffrey Fox\eol
California Inst. of Tech.\eol
Mail Code 158-79\eol
Pasadena, CA 91125}

\Eblock{Paul O. Frederickson\eol
NASA Ames Res. Center\eol
RIACS, Mail Stop 230-5\eol
Moffett Field, CA 94035}

\Eblock{Roland W. Freund\eol
NASA Ames Res. Center\eol
RIACS, Mail Stop 230-5\eol
Moffett Field, CA 94035}

\Eblock{Robert E. Funderlic\eol
North Carolina State Univ.\eol
Dept. of Comp. Sci.\eol
Raleigh, NC 27650}

\Eblock{Ralf Gaertner\eol
Max Plank Society\eol
Fritz Haber Inst.\eol
Faradayweg 4-6\eol
D-1000 Berlin 33\eol
Germany}

\Eblock{Patrick W. Gaffney\eol
Bergen Scientific Center\eol
Allegaten 36\eol
N-5000 Bergen\eol
Norway}

\Eblock{Kyle A. Gallivan\eol
Univ. of Illinois\eol
305 Talbot Lab.\eol
104 S Wright Street\eol
Urbana, IL 61801}

\Eblock{Dennis B. Gannon\eol
Indiana Univ.\eol
Dept. of Comp. Sci.\eol
Bloomington, IN 47405-6171}

\Eblock{Kevin E. Gates\eol
Univ. of Washington\eol
Dept. of Applied Math.\eol
Seattle, WA 98195}

\Eblock{David M. Gay\eol
AT\&T Bell Laboratories\eol
600 Mountain Av.\eol
Murray Hill, NJ 07974-2070}

\Eblock{C. W. Gear\eol
NEC Research Institute\eol
4 Independence Way\eol
Princeton, NJ 08540}

\Eblock{C. William Gear\eol
Univ. of Illinois\eol
Dept. of Comp. Sci.\eol
1304 W Springfield Av.\eol
Urbana, IL 61801}

\Eblock{George A. Geist\eol
Oak Ridge National Lab.\eol
Math. Sciences Section\eol
Building 9207-A\eol
P. O. Box 2009\eol
Oak Ridge, TN 37831}

\Eblock{J. Alan George\eol
Univ. of Waterloo\eol
Needles Hall\eol
Waterloo, Ontario N2L 3G1\eol
Canada}

\Eblock{Adam Gersztenkorn\eol
Amoco Production Company\eol
Geophysical Res. Dept.\eol
4502 E 41st Street\eol
P. O. Box 3385\eol
Tulsa, OK 74102}

\Eblock{John R. Gilbert\eol
Xerox Palo Alto Res. Center\eol
3333 Coyote Hill Road\eol
Palo Alto CA 94304}

\Eblock{Albert Gilg\eol
Siemens AG\eol
Corporate Res. and Tech.\eol
Otto Hahn Ring 6\eol
D-8000 Munchen 83\eol
Germany}

\Eblock{Roland Glowinski\eol
Univ. of Houston\eol
Dept. of Math.\eol
4800 Calhoun Road\eol
Houston, TX 77004}

\Eblock{Gene H. Golub\eol
Stanford Univ.\eol
Dept. of Comp. Sci.\eol
Stanford, CA 94305}

\Eblock{W. Gragg\eol
Naval Postgraduate School\eol
Dept. of Math.\eol
Mail Code 53ZH\eol
Monterey, CA 93943-5100}

\Eblock{Anne Greenbaum\eol
New York Univ.\eol
Courant Institute\eol
251 Mercer Street\eol
New York, NY 10012}

\Eblock{Roger Grimes\eol
Boeing Comp. Services\eol
Mail Stop 7L-21\eol
P. O. Box 24346\eol
Seattle, WA 98124-0346}

\Eblock{William D. Gropp\eol
Argonne National Lab.\eol
Math. and Comp. Sci. Div.\eol
9700 S Cass Av.\eol
Argonne, IL 60439}

\Eblock{Fred G. Gustavson\eol
IBM T. J. Watson Res. Center\eol
P. O. Box 218\eol
Yorktown Heights, NY 10598}

\Eblock{Martin H. Gutknecht\eol
ETH-Zentrum\eol
IPS, IFW D 25.1\eol
CH-8092 Zurich\eol
Switzerland}

\Eblock{Louis Hageman\eol
Bettis Atomic Power Lab.\eol
Box 79\eol
West Mifflin, PA 15122}

\Eblock{Charles A. Hall\eol
Univ. of Pittsburgh\eol
Dept. of Math. and Statistics\eol
Pittsburgh, PA 15260}

\Eblock{Sven J. Hammarling\eol
Numerical Algorithms Group Ltd.\eol
Wilkinson House\eol
Jordan Hill Road\eol
Oxford OX2 8DR\eol
United Kingdom}

\Eblock{Ken Hanson\eol
Los Alamos National Lab.\eol
Mail Stop P940\eol
Los Alamos, New Mexico 87545}

\Eblock{Michael T. Heath\eol
Oak Ridge National Lab.\eol
Math. Sciences Section\eol
Building 9207-A\eol
P. O. Box 2009\eol
Oak Ridge, TN 37831-8083}

\Eblock{Don E. Heller\eol
Shell Development Company\eol
Bellaire Res. Center\eol
P. O. Box 481\eol
Houston, TX 77001}

\Eblock{Nicholas J. Higham\eol
Cornell Univ.\eol
Dept. of Comp. Sci.\eol
Ithaca, NY 14853}

\Eblock{Mary Hill\eol
USGS-WRD\eol
Denver Federal Center\eol
Mail Stop 413\eol
P. O. Box 25046\eol
Lakewood, CO 80225}

\Eblock{Mike Holst\eol
Univ. of Illinois\eol
Dept. of Comp. Sci.\eol
1304 W Springfield Av.\eol
Urbana, IL 61801}

\Eblock{James M. Hyman\eol
Los Alamos National Lab.\eol
Group T-7\eol
Mail Stop B284\eol
Los Alamos, NM 87545}

\Eblock{Ilse Ipsen\eol
Yale University\eol
Dept. of Comp. Sci.\eol
P. O. Box 2158, Yale Station\eol
New Haven, CT 06520}

\Eblock{Doug James\eol
7305 Mill Ridge Road\eol
Raleigh, NC 27613}

\Eblock{S. Lennart Johnsson\eol
Thinking Machines Inc.\eol
245 First Street\eol
Cambridge, MA 02142-1214}

\Eblock{Tom Jordan\eol
Los Alamos National Lab.\eol
Group C-3\eol
Mail Stop B265\eol
Los Alamos, NM 87545}

\Eblock{Wayne Joubert\eol
Univ. of Texas\eol
Center for Numerical Analysis\eol
Moore Hall 13.150\eol
Austin, TX 78712}

\Eblock{E. F. Kaasschieter\eol
Eindhoven Univ. of Tech.\eol
Dept. of Math. and Comp. Sci.\eol
P. O. Box 513\eol
5600 MB Eindhoven\eol
The Netherlands}

\Eblock{Bo Kagstrom\eol
Univ. of Umea\eol
Inst. of Information Processing\eol
S-901 87 Umea\eol
Sweden}

\Eblock{W. M. Kahan\eol
Univ. of California, Berkeley\eol
Dept. of Math.\eol
Berkeley, CA 94720}

\Eblock{Hans Kaper\eol
Argonne National Lab.\eol
Math. and Comp. Sci. Div.\eol
9700 S Cass Av.\eol
Argonne, IL 60439}

\Eblock{Linda Kaufman\eol
AT\&T Bell Laboratories\eol
600 Mountain Av.\eol
Murray Hill, NJ 07974-2070}

\Eblock{C. T. Kelley\eol
North Carolina State Univ.\eol
Dept. of Math.\eol
Raleigh, NC 27695-8205}

\Eblock{David Keyes\eol
Yale University\eol
Dept. of Mechanical Eng.\eol
P. O. Box 2159, Yale Station\eol
New Haven, CT 06520}

\Eblock{David R. Kinkaid\eol
Univ. of Texas\eol
Center for Numerical Analysis\eol
Moore Hall 13.150\eol
Austin, TX 78712}

\Eblock{Virgina Klema\eol
Massachusetts Inst. of Tech.\eol
Statistics Center\eol
Cambridge, MA 02139}

\Eblock{Steven G. Kratzer\eol
Supercomputing Res. Center\eol
17100 Sci. Drive\eol
Bowie, MD 20715-4300}

\Eblock{Edward J. Kushner\eol
Intel Scientific Computers\eol
15201 NW Greenbrier Pkwy.\eol
Beaverton, OR 97006}

\Eblock{John Lavery\eol
NASA Lewis Res. Center\eol
Mail Stop 5-11\eol
Cleveland, OH 44135}

\Eblock{Kincho H. Law\eol
Stanford Univ.\eol
Dept. of Civil Eng.\eol
Stanford, CA 94305}

\Eblock{Charles Lawson\eol
Jet Propulsion Lab.\eol
Applied Math. Group\eol
Mail Stop 506-232\eol
4800 Oak Grove Drive\eol
Pasadena, CA 91109}

\Eblock{Yannick L. Le Coz\eol
Rensselaer Polytechnic Inst.\eol
Electrical Eng. Dept.\eol
Troy, NY 12180}

\Eblock{Randy LeVeque\eol
ETH-Zentrum\eol
Sem. Angew. Math.\eol
CH-8092 Zurich\eol
Switzerland}

\Eblock{Steve L. Lee\eol
Univ. of Illinois\eol
Dept. of Comp. Sci.\eol
1304 W Springfield Av.\eol
Urbana, IL 61801}

\Eblock{Ben Leimkuhler\eol
Univ. of Kansas\eol
Dept. of Math.\eol
Lawrence, Kansas 66045}

\Eblock{Steve Leon\eol
Southeastern Massachusetts Univ.\eol
Dept. of Math.\eol
North Dartmouth, MA 02747}

\Eblock{Michael R. Leuze\eol
Oak Ridge National Lab.\eol
Dept. of Math. Sciences\eol
P. O. Box 2009\eol
Oak Ridge, TN 37831}

\Eblock{Stewart A. Levin\eol
Mobil Res. and Dev. Corp.\eol
Dallas, TX 75381-9047}

\Eblock{John G. Lewis\eol
Boeing Comp. Services\eol
Mail Stop 7L-21\eol
P. O. Box 24346\eol
Seattle, WA 98124-0346}

\Eblock{Antonios Liakopoulos\eol
System Dynamics Inc.\eol
1211 NW 10th Av.\eol
Gainesville, FL 32601}

\Eblock{Joseph W. H. Liu\eol
York Univ.\eol
Dept. of Comp. Sci.\eol
North York, Ontario M3J 1P3\eol
Canada}

\Eblock{Peter Lory\eol
Tech. Univ. Munchen\eol
Dept. of Math.\eol
POB 20 24 20\eol
Germany\eol
D-8000 Munchen 2\eol
Germany}

\Eblock{Robert F. Lucas\eol
Supercomputing Res. Center\eol
17100 Sci. Drive\eol
Bowie, MD 20715-4300}

\Eblock{Franklin Luk\eol
Cornell Univ.\eol
School of Electrical Eng.\eol
Ithica, NY 14853}

\Eblock{Tom Manteuffel\eol
Univ. of Colorado\eol
Comp. Math. Group\eol
Campus Box 170\eol
1200 Larimer Street\eol
Denver, CO 80204}

\Eblock{Robert M. Mattheij\eol
Eindhoven Univ. of Tech.\eol
Dept. of Math.\eol
5600 MB Eindhoven\eol
The Netherlands}

\Eblock{Paul C. Messina\eol
California Inst. of Tech.\eol
Mail Code 158-79\eol
Pasadena, CA 91125}

\Eblock{Gerard A. Meurant\eol
Centre d'Etudes de Limeil\eol
Service Mathematiques Appliquees\eol
Boite Postale 27\eol
94190 Villeneuve St. Georges\eol
France}

\Eblock{Carl D. Meyer\eol
North Carolina State Univ.\eol
Dept. of Math.\eol
Raleigh, NC 27650}

\Eblock{Ignacy Misztal\eol
Univ. of Illinois\eol
Dept. of Animal Sciences\eol
1207 W Gregory Dr.\eol
Urbana, IL 61801}

\Eblock{Gautam Mitra\eol
Brunel Univ.\eol
Dept. of Math. and Statistics\eol
Uxbridge Middlesex UB8 3PH\eol
United Kingdom}

\Eblock{Hans Mittleman\eol
Arizona State Univ.\eol
Dept. of Math.\eol
Tempe, AZ 85287}

\Eblock{Cleve Moler\eol
The Mathworks\eol
325 Linfield Place\eol
Menlo Park, CA 94025}

\Eblock{Ronald B. Morgan\eol
Univ. of Missouri\eol
Dept. of Math.\eol
Columbia, MO 65211}

\Eblock{Rael Morris\eol
IBM Corp.\eol
Almaden Res. Center\eol
K08 / 282\eol
650 Garry Road\eol
San Jose, CA 95120-6099}

\Eblock{Edmond Nadler\eol
Wayne State Univ.\eol
Dept. of Math.\eol
Detroit, MI 48202}

\Eblock{N. Nandakumar\eol
Univ. of Nebraska\eol
Dept. of Math. and Comp. Sci.\eol
Omaha, NB 68182}

\Eblock{Olavi Nevanlinna\eol
Helsinki Univ. of Tech.\eol
Inst. of Math.\eol
SF-02150 Espoo\eol
Finland}

\Eblock{Esmond G. Y. Ng\eol
Oak Ridge National Lab.\eol
Math. Sciences Section\eol
Building 9207-A\eol
P. O. Box 2009\eol
Oak Ridge, TN 37831}

\Eblock{Viet-nam Nguyen\eol
Pratt \& Whitney\eol
Dept. of Eng. and Comp. Applications\eol
1000 Marie Victorian\eol
Longueuil, Quebec J4G 1A1\eol
Canada}

\Eblock{Nancy Nichols\eol
Reading Univ.\eol
Dept. of Math.\eol
Whiteknights Park\eol
Reading RG6 2AX\eol
United Kingdom}

\Eblock{W. Niethammer\eol
Univ. Karlsruhe\eol
Inst. fur Praktische Math.\eol
Engelerstrasse 2\eol
D-7500 Karlsruhe\eol
Germany}

\Eblock{Takashi Nodera\eol
Keio Univ.\eol
Dept. of Math.\eol
3-14-1 Hiyoshi Kohoku\eol
Yokohama 223\eol
Japan}

\Eblock{Wilbert Noronha\eol
Univ. of Tennessee\eol
310 Perkins Hall\eol
Knoxville, TN 37996}

\Eblock{Bahram Nour-Omid\eol
Lockheed Palo Alto Res. Lab.\eol
Comp. Mechanics Section\eol
Org. 93-30 Bldg. 251\eol
3251 Hanover Street\eol
Palo Alto, CA 94304}

\Eblock{Dianne P. O'Leary\eol
Univ. of Maryland\eol
Dept. of Comp. Sci.\eol
College Park, MD 20742}

\Eblock{Julia Olkin\eol
SRI International\eol
Building 301, Room 66\eol
333 Ravenswood Av.\eol
Menlo Park, CA 94025}

\Eblock{Steve Olson\eol
Supercomputer Systems Inc\eol
1414 W Hamilton Av.\eol
Eau Claire, WI 54701}

\Eblock{Elizabeth Ong\eol
Univ. of California, Los Angles\eol
Dept. of Math.\eol
Los Angeles, CA 90024}

\Eblock{James M. Ortega\eol
Univ. of Virginia\eol
Dept. of Applied Math.\eol
Charlottesville, VA 22903}

\Eblock{Christopher C. Paige\eol
McGill Univ.\eol
School of Comp. Sci.\eol
McConnell Eng. Building\eol
3480 Univ.\eol
Montreal, Quebec H3A 2A7\eol
Canada}

\Eblock{M. C. Pandian\eol
IBM Corp.\eol
Numerically Intensive Computing\eol
Dept. 41U/276\eol
Neighborhood Road\eol
Kingston, NY 12401}

\Eblock{Roy Pargas\eol
Clemson Univ.\eol
Dept. of Comp. Sci.\eol
Clemson, SC 29634-1906}

\Eblock{Beresford N. Parlett\eol
Univ. of California, Berkeley\eol
Dept. of Math.\eol
Berkeley, CA 94720}

\Eblock{Merrell Patrick\eol
Duke Univ.\eol
Dept. of Comp. Sci.\eol
Durham, NC 27706}

\Eblock{Victor Pereyra\eol
Weidlinger Associates\eol
Suite 110\eol
4410 El Camino Real\eol
Los Altos, CA 94022}

\Eblock{Michael Pernice\eol
Univ. of Utah\eol
Utah Supercomputing Inst.\eol
3334 Merrill Eng. Bldg.\eol
Salt Lake City, UT 84112}

\Eblock{Barry W. Peyton\eol
Oak Ridge National Lab.\eol
Math. Sciences Section\eol
Building 9207-A\eol
P. O. Box 2009\eol
Oak Ridge, TN 37831}

\Eblock{Robert J. Plemmons\eol
North Carolina State Univ.\eol
Dept. of Math.\eol
Raleigh, NC 27650}

\Eblock{Jim Purtilo\eol
Univ. of Maryland\eol
Dept. of Comp. Sci.\eol
College Park, MD 20742}

\Eblock{Giuseppe Radicati\eol
IBM Italia\eol
Via Giorgione 159\eol
00147 Roma\eol
Italy}

\Eblock{A. Ramage\eol
Univ. of Bristol\eol
Dept. of Math.\eol
Univ. Walk\eol
Bristol BS8 1TW\eol
United Kingdom}

\Eblock{Lothar Reichel\eol
Univ. of Kentucky\eol
Dept. of Math.\eol
Lexington, KY 40506}

\Eblock{John K. Reid\eol
Rutherford Appleton Lab.\eol
Numerical Analysis Group\eol
CC Dept.\eol
OXON OX11 OQX\eol
United Kingdom}

\Eblock{John R. Rice\eol
Purdue Univ.\eol
Dept. of Comp. Sci.\eol
Layfayette, IN 47907}

\Eblock{Jeff V. Richard\eol
Science Applications International\eol
Mail Stop 34\eol
10260 Campus Point Drive\eol
San Diego, CA 92121}

\Eblock{Donald J. Rose\eol
Duke Univ.\eol
Dept. of Comp. Sci.\eol
Durham, NC 27706}

\Eblock{Vona Bi Roubolo\eol
Univ. of Texas\eol
Center for Numerical Analysis\eol
Moore Hall 13.150\eol
Austin, TX 78712}

\Eblock{Axel Ruhe\eol
Chalmers Tekniska Hogskola\eol
Dept. of Comp. Sci.\eol
S-412 96 Goteborg\eol
Sweden}

\Eblock{Youcef Saad\eol
Univ. of Minnesota\eol
Dept. of Comp. Sci.\eol
200 Union Street SE\eol
Minneapolis, MN 55455-0159}

\Eblock{P. Sadayappan\eol
Ohio State Univ.\eol
Dept. of Comp. Sci.\eol
Columbus, OH 43210}

\Eblock{Joel Saltz\eol
NASA Langley Res. Center\eol
ICASE\eol
Mail Stop 132-C\eol
Hampton, VA 23665}

\Eblock{Ahmed H. Sameh\eol
Univ. of Illinois\eol
305 Talbot Lab.\eol
104 S Wright Street\eol
Urbana, IL 61801}

\Eblock{Michael Saunders\eol
Stanford Univ.\eol
Dept. of Operations Res.\eol
Stanford, CA 94305}

\Eblock{Paul E. Saylor\eol
Univ. of Illinois\eol
Dept. of Comp. Sci.\eol
1304 W Springfield Av.\eol
Urbana, IL 61801}

\Eblock{Mark Schaefer\eol
Texas A \& M Univ.\eol
Dept. of Math.\eol
College Station, TX 77843}

\Eblock{U. Schendel\eol
Freie Univ. Berlin\eol
Inst. fur Math.\eol
Arnimallee 2-6\eol
D-1000 Berlin 33\eol
Germany}

\Eblock{Robert S. Schreiber\eol
NASA Ames Res. Center\eol
RIACS, Mail Stop 230-5\eol
Moffett Field, CA 94035}

\Eblock{Martin Schultz\eol
Yale University\eol
Dept. of Comp. Sci.\eol
P. O. Box 2158, Yale Station\eol
New Haven, CT 06520}

\Eblock{David St. Clair Scott\eol
Intel Scientific Computers\eol
15201 NW Greenbrier Parkway\eol
Beaverton, OR 97006}

\Eblock{Jeffrey S. Scroggs\eol
NASA Langley Res. Center\eol
ICASE\eol
Mail Stop 132-C\eol
Hampton, VA 23665}

\Eblock{Steven M. Serbin\eol
Univ. of Tennessee\eol
Dept. of Math.\eol
Knoxville, TN 37996-1300}

\Eblock{Shahriar Shamsian\eol
MacNeal Schwendler Corp.\eol
175 S Madison\eol
Pasadena, CA 91101}

\Eblock{A. H. Sherman\eol
SCA Inc.\eol
Suite 307\eol
246 Church Street\eol
New Haven, CT 06510}

\Eblock{Kermit Sigmon\eol
Univ. of Florida\eol
Dept. of Math.\eol
Gainesville, FL 32611}

\Eblock{Horst D. Simon\eol
NASA Ames Res. Center\eol
Mail Stop 258-5\eol
Moffett Field, CA 94035}

\Eblock{Richard F. Sincovec\eol
NASA Ames Res. Center\eol
RIACS, Mail Stop 230-5\eol
Moffett Field, CA 94035}

\Eblock{Dennis Smolarski\eol
Univ. of Santa Clara\eol
Dept. of Math.\eol
Santa Clara, CA 95053}

\Eblock{Mitchell D. Smooke\eol
Yale University\eol
Dept. of Mechanical Eng.\eol
P. O. Box 2159, Yale Station\eol
New Haven, CT 06520}

\Eblock{P. Sonneveld\eol
Delft Univ. of Tech.\eol
Dept. of Math. and Informatics\eol
P. O. Box 356\eol
2600 AJ Delft\eol
The Netherlands}

\Eblock{Danny C. Sorensen\eol
Rice Univ.\eol
Dept. of Math. Sciences\eol
Houston, TX 77251-1892}

\Eblock{Alistair Spence\eol
Univ. of Bath\eol
School of Math. Sciences\eol
Bath BA2 7AY\eol
United Kingdom}

\Eblock{Trond Steihaug\eol
Statoil, Forus\eol
P. O. Box 300\eol
N-4001 Stravanger\eol
Norway}

\Eblock{G. W. Stewart\eol
Univ. of Maryland\eol
Dept. of Comp. Sci.\eol
College Park, MD 20742}

\Eblock{G. Strang\eol
Massachusetts Inst. of Tech.\eol
Dept. of Math.\eol
Cambridge, MA 02139}

\Eblock{Jin Su\eol
IBM Corp.\eol
41UC/276\eol
Neighborhood Road\eol
Kingston, NY 12498}

\Eblock{Uwe Suhl\eol
IBM Corp.\eol
Math. Sciences Dept.\eol
Bergen Gruen Strasse 17-19\eol
D-1000 Berlin 38\eol
Germany}

\Eblock{Daniel B. Szyld\eol
Temple Univ.\eol
Dept. of Math.\eol
Philadelphia, PA 19122}

\Eblock{Hillel Tal-Ezer\eol
Brown Univ.\eol
Div. of Applied Math.\eol
Box F\eol
Providence, RI 02912}

\Eblock{Wei Pai Tang\eol
Univ. of Waterloo\eol
Dept. of Comp. Sci.\eol
Waterloo, Ontario N2L 3G1\eol
Canada}

\Eblock{G. D. Taylor\eol
Colorado State Univ.\eol
Math. Dept.\eol
Fort Collins, CO 80523}

\Eblock{Charles Tong\eol
Univ. of California, Los Angeles\eol
Dept. of Comp. Sci.\eol
Los Angeles, CA 90024}

\Eblock{Toru Toyabe\eol
Hitachi Ltd.\eol
Central Res. Lab.\eol
7th Dept.\eol
Kokobunji, Tokyo 185\eol
Japan}

\Eblock{Eugene A. Trabka\eol
Eastman Kodak Company\eol
Res. Laboratories\eol
Rochester, NY 14650}

\Eblock{L. N. Trefethen\eol
Cornell Univ.\eol
Dept. of Comp. Sci.\eol
Ithica, NY 14853}

\Eblock{Donato Trigiante\eol
Univ. di Bari\eol
Inst. di Matematica\eol
Campus Universitario\eol
70125 Bari\eol
Italy}

\Eblock{Donald G. Truhlar\eol
Univ. of Minnesota\eol
Minnesota Supercomputer Inst.\eol
1200 Washington Av. South\eol
Minneapolis, MN 55415}

\Eblock{Kathryn L. Turner\eol
Utah State Univ.\eol
Dept. of Math.\eol
Logan, UT 84322-3900}

\Eblock{Charles Van Loan\eol
Cornell Univ.\eol
Dept. of Comp. Sci.\eol
Ithica, NY 14853}

\Eblock{Henk Van der Vorst\eol
Utrecht Univ.\eol
Math. Inst.\eol
P. O. Box 80.010\eol
3508 TA Utrecht\eol
The Netherlands}

\Eblock{James M. Varah\eol
Univ. of British Columbia\eol
Dept. of Comp. Sci.\eol
Vancouver, British Columbia V6T 1W5\eol
Canada}

\Eblock{Richard S. Varga\eol
Kent State Univ.\eol
Dept. of Math.\eol
Kent, OH 44242}

\Eblock{Anthony Vassiliou\eol
Mobil Res. and Dev. Corp.\eol
Dallas Reseach Lab.\eol
13777 Midway Road\eol
Dallas, TX 75244}

\Eblock{Robert G. Voigt\eol
NASA Langley Res. Center\eol
ICASE, Mail Stop 132-C\eol
Hampton, VA 23665}

\Eblock{Eugene L. Wachspress\eol
Univ. of Tennessee\eol
Dept. of Math.\eol
Knoxville, TN 37996-1300}

\Eblock{Robert C. Ward\eol
Oak Ridge National Lab.\eol
Math. Sciences Section\eol
Building 9207-A\eol
P. O. Box 2009\eol
Oak Ridge, TN 37831}

\Eblock{Daniel D. Warner\eol
Clemson Univ.\eol
Dept. of Math. Sciences\eol
Clemson, SC 29634-1907}

\Eblock{Andrew J. Wathen\eol
Univ. of Bristol\eol
Dept. of Math.\eol
Univ. Walk\eol
Bristol BS8 1TW\eol
United Kingdom}

\Eblock{Bruno Welfert\eol
Univ. of California, San Diego\eol
Dept. of Math.\eol
La Jolla, CA  92093}

\Eblock{Mary J. Wheeler\eol
Univ. of Houston\eol
Dept. of Math.\eol
Houston, TX 77204-3476}

\Eblock{Andrew B. White\eol
Los Alamos National Lab.\eol
Group C-3\eol
Mail Stop B265\eol
Los Alamos, NM 87545}

\Eblock{Jacob White\eol
Massachusetts Inst. of Tech.\eol
Dept. of Electrical Eng.\eol
Cambridge, MA 02139}

\Eblock{Robert. A. Whiteside\eol
Hypercube, Inc.\eol
875 Seminole Dr.\eol
Livermore, CA 94550}

\Eblock{Torbjorn Wiberg\eol
Univ. of Umea\eol
Inst. of Information Processing\eol
S-901 87 Umea\eol
Sweden}

\Eblock{Olaf Widlund\eol
New York Univ.\eol
Courant Institute\eol
251 Mercer Street\eol
New York, NY 10012}

\Eblock{Harry A. Wijshoff\eol
Univ. of Illinois\eol
305 Talbot Lab.\eol
104 S Wright Street\eol
Urbana, IL 61801}

\Eblock{Roy. S. Wikramaratna\eol
Winfrith Petroleum Tech.\eol
Dorchester\eol
Dorset DT2 8DH\eol
United Kingdom}

\Eblock{David S. Wise\eol
Indiana Univ.\eol
Dept. of Comp. Sci.\eol
Bloomington, IN 47405-6171}

\Eblock{P. H. Worley\eol
Oak Ridge National Lab.\eol
Math. Sciences Section\eol
Building 9207-A\eol
P. O. Box 2009\eol
Oak Ridge, TN 37831}

\Eblock{Margaret H. Wright\eol
AT\&T Bell Laboratories\eol
600 Mountain Av.\eol
Murray Hill, NJ 07974-2070}

\Eblock{Steven J. Wright\eol
Argonne National Lab.\eol
Math. and Comp. Sci. Div.\eol
9700 S Cass Av.\eol
Argonne, IL 60439}

\Eblock{Kuo W. Wu\eol
Cray Res. Inc\eol
1333 Northland Drive\eol
Mendota Heights, MN 55120-1095}

\Eblock{Chao W. Yang\eol
Cray Res. Inc\eol
Math. Software Group\eol
1408 Northland Drive\eol
Mendota Heights, MN 55120-1095}

\Eblock{Gung-Chung Yang\eol
Univ. of Illinois\eol
305 Talbot Lab.\eol
104 S Wright Street\eol
Urbana, IL 61801}

\Eblock{Elizabeth Yip\eol
Boeing Aerospace Corp.\eol
Mail Stop 8K-17\eol
P. O. Box 3999\eol
Seattle, WA 98124-2499}

\Eblock{David P. Young\eol
Boeing Comp. Services\eol
Mail Stop 7L-21\eol
P. O. Box 24346\eol
Seattle, WA 98124-0346}

\Eblock{David M. Young\eol
Univ. of Texas\eol
Center for Numerical Analysis\eol
Moore Hall 13.150\eol
Austin, TX 78712}

\Eblock{Earl Zmijewski\eol
Univ. of California, Santa Barbara\eol
Comp. Sci. Dept.\eol
Santa Barbara, CA 93106}

\Eblock{Qisu Zou\eol
Kansas State Univ.\eol
Dept. of Math.\eol
Manhattan, KS 66506}

}

 \Lblock {S. F. Ashby, LLNL, L-316\inc 1\eol
 J. B. Bell, LLNL, L-316\inc 1\eol
 J. H. Bolstad, LLNL, L-16\inc 1\eol
 P. N. Brown, LLNL, L-316\inc 1\eol
 M. R. Dorr, LLNL, L-316\inc 1\eol
 G. W. Hedstrom, LLNL, L-321\inc 1\eol
 A. C. Hindmarsh, LLNL, L-316\inc 1\eol
 D. S. Kershaw, LLNL, L-471\inc 1\eol
 A. E. Koniges, LLNL, L-561\inc 1\eol
 L. R. Petzold, LLNL, L-316\inc 1\eol
 G. H. Rodrigue, LLNL, L-306\inc 1\eol
 J. A. Trangenstein, LLNL, L-316\inc 1}

\vskip \parskip 
 \settabs \+\hglue 0.45in& \hglue 0.2in& \hglue 0.25in\cr
 \+1420&W. J. Camp\inc 1\cr
 \+1421&S. S. Dosanjh\inc 1\cr
 \+1422&R. C. Allen, Jr.\inc 1\cr
 \+1422&R. S. Tuminaro\inc 1\cr
 \+1422&D. E. Womble\inc 1\cr
 \+1423&E. F. Brickell\inc 1\cr
 \+1424&R. E. Brenner\inc 1\cr
 \+1553&W. L. Hermina\inc 1\cr
 \+2113&H. A. Watts\inc 1\cr
 \+8000& J. C. Crawford\inc 1\cr
 \+&&Attn: &E. E. Ives, 8100\cr
 \+&&&P. L. Mattern, 8300\cr
 \+&&&R. C. Wayne, 8400\cr
 \+&&&P. E. Brewer, 8500\cr
 \+8200& R. J. Detry\inc 1\cr
 \+&&Attn: &C. W. Robinson, 8240\cr
 \+&&&R. C. Dougherty, 8270\cr
 \+&&&R. A. Baroody, 8280\cr
 \+8210& J. C. Meza\inc 1\cr
 \+8210& W. D. Wilson\inc 1\cr
 \+&&Attn: &R. E. Cline, 8210\cr
 \+&&&J. M. Harris, 8210\cr
 \+&&&R. Y. Lee, 8210\cr
 \+8230& D. L. Crawford\inc 1\cr
 \+&&Attn: &R. E. Palmer, 8234\cr
 \+&&&P. W. Dean, 8236\cr
 \+&&&B. Stiefeld, 8237\cr
 \+8233&J. F. Grcar (40)\inc 40\cr
 \+8233&R. S. Judson\inc 1\cr
 \+8233&J. F. Lathrop\inc 1\cr
 \+8233&W. E. Mason\inc 1\cr
 \+8233&D. S. McGarrah\inc 1\cr
 \+8236&T. H. Jefferson\inc 1\cr
 \+8241&K. J. Perano\inc 1\cr
 \+8243&L. A. Bertram\inc 1\cr
 \+8245&D. S. Dandy\inc 1\cr
 \+8245&G. H. Evans\inc 1\cr
 \+8245&W. G. Houf\inc 1\cr
 \+8245&R. J. Kee\inc 1\cr
 \+8245&A. E. Lutz\inc 1\cr
 \+8245&W. S. Winters\inc 1\cr
 \+8364&J. K. Bechtold\inc 1\cr
 \+8364&S. B. Margolis\inc 1\cr
 \vskip \parskip
 \+8535&Publications/Tech. Lib. Processes, 3141\inc 1\cr
 \+3141&Tech. Library Process Division (2)\inc 2\cr
 \+8524-2& \ Central Tech. Files (3)\inc 3\cr

\relax \vfill \eject

}


 \end